\documentclass[12pt,a4paper]{article}
\usepackage{amsfonts,amsthm,amsmath,latexsym,bm,graphics,indentfirst}
\usepackage[hypertex]{hyperref}
\newtheorem{theorem}{Theorem}[section]
\newtheorem{lemma}[theorem]{Lemma}
\newtheorem{proposition}[theorem]{Proposition}
\newtheorem{corollary}[theorem]{Corollary}
\newtheorem{remark}[theorem]{Remark}

\newcommand{\D}{\displaystyle}
\newcommand{\F}[2]{\frac{#1}{#2}}
\newcommand{\DIV}{\operatorname{div}}
\newcommand{\BE}{\begin{equation}\begin{aligned}}
\newcommand{\BEN}{\begin{equation*}\begin{aligned}}
\newcommand{\EE}{\end{aligned}\end{equation}}
\newcommand{\EEN}{\end{aligned}\end{equation*}}
\newcommand{\BL}{\begin{lemma}}
\newcommand{\EL}{\end{lemma}}
\newcommand{\BT}{\begin{theorem}}
\newcommand{\ET}{\end{theorem}}
\newcommand{\BP}{\begin{proposition}}
\newcommand{\EP}{\end{proposition}}
\newcommand{\BC}{\begin{corollary}}
\newcommand{\EC}{\end{corollary}}
\allowdisplaybreaks

\title{Non-compactness of the Prescribed  $Q-\text{curvature}$ Problem in Large Dimensions}

\author{Juncheng Wei\thanks{Department of Mathematics, The Chinese University of Hong Kong, Shatin, Hong Kong. Email: wei@math.cuhk.edu.hk} \qquad
        Chunyi Zhao\thanks{Department of Mathematics, East China Normal University, Shanghai, 200241, China. Email: cyzhao@math.ecnu.edu.cn}}
        
\date{}

\begin{document}
\maketitle

\begin{abstract}
Let $(M, g)$ be a compact Riemannian manifold of dimension $N \geq 5$ and $Q_g$ be its $Q$ curvature.  The prescribed $Q$ curvature problem  is concerned with finding metric of constant $Q$ curvature in the conformal class of $g$. This amounts to finding a positive solution to
\[ P_g (u)= c u^{\frac{N+4}{N-4}},\quad u>0 \qquad \mbox{on} \ M\]
where $P_g$ is the Paneitz operator. We show that for dimensions $N \geq 25$, the set of all  positive solutions  to the prescribed $Q$ curvature problem is {\em non-compact}.
\end{abstract}

\section{Introduction}

Let $(M, g)$ be a Riemannian manifold of dimension $N$. A basic question in conformal geometry is
the following: can one change the original metric $g$ conformally into a new metric $g'$ with prescribed
properties? This means that one searches for some positive function $\psi$ (conformal factor)  such that $g' = \psi g$  and the new metric $g'$ has prescribed properties.

A best known example is the so-called Yamabe problem. For $N \geq 3$, let $L_g := -\F{4(N-1)}{N-2} \Delta_g + S_g$ be the conformal Laplacian, where $\Delta_g$ is the Laplace-Beltrami
operator and $S_g$ is the scalar curvature. If one sets the conformal
factor $ \psi = u^\F{4}{N-2}$ ($u > 0$), then it is well known that $L_g$ has the following conformal covariance property:
\begin{align*}
L_g(u\varphi)=u^\F{N+2}{N-2}L_{g'}(\varphi) \qquad \forall~\varphi\in C^\infty(M).
\end{align*}
If one prescribes the scalar curvature $S_{g'}$ for the metric $g'$ then $u$ has to satisfy the second-order
equation
\begin{equation}\label{1.2}
L_g(u)=u^\F{N+2}{N-2} L_{g'}(1)= S_{g'}u^\F{N+2}{N-2}.
\end{equation}
In the case when $S_{g'}$ is a constant, this is the  Yamabe problem. In the case when $S_{g'}$ is a prescribed function,
it is called the Nirenberg problem.

Precisely, the Yamabe equation is
\begin{align}\label{1.3}
\F{4(N-1)}{N-2} \Delta_g u - S_g u + cu^\F{N+2}{N-2}=0.
\end{align}
The question that whether the set of all solutions to
the Yamabe problem (\ref{1.3}) is compact in the $C^\infty-\text{topology}$ has been widely studied. It has been conjectured that this should be true unless $(M, g)$ is conformally equivalent to the round
sphere (see \cite{S1,S2,S3}).  The case of the round sphere ${\mathbb S}^N$ is special in that
(\ref{1.3}) is invariant under the action of the conformal group on ${\mathbb S}^N$, which is
non-compact.  The Compactness Conjecture has been verified in low dimensions and locally conformally flat by R. Schoen \cite{S2,S3}.
He also proposed a strategy to proving the conjecture in the non-locally conformally flat case.
Developing further this strategy, the conjecture is proved in low dimensions:
$N=3$ by Li-Zhu \cite{LZhu}, $N=4,5$ by Druet \cite{D}, $N=6,7$ by
Li-Zhang \cite{LZhang1} and Marques \cite{M},  $N=10, 11$ by Li-Zhang \cite{LZhang2} under the Positive Mass Theorem assumption.
Recently this conjecture is shown to be true  by Khuri-Marques-Schoen \cite{KMS} for dimensions $N \leq 24$ under the Positive Mass Theorem condition. On the other hand, the Compactness Conjecture is not true  for $N\geq 25$ in the recent papers  by Brendle \cite{B} ($N \geq 52$) and Brendle-Marques \cite{BM} ($25 \leq N \leq 51$).
More precisely, given any integer $N\geq 25$, there exists a
smooth Riemannian metric $g$ on ${\mathbb S}^N$ such that set of constant scalar curvature
metrics in the conformal class of $g$ is non-compact. Moreover, the blowing-up
sequences obtained in \cite{B,BM} form exactly one bubble. The construction
relies on a gluing procedure based on some local metric. The non-compactness of Yamabe problem in the $C^k-\text{topology}$ is studied by Ambrosetti-Malchiodi \cite{AM} and Berti-Malchiodi \cite{BeM}. A complete description of blow-up behavior of Yamabe type problems can be found in the book Druet-Hebey-Robert \cite{DHR}.

Besides  the conformal Laplacian $L_g$, there are many other  operators  which enjoy a conformal covariance property.
A particularly interesting one is the fourth order operator $P_g$ which was discovered by Paneitz in
1983 (\cite{P}), which can be written for $N \geq 5$ as follows:
\begin{align}
P_g=\Delta_g^2-\DIV_g(a_N S_g g+b_N {\mathcal R}_g)d + \F{N-4}{2}Q_g,
\end{align}
where $\DIV_g$ is the divergence of a vector-field, $d$ is the differential operator,
\begin{gather}
a_N=\F{(N-2)^2+4}{2(N-1)(N-2)}, \qquad b_N=-\F{4}{N-2},\nonumber\\
\intertext{and}
Q_g=-\F{1}{2(N-1)}\Delta_g S_g +\F{N^3-4N^2+16N-16}{8(N-1)^2(N-2)^2}S_g^2-\F{2}{(N-2)^2}|\mathcal R_g|^2. \label{1.1}
\end{gather}
Here $\mathcal R_g$ the Ricci tensor, $|\mathcal R_g|^2=\sum_{i,j}\mathcal R^{ij}\mathcal R_{ij}$ where $\mathcal R^{ij}=\sum_{s,t}g^{is}\mathcal R_{st}g^{tj}$.
$Q_g$ is the so called $Q-\text{curvature}$.
In the case $N > 4$, the conformal factor is usually chosen in the form $\psi = u^\F{4}{N-4}$ $(u > 0)$
and the conformal covariance property of the Paneitz operator reads as follows:
\begin{align*}
P_g(u\varphi)=u^\F{N+4}{N-4}P_{g'}(\varphi) \qquad \forall~\varphi\in C^\infty(M).
\end{align*}
If one prescribes the $Q-\text{curvature}$ for the metric $g'$ by a function $Q_{g'}$ this leads to the equation
\begin{align*}
P_g(u) = u^\F{N+4}{N-4}P_{g'}(1) = \F{N-4}{2}Q_{g'}u^\F{N+4}{N-4},
\end{align*}
which is a fourth-order analogue of (\ref{1.2}). We refer to the survey articles of Chang \cite{C1} and Chang-Yang \cite{CY1} and  the lecture notes \cite{C2,CEOY} for more background information on the Paneitz operator. Recently, there are more and more interests in using higher order partial differential
equations in the study of conformal geometry. See \cite{CY2,DM,WX} and references therein.\par

There is an analogue problem to the Yamabe problem, that is, to find metrics of constant $Q-\text{curvature}$ in the
conformal class of $g$. The problem can be transformed to solving the following $Q-\textit{curvature}$ \textit{equation}, for $N\geq 5$,
\begin{align}\label{1.4}
P_g u=\F{N-4}{2} u^\F{N+4}{N-4}, \ \ \ u>0 \qquad \text{on }M.
\end{align}
Clearly, every solution of (\ref{1.4}) is a critical point of the functional
\begin{align}
\mathcal E_g(u) =\F {\int_M (\Delta_g u)^2 + \sum_{i,j}(a_N S_g g^{ij}+b_N {\mathcal R}^{ij})\partial_i u \partial_j u +\F{N-4}{2}Q_g u^2 \mathrm{d}vol_g}
{\left(\int_M u^\F{2N}{N-4} \mathrm{d}vol_g\right)^\F{N-4}{N}}.
\end{align}
Consider
\begin{align*}
P(g) = \inf \Big\{\mathcal E_g(u):\ u\in H^2(M),\ u>0 \Big\}.
\end{align*}
We refer to $P(g)$ as the Paneitz energy. Clearly it is a conformal invariant.

A similar question is whether or not  the set of all positive solutions to the $Q-\text{curvature}$ equation is compact.
As far as we know, the only results in this direction are given by Hebey-Robert \cite{HR} and Qing-Raske \cite{QR}. 
Both papers  give   positive answers provided $M$ is locally conformally flat and satisfies some additional  assumptions. 
Compactness is also studied for non-geometric potentials of Paneitz operator in Hebey-Robert-Wen \cite{HRW}.

In this paper, we prove  the non-compactness of the set of solutions to the $Q-\text{curvature}$ problem in large dimensions. We construct a blowing-up sequence consisting of exactly one bubble. More precisely we prove the following theorem.

\BT
\label{t1.1}
Assume $N\geq 25$. Then there exists a $C^\infty$ Riemannian metric $g$ on ${\mathbb S}^N$ and a sequence of positive function $u_n$ ($n\in\mathbb N$)
with the following properties:
\begin{itemize}
\item[i)] $g$ is not conformally flat,
\item[ii)] $u_n$ is a positive solution of the $Q-\text{curvature}$ equation (\ref{1.4}) for all $n$,
\item[iii)] ${\mathcal E}_g(u_n)< P({\mathbb S}^N)$ for all $n$ and $\mathcal E_g (u_n)\to P({\mathbb S}^N)$ as $n\to \infty$,
\item[iv)] $\sup_{{\mathbb S}^N} u_n \to \infty$ as $n\to\infty$.
\end{itemize}
Here $P({\mathbb S}^N)$ denotes the Paneitz energy of the round metric on ${\mathbb S}^N$.
\ET

For convenience (and by stereo-graphic projection), we will work on $\mathbb R^N$ instead of  ${\mathbb S}^N$. Let $g$ be a smooth metric on $\mathbb R^N$
which agrees with the Euclidean metric outside a ball of radius $1$. We will assume throughout the paper that $\det g(x) = 1$ for all $x\in \mathbb R^N$, so that the volume form
associated with $g$ agrees with the Euclidean volume form. Precisely, we will consider
\begin{align}\label{main equation}
P_g u=\F{N-4}{2} u^\F{N+4}{N-4},\quad u>0 \qquad \text{in }\mathbb R^N.
\end{align}

Our goal is to construct solutions to the $Q-\text{curvature}$ equation  (\ref{main equation}) on $(\mathbb R^N,g)$.
Though we shall follow the main ideas in \cite{B} and \cite{BM}, there are some major difficulties in fourth order equations.
The main difficulty is that we need to ensure that $ u$ is strictly positive on $\mathbb R^N$. In the Yamabe problem case, one constructs solutions of the form
\[ u = u_0 + \phi\]
where $u_0$ is the standard bubble and $\phi$ is the error.
As long as $\| \phi\|_{H^1}$ is small, it can be shown easily that $u>0$.
In the prescribed $Q$ curvature problem, even if we can show that $ \| \phi \|_{H^2 (\mathbb R^N)}$ is small, it is not guaranteed that $u$ is positive. To overcome this difficulty, we need to use a weighted $L^\infty$ norm
\[ \| \phi\|_{*}= \| u_0^{-1} \phi \|_{L^\infty (\mathbb R^N)} \]
and we have to show that $\| \phi \|_{*}$ is small which then implies that $u$ is positive.
We use a  technical framework  which is more closely related to the finite dimensional Liapunov-Schmidt reduction procedure, as in  \cite{DFM}, \cite{RW}, and \cite{WY}.
Another difficulty is the choice of the auxiliary function $f(s)$.  
In the second order Yamabe problem case, a linear function is chosen to obtain $N \geq 52$ (\cite{B}) and a  cubic polynomial is chosen to obtain $N \geq 25$ (\cite{BM}). 
While  in the fourth order case, the best choice for $f$ seems to be fourth order polynomial. (Linear function only gives $ N \geq 52 $.) A surprising fact is that in both Yamabe and $Q-$curvature problems, the two dimensions $52$ and $25$ are the same.

The organization of the paper is as follows: In Section 2, we introduce the special metric $g$ in this paper. In Section 3, the first approximation of the solutions is given.
In Section 4, we calculate the corresponding $\Delta_g$, $S_g$, $R_g$, $Q_g$ under this special metric $g$, and then acquire
the estimate of the energy functional. In Section 5, the invertibility of the linearized operator is settled. In Section 6,
we solve the nonlinear problem. In Section 7, a variational reduction procedure is processed. In Section 8, we show that
the energy can be approximated by an  reduced energy functional. In Section 9, we compute the reduced energy functional in terms of an   auxiliary function $f$. In Section 10, we choose a linear auxilliary function to show that the reduced energy functional has a strict local minimum when $N \geq 52$. In Section 11, the dimension is reduced to $N \geq 25$ by choosing a fourth order polynomial. In Section 12, we prove the main theorem by a gluing method. Finally in the appendix, some inequalities used in the paper is proven.

\medskip

\noindent
{\bf Acknowledgment.} The first
author is supported by an Earmarked
Grant  from RGC  of Hong Kong and a Focused Research Scheme from CUHK.

\section{The Special Metric $\bm g$}

In this section we introduce the metric which will be used in this paper.\par

In what follows, we fix a multi-linear form $W$: $\mathbb R^N\times\mathbb R^N\times\mathbb R^N\times\mathbb R^N\to \mathbb R$.
We assume that $W_{ikj\ell}$ satisfy all the algebraic properties of the Weyl tensor.
Moreover, we assume that some components of $W$ are non-zero, so that
\begin{align*}
\sum_{i,j,k,\ell=1}^n (W_{ikj\ell}+W_{i\ell jk})^2>0.
\end{align*}
For simplicity, we put
\begin{align*}
H_{ij}(y)=\sum_{p,q=1}^n W_{ipjq}y_py_q
\intertext{and}
\overline H_{ij}(y)=f(|y|^2)H_{ij}(y),
\end{align*}
where $f(|y|^2) $ will be chosen later. (Specifically, we shall choose $ f(s)=\tau-12000s+2411s^2-135s^3+s^4$. The number $\tau$ depends only on the dimension $N$ and will be chosen later.)
It is easy to see that $H_{ij}(y)$ is trace-free, $\sum_{i=1}^n y_i H_{ij}(y)=0$ and $\sum_{i=1}^n \partial_i H_{ij}(y)=0$
for all $y\in\mathbb R^N$.\par

We consider a Riemannian metric of the form $g(x)=e^{h(x)}$, where $h(x)$ is a trace-free symmetric two-tensor on $\mathbb R^N$
satisfying $h(0)=0$, $h(x)=0$ for $|x|\geq 1$,
\begin{align*}
|h(x)|+|\partial h(x)|+|\partial^2 h(x)|+|\partial^3 h(x)|+|\partial^4 h(x)|\leq \alpha
\end{align*}
for all $x\in\mathbb R^N$, where $\alpha>0$ is a fixed small number, and
\begin{align*}
h_{ij}(x)=\mu\varepsilon^8 f(\varepsilon^{-2}|x|^2)H_{ij}(x)
\end{align*}
for $|x|\leq \rho$. We assume that the parameter $\varepsilon$, $\mu$ and $\rho$ are chosen such that $\mu\leq 1$ and $\varepsilon\leq\rho\leq 1$.
Note that $\sum_{i=1}^N x_i h_{ij}(x)=0$ and $\sum_{i=1}^N \partial_i h_{ij}(x)=0$ for $|x|\leq \rho$.\par

For later purpose, we need to understand the Green's function of $\Delta_g^2$.
Denote $\widetilde G(x,y)$ be the Green's function of $\Delta_g$ in $\mathbb R^N$. Then the Green's function of $\Delta^2_g$ is
\begin{align*}
G(x,y)=\int_{\mathbb R^N} \widetilde G(x,z) \widetilde G(y,z) \mathrm{d}z.
\end{align*}
Since $|\widetilde G(x,y)| \leq C[d(x,y)]^{2-N}$, $|\nabla \widetilde G(x,y)| \leq C[d(x,y)]^{1-N}$ and $|\nabla^2 \widetilde G(x,y)| \leq C[d(x,y)]^{-N}$,
we know that
\begin{align*}
|G(x,y)| &\leq C[d(x,y)]^{4-N},\\
|\nabla  G(x,y)| &\leq C[d(x,y)]^{3-N},\\
|\nabla^2  G(x,y)| &\leq C[d(x,y)]^{2-N}.
\end{align*}
Here $d(x,y)$ is the distance between $x$ and $y$ under the metric $g$. It is easy to see that
\begin{align*}
d(x,y)=[1+O(\alpha)]|x-y|
\end{align*}
since $g=e^h$ and $\alpha$ is sufficiently small.\par

\medskip

\noindent
\textit{Notations} In what follows, we use $C$ to denote the variable constant which is independent of $\alpha$ and $\varepsilon$.
$|O(A)|\leq CA$ and $o(A)/A\to 0$ as $\varepsilon\to 0$.

\section{First Approximation of the Solutions}

In this section we will provide an ansatz for solutions of Problem (\ref{main equation}).\par

Denote
\begin{align*}
u_0(x)=\gamma_N \left(\F{\lambda' \varepsilon}{\lambda'^2 \varepsilon^2+|x-\xi|^2}\right)^\F{N-4}{2},
\end{align*}
where $\gamma_N=\big[N(N-4)^2(N-2)(N+2)/2\big]^{-\F{N-4}{8}}$, $\lambda'>0$ and $\varepsilon>0$.
It is well known \cite{L} that $u_0$ is the only positive solution of
\begin{align*}
\Delta^2 u_0 = \F{N-4}{2} u_0^\F{N+4}{N-4} \qquad \text{in }\mathbb R^N.
\end{align*}

Observe that $u(x)$ satisfies (\ref{main equation}) if and only if $v(y)= \varepsilon^\F{N-4}{2} u(\varepsilon y)$
satisfies
\begin{align}\label{b}
P_{\tilde g} v =\F{N-4}{2} v^\F{N+4}{N-4}
\end{align}
where $\tilde g(y)=g(\varepsilon y)$. Denote
\begin{align*}
\tilde u_0(y) &=\varepsilon^\F{N-4}{2}u_0(\varepsilon y)=\gamma_N  \left(\F{\lambda' }{\lambda'^2 +|y-\xi'|^2}\right)^\F{N-4}{2},
\end{align*}
where $\xi'=\xi/\varepsilon$, $\lambda'=\lambda/\varepsilon$.
The configuration set for $(\xi',\lambda')$ is
\BEN
\Lambda = \left\{(\xi',\lambda')\in\mathbb R^N\times\mathbb R:\ |\xi'|\leq 1,\ \F{1}{2}<\lambda'<\F{3}{2}\right\}.
\EEN

We will look for a solution to (\ref{main equation}) with the form $\tilde u_0(y)+\phi(y)$.
It is easy to check that $\phi$ must be a solution of
\begin{align}\label{phi}
P_{\tilde g}\phi-\F{N+4}{2}\tilde u_0^\F{8}{N-4}\phi = -R+N(\phi) \qquad  \text{in }\mathbb R^N,
\end{align}
where
\begin{align}
R(y) &=P_{\tilde g}\tilde u_0- \F{N-4}{2} \tilde u_0^\F{N+4}{N-4}, \label{2.1}\\
N(\phi) &= \F{N-4}{2}(\tilde u_0+\phi)^\F{N+4}{N-4} -\F{N-4}{2} \tilde u_0^\F{N+4}{N-4} -\F{N+4}{2} \tilde u_0^\F{8}{N-4}\phi. \label{2.2}
\end{align}

A main step in solving (\ref{phi}) for small $\phi$ is that of a solvability theory for the linearized operator $P_{\tilde g}-\F{N+4}{2}\tilde u_0^\F{8}{N-4}$.

\section{Preliminary Estimates}

In this section we will mainly estimate the energy functional. From now on, denote $\tilde h(y)=h(\varepsilon y)$ and $\widetilde{\mathcal R}$,
$S_{\tilde g}$, $Q_{\tilde g}$ are the corresponding Ricci tensor, scaler curvature, $Q-\text{curvature}$ under the metric $\tilde g$.

\BL\label{l7.1}
For $|y|\leq \rho/\varepsilon$,  it holds
\begin{align*}
\widetilde{\mathcal R}_{ij}=&\ -\sum_m \F{1}{2}\partial_{mm}\tilde h_{ij}+\sum_{m,s}\F{1}{2}[\tilde h_{ms}(\partial_{ms}\tilde h_{ij})-(\partial_s \tilde h_{mj})(\partial_m \tilde h_{si})]\\
&\ +\sum_{m,s}\F{1}{4}\bigg[(\partial_i \tilde h_{ms})(\partial_m \tilde h_{sj})+(\partial_j \tilde h_{ms})(\partial_m \tilde h_{si})-(\partial_j \tilde h_{ms})(\partial_i \tilde h_{sm})\\
&\qquad\quad -\tilde h_{ms}(\partial_{mi}\tilde h_{sj})-\tilde h_{ms}(\partial_{mj}\tilde h_{si})-(\partial_{mm}\tilde h_{js})\tilde h_{si}-\tilde h_{js}(\partial_{mm}\tilde h_{si})\bigg]\\
&\ +O(|\tilde h|^2|\partial^2 \tilde h|+|\tilde h||\partial \tilde h|^2).
\end{align*}
For $|y|\geq \rho/\varepsilon$, we have
\begin{align*}
\widetilde{\mathcal R}_{ij}= O(\alpha\varepsilon^2).
\end{align*}
\EL

\begin{proof}
Recall that $\widetilde\Gamma_{k\ell}^i=\sum_{m}\F{1}{2}\tilde g^{im}(\partial_\ell \tilde g_{mk}+\partial_k \tilde g_{m\ell}-\partial_m \tilde g_{k\ell})$.
Since $\tilde h$ is trace-free, we have $\det\tilde g=1$ for all $y\in\mathbb R^N$. This implies $\sum_{i}\widetilde\Gamma_{ik}^i=\sum_{i,\ell}\F{1}{2}\tilde g^{i\ell}\partial_k \tilde g_{i\ell}=\F{1}{2}\partial_k\log |\tilde g|=0$.
Therefore, we obtain
\begin{align}
\widetilde{\mathcal R}_{ij} &= \sum_m \partial_m \widetilde\Gamma_{ji}^m - \sum_m \partial_j \widetilde\Gamma_{mi}^m +\sum_{\ell,m}\widetilde\Gamma_{m\ell}^m \widetilde\Gamma_{ji}^\ell - \sum_{\ell,m}\widetilde\Gamma_{j\ell}^m \widetilde\Gamma_{mi}^\ell\nonumber\\
&= \sum_m \partial_m \widetilde\Gamma_{ji}^m-\sum_{\ell,m}\widetilde\Gamma_{j\ell}^m \widetilde\Gamma_{mi}^\ell.\label{Ricci}
\end{align}
Direct calculation shows
\begin{align*}
\sum_m \partial_m \widetilde\Gamma_{ji}^m =
&\ \sum_m\F{1}{2}\Big[\partial_{mi}\tilde h_{mj}+\partial_{mj}\tilde h_{mi}-\partial_{mm}\tilde h_{ij}\Big]\\
&\ +\sum_{m,s} \F{1}{2}\Big[(\partial_m \tilde h_{ms})(\partial_s \tilde h_{ij})+\tilde h_{ms}(\partial_{ms}\tilde h_{ij})\Big]\\
&\ +\sum_{m,s} \F{1}{4}\Big[(\partial_i \tilde h_{ms})(\partial_m \tilde h_{sj})-(\partial_m \tilde h_{ms})(\partial_i \tilde h_{sj})+(\partial_j \tilde h_{ms})(\partial_m \tilde h_{si})\\
&\qquad\quad -(\partial_m \tilde h_{ms})(\partial_j \tilde h_{si})-2(\partial_m \tilde h_{js})(\partial_m \tilde h_{si})+(\partial_{mi}\tilde h_{ms})\tilde h_{sj}\\
&\quad\qquad -(\partial_{mi}\tilde h_{sj})\tilde h_{ms}+(\partial_{mj}\tilde h_{ms})\tilde h_{si}-(\partial_{mj}\tilde h_{si})\tilde h_{ms}-(\partial_{mm}\tilde h_{js})\tilde h_{si}\\
&\quad\qquad -(\partial_{mm}\tilde h_{si})\tilde h_{js}\Big]+O(|\tilde h|^2|\partial^2 \tilde h|+|\tilde h||\partial \tilde h|^2),
\end{align*}
and
\begin{align*}
\sum_{\ell,m}\widetilde\Gamma_{j\ell}^m \widetilde\Gamma_{mi}^\ell =&\ \sum_{\ell,m}\F{1}{2}\Big[(\partial_\ell \tilde h_{mj})(\partial_m \tilde h_{\ell i})-(\partial_\ell \tilde h_{mj})(\partial_\ell \tilde h_{mi})\Big]
+\F{1}{4}(\partial_j \tilde h_{m\ell})(\partial_i \tilde h_{\ell m}) \\
&\ +O(|\tilde h||\partial \tilde h|^2).
\end{align*}
Since $\sum_m \partial_m \tilde h_{mk}=0$ for $|y|\leq \rho/\varepsilon$, the lemma follows from (\ref{Ricci}).
\end{proof}

Thus we have the following calculations  for the Ricci tensor $\widetilde{\mathcal R}^{ij}=\sum_{s,t}\tilde g^{is}\widetilde{\mathcal R}_{st}\tilde g^{t j}$.

\BC\label{c7.1}
For $|y|\leq \rho/\varepsilon$,  we have
\begin{align*}
\widetilde{\mathcal R}^{ij}=&\ -\sum_m \F{1}{2}\partial_{mm}\tilde h_{ij}+\sum_{m,s}\F{1}{2}[\tilde h_{ms}(\partial_{ms}\tilde h_{ij})-(\partial_s\tilde h_{mj})(\partial_m\tilde h_{si})]\\
&\ +\sum_{m,s}\F{1}{4}\bigg[(\partial_i\tilde h_{ms})(\partial_m\tilde h_{sj})+(\partial_j\tilde h_{ms})(\partial_m\tilde h_{si})-(\partial_j\tilde h_{ms})(\partial_i\tilde h_{sm})\\
&\qquad\quad -\tilde h_{ms}(\partial_{mi}\tilde h_{sj})-\tilde h_{ms}(\partial_{mj}\tilde h_{si})+(\partial_{mm}\tilde h_{js})\tilde h_{si}+h_{js}(\partial_{mm}\tilde h_{si})\bigg]\\
&\ +O(|\tilde h|^2|\partial^2\tilde h|+|\tilde h||\partial\tilde h|^2).
\end{align*}
For $|y|\geq \rho/\varepsilon$, it holds
\begin{align*}
\widetilde{\mathcal R}^{ij}= O(\alpha\varepsilon^2).
\end{align*}
\EC

\BL\label{l7.5}
For $|y|\leq \rho/\varepsilon$, there holds
\begin{align*}
S_{\tilde g}=-\F{1}{4}\sum_{k,\ell,m}(\partial_\ell \tilde h_{mk})^2+O(|\tilde h|^2|\partial^2 \tilde h|+|\tilde h||\partial\tilde h|^2).
\end{align*}
For $|y|\geq \rho/\varepsilon$,  we have
\begin{align*}
S_{\tilde g}=O(\alpha\varepsilon^2).
\end{align*}
\EL

\begin{proof}
The  detailed proof can be found  in \cite[Prop. 26]{B}, noting  that $\sum_m \partial_m \tilde h_{mk}=0$ in $|y|\leq \rho/\varepsilon$.
\end{proof}

The above lemma and a direct computation show the following conclusion.
\BC\label{c7.2}
For $|y|\leq \rho/\varepsilon$, we have
\begin{align*}
\Delta_{\tilde g} S_{\tilde g}=&\ -\F{1}{2}\sum_{i,k,\ell,m}(\partial_{i\ell}\tilde h_{mk})^2-\F{1}{2}\sum_{i,k,\ell,m}(\partial_\ell\tilde h_{mk})(\partial_{ii\ell}\tilde h_{mk}) \\
&\ +O(|\partial\tilde h|^2|\partial^2\tilde h|+|\tilde h||\partial^2\tilde h|^2+|\tilde h||\partial\tilde h||\partial^3\tilde h|+|\tilde h|^2|\partial^4\tilde h|).
\end{align*}
For $|y|\geq \rho/\varepsilon$, it holds
\begin{align*}
\Delta_{\tilde g} S_{\tilde g} =O(\alpha\varepsilon^4).
\end{align*}
\EC

Now it is ready to estimate $Q_{\tilde g}$.
\BL\label{l7.4}
For $|y|\leq\rho/\varepsilon$, we have
\begin{align*}
Q_{\tilde g}=&\ \F{1}{4(N-1)}\sum_{i,k,\ell,m}\Big[(\partial_{i\ell}\tilde h_{mk})^2+(\partial_\ell\tilde h_{mk})(\partial_{ii\ell}\tilde h_{mk})\Big]\\
&\ -\F{1}{2(N-2)^2}\sum_{i,j,m,s}(\partial_{mm}\tilde h_{ij})(\partial_{ss}\tilde h_{ij})\\
&\ +O(|\partial\tilde h|^2|\partial^2\tilde h|+|\tilde h||\partial^2\tilde h|^2+|\tilde h||\partial\tilde h||\partial^3\tilde h|+|\tilde h|^2|\partial^4\tilde h|).
\end{align*}
For $|y|\geq\rho/\varepsilon$, there holds
\begin{align*}
Q_{\tilde g}=O(\alpha\varepsilon^4).
\end{align*}
\EL

\begin{proof}
This is a direct result of (\ref{1.1}), Lemma \ref{l7.1}, Corollary \ref{c7.1} and \ref{c7.2}.
\end{proof}

Our next goal is  to estimate $R(y)$ defined in (\ref{2.1}).

\BL
For $|y|\leq \F{\rho}{\varepsilon}$,
\begin{align*}
\Delta_{\tilde g}^2\tilde u_0-\Delta^2 \tilde u_0=&\ -(\partial_{ss}\tilde h_{ij})(\partial_{ij}\tilde u_0)-2(\partial_s\tilde h_{ij})(\partial_{sij}\tilde u_0)-2\tilde h_{ij}(\partial_{ssij}\tilde u_0)\\
&\ +O(|\tilde h||\partial^2\tilde h|)|\partial^2\tilde u_0|+O(|\tilde h||\partial\tilde h|)|\partial^3\tilde u_0|+O(|\tilde h|^2)|\partial^4\tilde u_0|.
\end{align*}
For $|y|\geq \F{\rho}{\varepsilon}$,
\begin{align*}
\Delta_{\tilde g}^2\tilde u_0-\Delta^2 \tilde u_0= O\left(\F{\alpha\varepsilon}{(1+|y-\xi'|)^{N-1}}\right).
\end{align*}
\EL

\begin{proof}
The computations follow easily  from  the definition of $\Delta_{\tilde g}$ and the properties of $h$.
\end{proof}

By Lemma \ref{l7.5}, Corollary \ref{c7.1} and the properties of $h$, it is also not difficult to verify the following two lemmas.

\BL
For $|y|\leq \F{\rho}{\varepsilon}$,
\begin{align*}
\sum_{i,j}\partial_j(S_{\tilde g}\tilde g^{ij} \partial_i\tilde u_0)=O\left(\F{\mu^2\varepsilon^{20}}{(1+|y-\xi'|)^{N-20}}\right).
\end{align*}
For $|y|\geq \F{\rho}{\varepsilon}$,
\begin{align*}
\sum_{i,j}\partial_j(S_{\tilde g}\tilde g^{ij} \partial_i\tilde u_0)= O\left(\F{\alpha\varepsilon^2}{(1+|y-\xi'|)^{N-2}}\right).
\end{align*}
\EL

\BL
For $|y|\leq \F{\rho}{\varepsilon}$,
\begin{align*}
\sum_{i,j}\partial_j(\mathcal{\widetilde R}^{ij}\partial_i\tilde u_0)=&\ -\F{1}{2}(\partial_{jmm}\tilde h_{ij})(\partial_i\tilde u_0)-\F{1}{2}(\partial_{mm}\tilde h_{ij})(\partial_{ij}\tilde u_0)\\
&\ +O\left(\F{\mu^2\varepsilon^{20}}{(1+|y-\xi'|)^{N-20}}\right).
\end{align*}
For $|y|\geq \F{\rho}{\varepsilon}$,
\begin{align*}
\sum_{i,j}\partial_j(\mathcal{\widetilde R}^{ij}\partial_i\tilde u_0)= O\left(\F{\alpha\varepsilon^2}{(1+|y-\xi'|)^{N-2}}\right).
\end{align*}
\EL

Combining the above results, we have the following estimate for $R(y)$.

\BP\label{p4.3}
It holds
\begin{align*}
R(y)\leq
\begin{cases}
\D C\F{\mu\varepsilon^{10}}{(1+|y-\xi'|)^{N-10}} \qquad &\D \text{for }|y|\leq \F{\rho}{\varepsilon},\medskip\\
\D C\F{\alpha\varepsilon}{(1+|y-\xi'|)^{N-1}}  &\D \text{for }\F{\rho}{\varepsilon}\leq |y|\leq \F{1}{\varepsilon} ,\medskip \\
\D 0 & \D \text{for }|y|\geq \F{1}{\varepsilon}.
\end{cases}
\end{align*}
\EP

Let us consider the energy functional $E_{\tilde g}(v)$ associated
to Problem (\ref{b}), namely
\begin{align*}
E_{\tilde g}(v)=&\ \F{1}{2}\int_{\mathbb R^N} (\Delta_{\tilde g} v)^2 +\sum_{i,j}(a_N S_{\tilde g} {\tilde g}^{ij}+b_N \widetilde{\mathcal R}_{\tilde g}^{ij})\partial_i v\partial_j v +\F{N-4}{2}Q_{\tilde g} v^2 \mathrm{d}y\\
&\ -\F{(N-4)^2}{4N}\int_{\mathbb R^N} v^\F{2N}{N-4}\mathrm{d}y.
\end{align*}
In what follows, we will calculate the energy $E_{\tilde g}(\tilde u_0)$,
 which is a important step  for the existence of the solutions of our equation. First we have
\BL\label{l7.2}
For $|y|\leq \rho/\varepsilon$,
\begin{align}\label{6.2}
&\ (\Delta_{\tilde g} \tilde u_0)^2-(\Delta \tilde u_0)^2 \nonumber\\
=&\ \sum_{i,j,k,\ell}\tilde h_{i\ell}(\partial_i \tilde h_{j\ell})(\partial_{kk}\tilde u_0)(\partial_j \tilde u_0)-\sum_{i,j,k}2 \tilde h_{ij}(\partial_{kk}\tilde u_0)(\partial_{ij}\tilde u_0)\nonumber\\
&\ +\sum_{i,j,k,\ell} \tilde h_{i\ell} \tilde h_{j\ell}(\partial_{kk}\tilde u_0)(\partial_{ij}\tilde u_0)+\left(\sum_{i,j} \tilde h_{ij}(\partial_{ij}\tilde u_0)\right)^2\nonumber\\
&\ +O(|\tilde h||\partial \tilde h|^2)|\partial \tilde u_0|^2 + O(|\tilde h|^2|\partial \tilde h|)|\partial \tilde u_0||\partial^2 \tilde u_0| +O(|\tilde h|^3)|\partial^2 \tilde u_0|^2.
\end{align}
While for $|y|\geq \rho/\varepsilon$, we have
\begin{align}\label{6.3}
(\Delta_{\tilde g} \tilde u_0)^2-(\Delta \tilde u_0)^2=O(\alpha(1+|y-\xi'|)^{4-2N}).
\end{align}
\EL

\begin{proof}
It is easy to check that
\begin{align*}
(\Delta_{\tilde g} \tilde u_0)^2-(\Delta \tilde u_0)^2 =&\ \left(\sum_{i,j} (\partial_i \tilde g^{ij}) (\partial_j \tilde u_0)\right)^2 + \sum_{i,j,k,\ell}2(\partial_k \tilde g^{k\ell})\tilde g^{ij}(\partial_\ell \tilde u_0)(\partial_{ij}\tilde u_0)\\
&\ +\sum_{i,j,k,\ell}(\tilde g^{ij}-\delta_{ij})(\tilde g^{k\ell}+\delta_{k\ell})(\partial_{ij}\tilde u_0) (\partial_{k\ell}\tilde u_0).
\end{align*}
By direct calculation, we have
\begin{align*}
&\ \left(\sum_{i,j} \partial_i \tilde g^{ij} \partial_j \tilde u_0\right)^2 \\
=&\ \sum_{i,j,k,\ell}(\partial_k \tilde h_{ik})(\partial_\ell \tilde h_{j\ell})\partial_i \tilde u_0 \partial_j \tilde u_0 + O(|\tilde h||\partial \tilde h|^2)|\partial\tilde u_0|^2,
\end{align*}
\begin{align*}
&\ \sum_{i,j,k,\ell}2\tilde g^{ij}(\partial_k \tilde g^{k\ell})\partial_\ell \tilde u_0\partial_{ij}\tilde u_0 \\
=&\ -\sum_{i,k,\ell}2(\partial_k \tilde h_{k\ell})\partial_\ell \tilde u_0 \partial_{ii}\tilde u_0+\sum_{i,k,\ell,m}\Big[(\partial_k\tilde h_{km})\tilde h_{k\ell}+(\partial_k \tilde h_{m\ell})\tilde h_{km}\Big]\partial_\ell \tilde u_0 \partial_{ii}\tilde u_0\\
&\ +\sum_{i,j,k,\ell}2\tilde h_{ij}(\partial_k \tilde h_{k\ell})\partial_\ell \tilde u_0 \partial_{ij}\tilde u + O(|\tilde h|^2|\partial \tilde h|)|\partial \tilde u_0||\partial^2 \tilde u_0|
\end{align*}
and
\begin{align*}
&\ \sum_{i,j,k,\ell}(\tilde g^{ij}-\delta_{ij})(\tilde g^{k\ell}+\delta_{k\ell})(\partial_{ij}\tilde u_0) (\partial_{k\ell}\tilde u_0) \\
=&\ -\sum_{i,j,k}2\tilde h_{ij}\partial_{ij}\tilde u_0\partial_{kk}\tilde u_0+\sum_{i,j,k,m}\tilde h_{im}\tilde h_{mj}\partial_{ij}\tilde u_0\partial_{kk}\tilde u_0\\
&\ +\left(\sum_{i,j} \tilde h_{ij}\partial_{ij}\tilde u_0\right)^2 +O(|\tilde h|^3)|\partial^2 \tilde u_0|^2.
\end{align*}
Therefore, for $|y|\leq \rho/\varepsilon$, $\sum_i \partial_i \tilde h_{ij}=0$ yields (\ref{6.2}).
Since $\tilde h=0$ for $|y|\geq 1/\varepsilon$,
(\ref{6.3}) can be easily gotten. This concludes the proof.
\end{proof}

\BL\label{l7.3}
It holds
\begin{align*}
&\ \int_{\mathbb R^N}(\Delta_{\tilde g} \tilde u_0)^2-(\Delta \tilde u_0)^2\\
=&\ -\int_{B_\F{\rho}{\varepsilon}}\sum_{i,j,k,\ell}\tilde h_{i\ell} \tilde h_{j\ell}(\partial_{ikk}\tilde u_0)(\partial_{j}\tilde u_0)+\int_{B_\F{\rho}{\varepsilon}}\left(\sum_{i,j} \tilde h_{ij}(\partial_{ij}\tilde u_0)\right)^2\\
&\ +O(\mu^3\varepsilon^\F{20N}{N-1})+O\left(\alpha \Big(\F{\varepsilon}{\rho}\Big)^{N-4}\right).
\end{align*}
\EL

\begin{proof}
Since for $|y|\leq \rho/\varepsilon$,
\begin{align*}
&\ |\tilde h||\partial \tilde h|^2|\partial \tilde u_0|^2 + |\tilde h|^2|\partial \tilde h||\partial \tilde u_0||\partial^2 \tilde u_0| +|\tilde h|^3|\partial^2 \tilde u_0|^2\\
\leq &\ C\mu^3\varepsilon^{30}(1+|y-\xi'|)^{34-2N}\\
\leq &\ C\mu^3 \varepsilon^\F{20N}{N-1} (1+|y-\xi'|)^{\F{20N}{N-1}+4-2N},
\end{align*}
from Lemma \ref{l7.2} we have
\begin{align}\label{7.3}
&\ \int_{\mathbb R^N}(\Delta_{\tilde g} \tilde u_0)^2-(\Delta \tilde u_0)^2 \nonumber\\
=&\ \int_{B_\F{\rho}{\varepsilon}}\sum_{i,j,k,\ell}\tilde h_{i\ell}(\partial_i \tilde h_{j\ell})(\partial_{kk}\tilde u_0)(\partial_j \tilde u_0)-\sum_{i,j,k}2 \tilde h_{ij}(\partial_{kk}\tilde u_0)(\partial_{ij}\tilde u_0)\nonumber\\
&\ +\sum_{i,j,k,\ell} \tilde h_{i\ell} \tilde h_{j\ell}(\partial_{kk}\tilde u_0)(\partial_{ij}\tilde u_0)+\left(\sum_{i,j} \tilde h_{ij}(\partial_{ij}\tilde u_0)\right)^2\mathrm{d}y\nonumber\\
&\ +O(\mu^3\varepsilon^\F{20N}{N-1})+O\left(\alpha \Big(\F{\varepsilon}{\rho}\Big)^{N-4}\right).
\end{align}
On the other hand,  integrating by parts and using $\sum_i \partial_i \tilde h_{ij}=0$ for $|y|\leq \rho/\varepsilon$ and $\tilde h(y)=0$ for $|y|\geq 1/\varepsilon$, we know
\begin{align}\label{7.1}
&\ \sum_{i,j,k,\ell}\int_{B_\F{\rho}{\varepsilon}}\tilde h_{i\ell}(\partial_i \tilde h_{j\ell})(\partial_{kk}\tilde u_0)(\partial_j \tilde u_0)+\tilde h_{i\ell} \tilde h_{j\ell}(\partial_{kk}\tilde u_0)(\partial_{ij}\tilde u_0)\nonumber\\
=&\ \sum_{i,j,k,\ell}\int_{\mathbb R^N}\partial_i (\tilde h_{i\ell}\tilde h_{j\ell})(\partial_{kk}\tilde u_0)(\partial_j \tilde u_0)+\tilde h_{i\ell} \tilde h_{j\ell}(\partial_{kk}\tilde u_0)(\partial_{ij}\tilde u_0)+O\left(\alpha\Big(\F{\varepsilon}{\rho}\Big)^{N-4}\right)\nonumber\\
=&\ -\sum_{i,j,k,\ell}\int_{\mathbb R^N}\tilde h_{i\ell} \tilde h_{j\ell}(\partial_{ikk}\tilde u_0)(\partial_{j}\tilde u_0)+O\left(\alpha\Big(\F{\varepsilon}{\rho}\Big)^{N-4}\right)\nonumber\\
=&\ -\sum_{i,j,k,\ell}\int_{B_\F{\rho}{\varepsilon}}\tilde h_{i\ell} \tilde h_{j\ell}(\partial_{ikk}\tilde u_0)(\partial_{j}\tilde u_0)+O\left(\alpha\Big(\F{\varepsilon}{\rho}\Big)^{N-4}\right).
\end{align}
Next, direct computation shows
\begin{align*}
\tilde u_0(\partial_{ijkk} \tilde u_0)=&\ \frac{N}{(N-3) \left(N^2-4 N+8\right)}(\partial_{ijkk} \tilde u_0^2)+\frac{N^2+4 N}{N^2-4 N+8}(\partial_{ij} \tilde u_0)(\partial_{kk}\tilde u_0)\\
&\ +\frac{4 (N-4)^2 (N-2) N}{N^2-4 N+8} \tilde u_0^2\F{|y-\xi'|^2\delta_{ij}}{(\varepsilon^2+|y-\xi'|^2)^3}\\
&\ -\frac{4 (N-4)^2 \left(N^2-2\right)}{N^2-4 N+8}\tilde u_0^2\F{\delta_{ij}}{(\varepsilon^2+|y-\xi'|^2)^2}.
\end{align*}
Recalling $\tilde h$ is divergence-free, we get
\begin{align*}
&\ \sum_{i,j,k}\int_{B_\F{\rho}{\varepsilon}}\tilde h_{ij}(\partial_{kk}\tilde u_0)(\partial_{ij}\tilde u_0)\\
=&\ \sum_{i,j,k}\int_{\mathbb R^N}\partial_{ij}(\tilde h_{ij}\tilde u_0)(\partial_{kk}\tilde u_0) + O\left(\alpha\Big(\F{\varepsilon}{\rho}\Big)^{N-4}\right)\\
=&\ \sum_{i,j,k}\int_{\mathbb R^N}\tilde h_{ij}\tilde u_0(\partial_{ijkk}\tilde u_0) + O\left(\alpha\Big(\F{\varepsilon}{\rho}\Big)^{N-4}\right).
\end{align*}
Since $\sum_{i=1}^N \partial_i \tilde h_{ij}=0$ for $|y|\leq\F{\rho}{\varepsilon}$, it follows that
\begin{align*}
&\ \int_{\mathbb R^N}\tilde h_{ij}(\partial_{ijkk} \tilde u_0^2) =\int_{\mathbb R^N}(\partial_{ijkk}\tilde h_{ij}) \tilde u_0^2\\
=&\ \int_{\mathbb R^N\setminus B_{\F{\rho}{\varepsilon}}}(\partial_{ijkk}\tilde h_{ij}) \tilde u_0^2 =O\left(\alpha\Big(\F{\varepsilon}{\rho}\Big)^{N-4}\right).
\end{align*}
Thus
\begin{align}\label{7.2}
\sum_{i,j,k}\int_{B_\F{\rho}{\varepsilon}}\tilde h_{ij}(\partial_{kk}\tilde u_0)(\partial_{ij}\tilde u_0)=O\left(\alpha\Big(\F{\varepsilon}{\rho}\Big)^{N-4}\right).
\end{align}
The proof of the lemma is completed  by (\ref{7.3}), (\ref{7.1}) and (\ref{7.2}).
\end{proof}

\BL\label{c7.3}
For $|y|\leq \rho/\varepsilon$,
\begin{align*}
\sum_{i,j}S_{\tilde g} \tilde g^{ij}\partial_i \tilde u_0\partial_j\tilde u_0=-\F{1}{4}\sum_{i,k,\ell,m}(\partial_\ell\tilde h_{mk})^2 (\partial_i\tilde u_0)^2 +O(|\tilde h||\partial\tilde h|^2)|\partial\tilde u_0|^2.
\end{align*}
For $|y|\geq \rho/\varepsilon$,
\begin{align*}
\sum_{i,j}S_{\tilde g} \tilde g^{ij}\partial_i \tilde u_0\partial_j\tilde u_0=O(\alpha\varepsilon^2 (1+|y-\xi'|)^{6-2N}).
\end{align*}
\EL
\begin{proof}
Recalling that $g=e^h$,  this lemma is an easy consequence  of Lemma \ref{l7.5}.
\end{proof}

\BL
We have
\begin{align*}
&\ \int_{\mathbb R^N} \sum_{i,j} a_N S_{\tilde g} \tilde g^{ij} \partial_i\tilde u_0 \partial_j\tilde u_0\\
=&\ -\F{a_N}{4}\int_{B_\F{\rho}{\varepsilon}} \sum_{i,k,\ell,m}(\partial_\ell \tilde h_{mk})^2(\partial_i\tilde u_0)^2
+O(\mu^3\varepsilon^\F{20N}{N-1}) +O\left(\alpha\Big(\F{\varepsilon}{\rho}\Big)^{N-4}\right).
\end{align*}
\EL

\begin{proof}
This follows from Lemma \ref{c7.3} by direct calculation.
\end{proof}

\BL
It holds
\begin{align*}
&\ \int_{\mathbb R^N} \sum_{i,j} b_N \widetilde{\mathcal R}^{ij}\partial_i\tilde u_0 \partial_j\tilde u_0\\
=&\ -\F{b_N}{4}\int_{B_\F{\rho}{\varepsilon}}\sum_{i,j,m,s}(\partial_j\tilde h_{ms})(\partial_i\tilde h_{sm})(\partial_i\tilde u_0)(\partial_j\tilde u_0)\\
&\ -\F{b_N}{2}\int_{B_\F{\rho}{\varepsilon}} \sum_{i,j,m,s}\Big[\tilde h_{ms}(\partial_{s}\tilde h_{ij})-\tilde h_{si}(\partial_s\tilde h_{mj}) +\tilde h_{sj}(\partial_i\tilde h_{ms})\\
&\hspace{15em} -\tilde h_{ms}(\partial_{i}\tilde h_{sj})\Big]\partial_m (\partial_i\tilde u_0\partial_j\tilde u_0)\\
&\ +\F{b_N}{2}\int_{B_\F{\rho}{\varepsilon}} \sum_{i,j,m,s}\tilde h_{is}(\partial_{mm}\tilde h_{js})(\partial_i\tilde u_0)(\partial_j\tilde u_0)\\
&\ +O(\mu^3\varepsilon^\F{20N}{N-1}) +O\left(\alpha\Big(\F{\varepsilon}{\rho}\Big)^{N-4}\right).
\end{align*}
\EL

\begin{proof}
From Corollary \ref{c7.1}, we have
\begin{align*}
&\ \int_{\mathbb R^N} \sum_{i,j}b_N \widetilde{\mathcal R}^{ij}\partial_i\tilde u_0 \partial_j\tilde u_0\\
=&\ -\F{b_N}{2}\int_{B_\F{\rho}{\varepsilon}}\sum_{i,j,m}(\partial_{mm}\tilde h_{ij})\partial_i\tilde u_0\partial_j\tilde u_0 \\
&\ -\F{b_N}{4}\int_{B_\F{\rho}{\varepsilon}}\sum_{i,j,m,s}(\partial_j\tilde  h_{ms})(\partial_i\tilde  h_{sm})(\partial_i u_0)(\partial_j u_0)\\
&\ +\F{b_N}{2}\int_{B_\F{\rho}{\varepsilon}}\sum_{i,j,m,s}\Big[\tilde h_{ms}(\partial_{ms}\tilde h_{ij})-(\partial_s\tilde  h_{mj})(\partial_m\tilde  h_{si})+(\partial_i\tilde  h_{ms})(\partial_m\tilde  h_{sj})\\
&\hspace{17em}-\tilde h_{ms}(\partial_{mi}\tilde h_{sj})\Big](\partial_i\tilde  u_0)(\partial_j\tilde  u_0)\\
&\ +\F{b_N}{2}\int_{B_\F{\rho}{\varepsilon}}\sum_{i,j,m,s}\tilde  h_{is}(\partial_{mm}\tilde h_{js})(\partial_i\tilde  u_0)(\partial_j\tilde  u_0)\\
&\ +O(\mu^3\varepsilon^\F{20N}{N-1}) +O\left(\alpha\Big(\F{\varepsilon}{\rho}\Big)^{N-4}\right).
\end{align*}
Since
\begin{align*}
\partial_i\tilde u_0 \partial_j \tilde u_0 - \F{(N-4)}{4(N-3)}\partial_{ij} \tilde u_0^2
=\F{(N-4)^2}{2(N-3)} \tilde u_0^2 \F{\delta_{ij}}{\lambda'^2+|y-\xi'|^2},
\end{align*}
it is easy to check, noting that $\sum_{i=1}^N \partial_i \tilde h_{ij}=0$,
\begin{align*}
\left|\sum_{i,j,m}\int_{B_\F{\rho}{\varepsilon}}(\partial_{mm}\tilde h_{ij})\partial_i\tilde u_0\partial_j\tilde u_0\right|
&\leq C\int_{\mathbb R^n\setminus B_\F{\rho}{\varepsilon}}|\partial^4\tilde h|\tilde u_0^2 +O\left(\alpha\rho^2\Big(\F{\varepsilon}{\rho}\Big)^{N-4}\right) \\
&= O\left(\alpha\rho^2\Big(\F{\varepsilon}{\rho}\Big)^{N-4}\right).
\end{align*}
Integrating by parts, we have
\begin{align*}
&\ \sum_{i,j,m,s}\int_{B_\F{\rho}{\varepsilon}}\tilde h_{ms} (\partial_{ms}\tilde h_{ij})(\partial_i\tilde u_0)(\partial_j\tilde u_0)\\
=&\ \sum_{i,j,m,s}\int_{\mathbb R^n} \partial_m [\tilde h_{ms} (\partial_{s}\tilde h_{ij})]\partial_i\tilde u_0\partial_j\tilde u_0+O(\alpha^2\rho^2(\F{\varepsilon}{\rho})^{N-4})\\
=&\ -\sum_{i,j,m,s}\int_{B_\F{\rho}{\varepsilon}}\tilde h_{ms} (\partial_s\tilde h_{ij}) \partial_{m}(\partial_i\tilde u_0\partial_j\tilde u_0) +O(\alpha^2\rho(\F{\varepsilon}{\rho})^{N-4}),
\end{align*}
\begin{align*}
&\ \sum_{i,j,m,s}\int_{B_\F{\rho}{\varepsilon}} (\partial_s\tilde h_{jm}) (\partial_{m}\tilde h_{is})(\partial_i \tilde u_0)(\partial_j \tilde u_0)\\
=&\ \sum_{i,j,m,s}\int_{\mathbb R^n} \partial_{m}[\tilde h_{is} (\partial_s \tilde h_{jm})] \partial_i\tilde  u_0\partial_j \tilde u_0+O(\alpha^2\rho^2(\F{\varepsilon}{\rho})^{N-4})\\
=&\ -\sum_{i,j,m,s}\int_{B_\F{\rho}{\varepsilon}} \tilde h_{is} (\partial_s \tilde h_{jm})  \partial_{m}(\partial_i \tilde u_0\partial_j \tilde u_0) +O(\alpha\rho(\F{\varepsilon}{\rho})^{N-4}),
\end{align*}
\begin{align*}
&\ \sum_{i,j,m,s}\int_{B_\F{\rho}{\varepsilon}} (\partial_i \tilde h_{ms})(\partial_m \tilde h_{sj}) (\partial_i \tilde u_0)(\partial_j \tilde u_0)\\
=&\ -\sum_{i,j,m,s}\int_{B_\F{\rho}{\varepsilon}} \tilde h_{sj}  (\partial_i \tilde h_{ms}) (\partial_m(\partial_i \tilde u_0)(\partial_j \tilde u_0))+O(\alpha\rho(\F{\varepsilon}{\rho})^{N-4})
\end{align*}
and
\begin{align*}
&\ \sum_{i,j,m,s}\int_{B_\F{\rho}{\varepsilon}} \tilde h_{ms}(\partial_{mi}\tilde h_{sj}) (\partial_i \tilde u_0)(\partial_j \tilde u_0)\\
=&\ -\sum_{i,j,m,s}\int_{B_\F{\rho}{\varepsilon}} \tilde h_{ms}  (\partial_i \tilde h_{js}) (\partial_m(\partial_i \tilde u_0)(\partial_j \tilde u_0))+O(\alpha\rho(\F{\varepsilon}{\rho})^{N-4}).
\end{align*}
The proof is complete.
\end{proof}

\BL
\begin{align*}
&\ \F{N-4}{2}\int_{\mathbb R^N}Q_{\tilde g} \tilde u_0^2\\
=&\ \F{N-4}{8(N-1)}\int_{B_\F{\rho}{\varepsilon}} \sum_{i,k,\ell,m}\Big[(\partial_{i\ell}\tilde h_{mk})^2+(\partial_\ell\tilde h_{mk})(\partial_{ii\ell}\tilde h_{mk})\Big]\tilde u_0^2\\
&\ -\F{N-4}{4(N-2)^2}\int_{B_\F{\rho}{\varepsilon}} \sum_{i,j,m,s}(\partial_{mm}\tilde h_{ij})(\partial_{ss}\tilde h_{ij})u_0^2\\
& +O\left(\mu^3 \varepsilon^\F{20N}{N-1}\right)+O\left(\alpha\rho^4\left(\F{\varepsilon}{\rho}\right)^{N-4}\right).
\end{align*}
\EL

\begin{proof}
This is an easy consequence of Lemma \ref{l7.4}.
\end{proof}

Now we have the following estimate of $E_{\tilde g}(\tilde u_0)$.

\BP\label{p6.1}
\begin{align*}
2E_{\tilde g}(\tilde u_0)=&\ 2E-\int_{B_\F{\rho}{\varepsilon}}\sum_{i,j,k,\ell}\tilde h_{i\ell} \tilde h_{j\ell}(\partial_{ikk}\tilde u_0)(\partial_{j}\tilde u_0)+\int_{B_\F{\rho}{\varepsilon}}\Big(\sum_{i,j} \tilde h_{ij}(\partial_{ij}\tilde u_0)\Big)^2\\
&\ -\F{a_N}{4}\int_{B_\F{\rho}{\varepsilon}} \sum_{i,k,\ell,m}(\partial_\ell \tilde h_{mk})^2(\partial_i\tilde u_0)^2\\
&\ -\F{b_N}{4}\int_{B_\F{\rho}{\varepsilon}}\sum_{i,j,m,s}(\partial_j\tilde h_{ms})(\partial_i\tilde h_{sm})(\partial_i\tilde u_0)(\partial_j\tilde u_0)\\
&\ -\F{b_N}{2}\int_{B_\F{\rho}{\varepsilon}} \sum_{i,j,m,s}\Big[\tilde h_{ms}(\partial_{s}\tilde h_{ij})-\tilde h_{si}(\partial_s\tilde h_{mj}) +\tilde h_{sj}(\partial_i\tilde h_{ms})\\
&\hspace{15em} -\tilde h_{ms}(\partial_{i}\tilde h_{sj})\Big]\partial_m (\partial_i\tilde u_0\partial_j\tilde u_0)\\
&\ +\F{b_N}{2}\int_{B_\F{\rho}{\varepsilon}} \sum_{i,j,m,s}\tilde h_{is}(\partial_{mm}\tilde h_{js})(\partial_i\tilde u_0)(\partial_j\tilde u_0)\\
&\ +\F{N-4}{8(N-1)}\int_{B_\F{\rho}{\varepsilon}} \sum_{i,k,\ell,m}\Big[(\partial_{i\ell}\tilde h_{mk})^2+(\partial_\ell\tilde h_{mk})(\partial_{ii\ell}\tilde h_{mk})\Big]\tilde u_0^2\\
&\ -\F{N-4}{4(N-2)^2}\int_{B_\F{\rho}{\varepsilon}} \sum_{i,j,m,s}(\partial_{mm}\tilde h_{ij})(\partial_{ss}\tilde h_{ij})\tilde u_0^2\\
& +O\left(\mu^3 \varepsilon^\F{20N}{N-1}\right)+O\left(\alpha\big(\F{\varepsilon}{\rho}\big)^{N-4}\right),
\end{align*}
where $E$ is the constant such that
\begin{align*}
E= \F{N-4}{N}\int_{\mathbb R^N} \left(\F{1}{1+|y|^2}\right)^N \mathrm{d}y.
\end{align*}

\noindent
{\bf Remark:} Note that
\begin{equation}
E= P({\mathbb S}^N).
\end{equation}
\EP

\section{Linearized Operator}

In this section we develop the invertibility theory  for the linearized operator $P_{\tilde g}-\F{N+4}{2}\tilde u_0^\F{8}{N-4}$ in suitable weighted $L^\infty$ spaces.

We define two norms
\begin{gather*}
\|\phi\|_* =\sup_{y\in\mathbb R^N} \sum_{i=0}^2 \left[\F{1}{\F{\mu \varepsilon^{10}}{(1+|y-\xi'|)^{N-14+i}}+\alpha(\F{\varepsilon}{\rho})^{N-4+i}}+\F{(1+|y-\xi'|)^{N-4+i}}{\alpha}\right]|\partial^i\phi(y)|,\\
\|\zeta\|_{**}  =\sup_{y\in\mathbb R^N} \Bigg[\F{\chi_{\{|y-\xi'|\leq \F{\rho}{\varepsilon}\}}(1+|y-\xi'|)^{N-10}}{\mu \varepsilon^{10}}
+\F{\chi_{\{\F{\rho}{\varepsilon}\leq |y-\xi'|\leq\F{1}{\varepsilon} \}}(1+|y-\xi'|)^{N-1}}{\alpha\varepsilon}\\
+\F{ \chi_{\{|y-\xi'|\geq\F{1}{\varepsilon} \}}(1+|y-\xi'|)^{N+\sigma}}{\alpha}\Bigg]|\zeta(y)|,
\end{gather*}
where $\chi_S$ is the characteristic function on the set $S$, and $0<\sigma<1$ is a small constant.

Denote
\begin{align*}
Z_0=\F{\partial \tilde u_0}{\partial \lambda'},\qquad Z_j=\F{\partial \tilde u_0}{\partial \xi'_j}\quad j=1,\cdots,N.
\end{align*}

First, we consider the following problem. Given $\zeta\in C^\alpha(\mathbb R^N)$, find a function $\phi$ such that
for certain constants $c_{i}$, $i=0,1,\cdots,N$,
\begin{align}\label{first}
\begin{cases}
\D P_{\tilde g}\phi-\F{N+4}{2}\tilde u_0^\F{8}{N-4}\phi=\zeta+\sum_{i=0}^N c_i \chi Z_i \qquad \text{in }\mathbb R^N,\\
\D \int_{\mathbb R^N} \phi \chi Z_i =0,
\end{cases}
\end{align}
where $\chi(y)=\chi(|y-\xi'|)$ is a cut-off function satisfying $\chi(y)=1$ for $|y-\xi'|\leq r_0$, $\chi(y)=0$ for $|y-\xi'|\geq r_0+1$.
Here $r_0>0$ is large but fixed.

\BP\label{p3.1}
Assume $N\geq 18$,  $(\xi',\lambda')\in \Lambda$ and $\alpha$ is small and fixed. Then for small $\varepsilon$, there is a unique solution $\phi$ to (\ref{first}). Moreover
\begin{align*}
\|\phi\|_* \leq C\|\zeta\|_{**}
\end{align*}
where $C$ is independent of $\alpha$ and $\varepsilon$.
\EP

To prove the above proposition, we need the following priori estimate.

\BL\label{l3.2}
Under the assumptions of Proposition \ref{p3.1}, for any solution $\phi$ to (\ref{first}), there exists a constant $C$ such that
\begin{align*}
\|\phi\|_* \leq C\|\zeta\|_{**}.
\end{align*}
\EL

\begin{proof}
We use the  contradiction argument as in \cite{DFM} and \cite{RW}. Assume there are sequences $\varepsilon_n\to 0$ and the corresponding $\zeta_n$, $\phi_n$ such that
$\|\zeta_n\|_{**}\to 0$ but $\|\phi_n\|_* =1$. For abbreviation, we omit the subscript $n$ in the following proof.
Testing the equation against $\bar\chi Z_j$ and integrating by parts four times,
where $\bar\chi(y)$ is a smooth cut-off function satisfying $\bar\chi(y)=1$ for $|y-\xi'|\leq \F{\rho}{4\varepsilon}$,
$\bar\chi(y)=0$ for $|y-\xi'|\geq \F{\rho}{2\varepsilon}$ and $|\nabla^i\bar\chi|\leq C(\F{\varepsilon}{\rho})^i$, $1\leq i\leq 4$.
we get
\begin{align*}
\sum_i c_i \int_{\mathbb R^N}  \chi Z_i Z_j = \int_{\mathbb R^N} \left(P_{\tilde g_n}(\bar\chi Z_j)-\F{N+4}{2}\tilde u_0^\F{8}{N-4}\bar\chi Z_j\right)\phi -\int_{\mathbb R^N} \zeta \bar\chi Z_j.
\end{align*}
This defines a linear system in the $c_i$ which is ``almost diagonal'' as $\varepsilon$ approaches zero,
since we have
\begin{align*}
\int_{\mathbb R^N} \chi Z_0 Z_j &= \delta_{0j} \int \chi \left(\F{\partial \tilde u_0}{\partial \lambda'}\right)^2,\\
\int_{\mathbb R^N} \chi Z_i Z_j &= \delta_{ij} \int \chi\left(\F{\partial \tilde u_0}{\partial y_j}\right)^2 \qquad \forall~i=1,\cdots,N.
\end{align*}
On the other hand, using $\Delta^2 Z_j-\F{N+4}{2}\tilde u_0^\F{8}{N-4}Z_j=0$ and the estimates of $S_{\tilde g}$, $\widetilde{\mathcal R}^{st}$, $Q_{\tilde g}$
in the previous section, we have
\begin{align*}
&\ \int_{\mathbb R^N} \left(P_{\tilde g}(\bar\chi Z_j)-\F{N+4}{2}\tilde u_0^\F{8}{N-4}\bar\chi Z_j\right)\phi \\
= &\ \int_{\mathbb R^N} \Bigg\{\Delta^2_{\tilde g} (\bar\chi Z_j)-\bar\chi\Delta^2 ( Z_j)-\sum_{s,t}\partial_s\left[\left(a_N S_{\tilde g} {\tilde g}^{st}+b_N\widetilde{\mathcal R}^{st}\right)\partial_t(\bar\chi Z_j)\right]\\
&\qquad\qquad\qquad +\F{N-4}{2}Q_{\tilde g} \bar\chi Z_j\Bigg\}\phi\\
=&\ o(\mu\varepsilon^{10})\|\phi\|_*.
\end{align*}
It is also easy to get
\begin{align*}
&\ \int_{\mathbb R^N} \zeta \bar\chi Z_j \leq C\mu\varepsilon^{10}\|\zeta\|_{**}.
\end{align*}
Thus we conclude
\begin{align*}
|c_i| \leq&\ o(\mu\varepsilon^{10})\|\phi\|_*+C\mu\varepsilon^{10}\|\zeta\|_{**} \qquad \forall~ i=0,\cdots, N,
\end{align*}
so $c_i=o(\mu\varepsilon^{10})$.\par

Next we claim that, for any fixed $R>0$,
\begin{align*}
\|\phi\|_{L^\infty(B_R(\xi'))}=o(\mu\varepsilon^{10}).
\end{align*}
Indeed, by elliptic regularity we can get a $\hat\phi$ such that $\F{\phi}{\mu\varepsilon^{10}}\to\hat\phi$ in $C^4_\text{loc}(\mathbb R^N)$ and
\begin{align*}
\Delta^2 \hat\phi - \F{N+4}{2}\tilde u_0^\F{8}{N-4}\hat\phi=0 \qquad \text{in }\mathbb R^N.
\end{align*}
This implies $\hat\phi$ is a linear combination of the functions $Z_j$, $j=0,1,\cdots,N$, see \cite{LW}.
On the other hand, the assumed orthogonality conditions on $\phi$
yields $\int_{\mathbb R^N}\hat\phi \chi Z_j=0$ for all $j$. Hence $\hat\phi\equiv 0$, which concludes the claim.\par

Now rewrite the equation in the following form
\begin{align}\label{potential}
\phi(y) =&\ \int_{\mathbb R^N} G(y,z) \sum_{i,j} \partial_j\left[(a_N S_{\tilde g} \tilde g^{ij} +b_N \widetilde{\mathcal R}^{ij})\partial_i \phi\right](z)\mathrm{d}z\nonumber\\
&\ -\int_{\mathbb R^N} \F{N-4}{2} G(y,z) Q_{\tilde g}(z) \phi(z)\mathrm{d}z + \int_{\mathbb R^N}G(y,z) \F{N+4}{2}\tilde u_0^\F{8}{N-4}\phi(z)\mathrm{d}z\nonumber\\
&\ +\int_{\mathbb R^N} G(y,z) \zeta(z)\mathrm{d}z + \sum_i c_i \int_{\mathbb R^N} G(y,z) \chi Z_i(z)\mathrm{d}z.
\end{align}
We make now the following observations:\par

Owing to $S_{\tilde g}(z)=O(\alpha \varepsilon^2)$, $\partial S_{\tilde g}=O(\alpha \varepsilon^3)$ and $\widetilde{\mathcal R}^{ij}(z)=O(\alpha \varepsilon^2)$, $\partial \widetilde{\mathcal R}^{ij}(z)=O(\alpha \varepsilon^3)$,
we have
\begin{align}\label{2.13}
&\ \left|\int_{\mathbb R^N} G(y,z) \sum_{i,j} \partial_j\left[(a_N S_{\tilde g} {\tilde g}^{ij} +b_N \widetilde{\mathcal R}^{ij})\partial_i \phi\right](z)\mathrm{d}z\right|\nonumber\\
=&\ O(\alpha\varepsilon^3) \int_{B_{1/\varepsilon}} G(y,z) |\partial\phi(z)|\mathrm{d}z + O(\alpha \varepsilon^2) \int_{B_{1/\varepsilon}} G(y,z) |\partial^2\phi(z)|\mathrm{d}z\nonumber\\
:=&\ I+II.
\end{align}
For $|y-\xi'|\leq \F{\rho}{2\varepsilon}$,
\begin{align}\label{2.14}
|I|\leq &\ C\alpha\varepsilon^3 \|\phi\|_*
\Bigg\{\int_{|z|\leq\F{\rho}{\varepsilon}}\F{1}{|y-z|^{N-4}} \left[\F{\mu \varepsilon^{10}}{(1+|z-\xi'|)^{N-13}}+\alpha(\F{\varepsilon}{\rho})^{N-3}\right]\mathrm{d}z \nonumber\\
&\qquad\qquad\qquad +\int_{\F{\rho}{\varepsilon}\leq|z|\leq \F{1}{\varepsilon}}\F{1}{|y-z|^{N-4}} \F{\alpha}{(1+|z-\xi'|)^{N-3}}\mathrm{d}z\Bigg\}\nonumber\\
\leq&\ C\alpha\|\phi\|_*\left[\F{\mu \varepsilon^{10}\varepsilon^3(1+|y-\xi'|)^3}{(1+|y-\xi'|)^{N-14}}+\alpha\rho^3(\F{\varepsilon}{\rho})^{N-4}\right]\nonumber\\
\leq &\ C\alpha\rho^3\|\phi\|_* \left[\F{\mu \varepsilon^{10}}{(1+|y-\xi'|)^{N-14}}+\alpha(\F{\varepsilon}{\rho})^{N-4}\right],
\end{align}
where we use, for any $0<s , k< N$ such that $s+k<N$,
\begin{align}\label{2.11}
\int_{\mathbb R^N} \F{1}{|y-z|^{N-s}}\F{1}{(1+|z-\xi'|)^{N-k}}\mathrm{d}z \leq C (1+|y-\xi'|)^{k+s-N}.
\end{align}
The proof of (\ref{2.11}) is standard and is given in the appendix.
Similarly, for $|y-\xi'|\leq \F{\rho}{2\varepsilon}$,
\begin{align}\label{2.15}
|II| \leq C\alpha\rho^2\|\phi\|_* \left[\F{\mu \varepsilon^{10}}{(1+|y-\xi'|)^{N-14}}+\alpha(\F{\varepsilon}{\rho})^{N-4}\right].
\end{align}
On the other hand, for $|y-\xi'|\geq \F{\rho}{2\varepsilon}$,
\begin{align}\label{2.16}
|I|\leq &\ C\alpha\varepsilon^3 \|\phi\|_*
\Bigg\{\int_{|z|\leq\F{\rho}{\varepsilon}}\F{1}{|y-z|^{N-4}} \left[\F{\mu \varepsilon^{10}}{(1+|z-\xi'|)^{N-13}}+\alpha(\F{\varepsilon}{\rho})^{N-3}\right]\mathrm{d}z \nonumber\\
&\qquad\qquad\qquad +\int_{\F{\rho}{\varepsilon}\leq|z|\leq \F{1}{\varepsilon}}\F{1}{|y-z|^{N-4}} \F{\alpha}{(1+|z-\xi'|)^{N-3}}\mathrm{d}z\Bigg\}\nonumber\\
\leq&\ C\rho^3\|\phi\|_* \F{\alpha}{(1+|y-\xi'|)^{N-4}},
\end{align}
where we use, for any $0<s , k< N$ such that $s+k<N$,
\begin{align}\label{2.12}
\int_{B_r} \F{1}{|y-z|^{N-s}}\F{1}{(1+|z-\xi'|)^{N-k}}\mathrm{d}z \leq C r^k (1+|y-\xi'|)^{s-N},
\end{align}
which is a direct result of (\ref{2.11}) and the proof is also given in the appendix.
A similar proof also gives
\begin{align}\label{2.17}
|II|\leq C\rho^2\|\phi\|_* \F{\alpha}{(1+|y-\xi'|)^{N-4}} \qquad \text{for }|y-\xi'|\geq \F{\rho}{2\varepsilon}.
\end{align}
Since $Q_{\tilde g}(z)=O(\alpha \varepsilon^4)$, we similarly have
\begin{align}\label{2.18}
&\ \left|\int_{\mathbb R^N} \F{N-4}{2}G(y,z) Q_{\tilde g}(z) \phi(z)\mathrm{d}z\right|\nonumber\\
\leq &\  C\alpha\rho^4\|\phi\|_* \left[\F{\mu \varepsilon^{10}}{(1+|y-\xi'|)^{N-14}}+\alpha(\F{\varepsilon}{\rho})^{N-4}\right] \qquad \text{for }|y-\xi'|\leq \F{\rho}{2\varepsilon},
\end{align}
and for $|y-\xi'|\geq \F{\rho}{2\varepsilon}$,
\begin{align}\label{2.19}
\left|\int_{\mathbb R^N} \F{N-4}{2}G(y,z) Q_{\tilde g}(z) \phi(z)\mathrm{d}z\right|
\leq C\rho^4\|\phi\|_*  \F{\alpha}{(1+|y-\xi'|)^{N-4}} .
\end{align}
Therefore, combining (\ref{2.13}), (\ref{2.14}), (\ref{2.15}), (\ref{2.16}), (\ref{2.17})-(\ref{2.19}), we finally have
\begin{align}\label{2.23}
&\ \Bigg|\int_{\mathbb R^N} G(y,z) \sum_{i,j} \partial_i\left[(a_N S_{\tilde g} \tilde g^{ij} +b_N\widetilde{\mathcal R}^{ij})\partial_j \phi\right](z)\mathrm{d}z \nonumber\\
&\qquad\qquad\qquad\qquad -\int_{\mathbb R^N} \F{N-4}{2} G(y,z) Q_{\tilde g}(z) \phi(z)\mathrm{d}z\Bigg| \nonumber\\
\leq&\
\begin{cases}
\D C\alpha\rho^2\|\phi\|_* \left[\F{\mu \varepsilon^{10}}{(1+|y-\xi'|)^{N-14}}+\alpha(\F{\varepsilon}{\rho})^{N-4}\right] \quad & \D\text{for }|y-\xi'|\leq \F{\rho}{2\varepsilon},\bigskip\\
\D C\rho^2\|\phi\|_* \F{\alpha}{(1+|y-\xi'|)^{N-4}} & \D\text{for }|y-\xi'|\geq \F{\rho}{2\varepsilon}.
\end{cases}
\end{align}\par

Next note that
\begin{align*}
&\ \left|\int_{\mathbb R^N} G(y,z) \F{N+4}{2}\tilde u_0(z)^\F{8}{N-4}\phi(z)\mathrm{d}z\right|\\
\leq &\ C \left\{\int_{B_R(\xi')}  + \int_{R<|z-\xi'|<\F{\rho}{\varepsilon}}+ \int_{|z-\xi'|> \F{\rho}{\varepsilon}}\right\} \F{1}{|y-z|^{N-4}}\F{\phi(z)}{(1+|z-\xi'|)^8}.
\end{align*}
Since $\|\phi\|_{L^\infty(B_R(\xi'))}=o(\mu\varepsilon^{10})$, it is easy to check that
\begin{align*}
& \int_{B_R(\xi')}\F{1}{|y-z|^{N-4}}\F{\phi(z)}{(1+|z-\xi'|)^8} \mathrm{d}z\\
=&\
\begin{cases}
\D o(1) \left[\F{\mu \varepsilon^{10}}{(1+|y-\xi'|)^{N-14}}+\alpha(\F{\varepsilon}{\rho})^{N-4}\right] \quad & \D\text{for }|y-\xi'|\leq \F{\rho}{2\varepsilon},\bigskip\\
\D o(\mu\varepsilon^{10}) \F{\alpha}{(1+|y-\xi'|)^{N-4}} & \D\text{for }|y-\xi'|\geq \F{\rho}{2\varepsilon}.
\end{cases}
\end{align*}
\begin{align*}
& \int_{R<|z-\xi'|<\F{\rho}{\varepsilon}}\F{1}{|y-z|^{N-4}}\F{\phi(z)}{(1+|z-\xi'|)^8} \mathrm{d}z\\
\leq &\
\begin{cases}
\D \F{C}{R^3} \|\phi\|_* \left[\F{\mu \varepsilon^{10}}{(1+|y-\xi'|)^{N-14}}+\alpha(\F{\varepsilon}{\rho})^{N-4}\right] \quad & \D\text{for }|y-\xi'|\leq \F{\rho}{2\varepsilon},\bigskip\\
\D C\alpha(\F{\varepsilon}{\rho})^4 \|\phi\|_* \F{\alpha}{(1+|y-\xi'|)^{N-4}} & \D\text{for }|y-\xi'|\geq \F{\rho}{2\varepsilon}.
\end{cases}
\end{align*}
\begin{align*}
& \int_{|z-\xi'|\geq\F{\rho}{\varepsilon}}\F{1}{|y-z|^{N-4}}\F{\phi(z)}{(1+|z-\xi'|)^8} \mathrm{d}z\\
\leq &\
\begin{cases}
\D C(\F{\varepsilon}{\rho})^4 \|\phi\|_* \left[\F{\mu \varepsilon^{10}}{(1+|y-\xi'|)^{N-14}}+\alpha(\F{\varepsilon}{\rho})^{N-4}\right] \quad & \D\text{for }|y-\xi'|\leq \F{\rho}{2\varepsilon},\bigskip\\
\D C(\F{\varepsilon}{\rho})^3 \|\phi\|_* \F{\alpha}{(1+|y-\xi'|)^{N-4}} & \D\text{for }|y-\xi'|\geq \F{\rho}{2\varepsilon}.
\end{cases}
\end{align*}
Thus
\begin{align}\label{2.20}
&\ \left|\int_{\mathbb R^N} G(y,z) \F{N+4}{2}\tilde u_0(z)^\F{8}{N-4}\phi(z)\mathrm{d}z\right| \nonumber\\
\leq &\
\begin{cases}
\D \left(\F{C}{R^3} \|\phi\|_*+o(1)\right) \left[\F{\mu \varepsilon^{10}}{(1+|y-\xi'|)^{N-14}}+\alpha(\F{\varepsilon}{\rho})^{N-4}\right]  & \D\text{for }|y-\xi'|\leq \F{\rho}{2\varepsilon},\bigskip\\
\D C(\F{\varepsilon}{\rho})^3 \|\phi\|_* \F{\alpha}{(1+|y-\xi'|)^{N-4}} & \D\text{for }|y-\xi'|\geq \F{\rho}{2\varepsilon}.
\end{cases}
\end{align}\par

Similarly,
\begin{align*}
\int_{\mathbb R^N} G(y,z) \zeta(z)\mathrm{d}z
 =&\  \left\{\int_{|z-\xi'|<\F{\rho}{\varepsilon}}  + \int_{\F{\rho}{\varepsilon}<|z-\xi'|<\F{1}{\varepsilon}}+ \int_{|z-\xi'|> \F{1}{\varepsilon}}\right\} G(y,z) \zeta(z).
\end{align*}
Using (\ref{2.11}) and (\ref{2.12}), we have
\begin{align*}
\left|\int_{|z-\xi'|<\F{\rho}{\varepsilon}}G(y,z) \zeta(z)\right|\leq
\begin{cases}
\D C\|\zeta\|_{**}\F{\mu\varepsilon^{10}}{(1+|y-\xi'|)^{N-14}} \quad &  \D \text{for }|y-\xi'|\leq \F{\rho}{2\varepsilon},\medskip\\
\D C\|\zeta\|_{**}\F{\mu\rho^{10}}{(1+|y-\xi'|)^{N-4}}  &  \D \text{for }|y-\xi'|\geq \F{\rho}{2\varepsilon},
\end{cases}
\end{align*}
\begin{align*}
\left|\int_{\F{\rho}{\varepsilon}<|z-\xi'|<\F{1}{\varepsilon}}G(y,z) \zeta(z)\right|\leq
\begin{cases}
\D C\|\zeta\|_{**}\alpha(\F{\varepsilon}{\rho})^{N-4} \quad &  \D \text{for }|y-\xi'|\leq \F{\rho}{2\varepsilon},\medskip\\
\D C\|\zeta\|_{**}\F{\alpha}{(1+|y-\xi'|)^{N-4}}  &  \D \text{for }|y-\xi'|\geq \F{\rho}{2\varepsilon}
\end{cases}
\end{align*}
and
\begin{align*}
\left|\int_{|z-\xi'|\geq \F{1}{\varepsilon}}G(y,z) \zeta(z)\right|\leq
\begin{cases}
\D C\|\zeta\|_{**}\alpha(\F{\varepsilon}{\rho})^{N-4} \quad &  \D \text{for }|y-\xi'|\leq \F{\rho}{2\varepsilon},\medskip\\
\D C\|\zeta\|_{**}\F{\alpha}{(1+|y-\xi'|)^{N-4}}  &  \D \text{for }|y-\xi'|\geq \F{\rho}{2\varepsilon}.
\end{cases}
\end{align*}
So
\begin{align}\label{2.21}
&\ \left|\int_{\mathbb R^N} G(y,z) \zeta(z)\mathrm{d}z\right|\nonumber\\
\leq &\
\begin{cases}
\D C \|\zeta\|_{**} \left[\F{\mu \varepsilon^{10}}{(1+|y-\xi'|)^{N-14}}+\alpha(\F{\varepsilon}{\rho})^{N-4}\right] \quad & \D\text{for }|y-\xi'|\leq \F{\rho}{2\varepsilon},\bigskip\\
\D C\|\zeta\|_{**} \F{\alpha}{(1+|y-\xi'|)^{N-4}} & \D\text{for }|y-\xi'|\geq \F{\rho}{2\varepsilon}.
\end{cases}
\end{align}\par

Since we have know $c_i=o(\mu\varepsilon^{10})$, it holds
\begin{align}\label{2.22}
&\ \sum_i c_i \int_{\mathbb R^N} G(y,z) \chi(z) Z_i(z)\mathrm{d}z \nonumber\\
\leq &\ C |c_i| \int_{\mathbb R^N} \F{1}{|y-z|^{N-4}}\chi(z) \F{1}{(1+|z-\xi'|)^{N-4}}\mathrm{d}z \nonumber\\
= &\
\begin{cases}
\D o(1) \left[\F{\mu \varepsilon^{10}}{(1+|y-\xi'|)^{N-14}}+\alpha(\F{\varepsilon}{\rho})^{N-4}\right] \quad & \D\text{for }|y-\xi'|\leq \F{\rho}{2\varepsilon},\bigskip\\
\D o(\mu\varepsilon^{10}) \F{\alpha}{(1+|y-\xi'|)^{N-4}} & \D\text{for }|y-\xi'|\geq \F{\rho}{2\varepsilon}.
\end{cases}
\end{align}
Now we obtain that, by combining (\ref{2.23})-(\ref{2.22}) and choosing $R$ large enough,
\begin{align*}
\left[\F{1}{\F{\mu \varepsilon^{10}}{(1+|y-\xi'|)^{N-14}}+\alpha(\F{\varepsilon}{\rho})^{N-4}}+\F{(1+|y-\xi'|)^{N-4}}{\alpha}\right] |\phi| \leq  C\|\zeta\|_{**}+o(1).
\end{align*}\par

Taking the derivative in (\ref{potential}), we have
\begin{align*}
\partial_{y_i}\phi(y) =&\ \int_{\mathbb R^N} \partial_{y_i}G(y,z) \sum_{s,t} \partial_s\left[(a_N S_{\tilde g}\tilde g^{st} +b_N \widetilde{\mathcal R}^{st})\partial_t \phi\right](z)\mathrm{d}z\nonumber\\
&\ -\int_{\mathbb R^N} \F{N-4}{2} \partial_{y_i}G(y,z) Q_{\tilde g}(z)  \phi(z)\mathrm{d}z \\
&\ + \int_{\mathbb R^N}\partial_{y_i}G(y,z)\F{N+4}{2} \tilde u_0^\F{8}{N-4}\phi(z)\mathrm{d}z\nonumber\\
&\ +\int_{\mathbb R^N} \partial_{y_i}G(y,z) \zeta(z)\mathrm{d}z \\
&\ + \sum_i c_i \int_{\mathbb R^N} \partial_{y_i}G(y,z) Z_i(z)\mathrm{d}z.
\end{align*}
Since $|\partial_{y_i}G(y,z)| \leq C[1+O(\alpha)]|x-y|^{3-N}$, similarly we can prove that
\begin{align*}
\left[\F{1}{\F{\mu \varepsilon^{10}}{(1+|y-\xi'|)^{N-13}}+\alpha(\F{\varepsilon}{\rho})^{N-3}}+\F{(1+|y-\xi'|)^{N-3}}{\alpha}\right] |\partial \phi|\leq C\|\zeta\|_{**}+o(1).
\end{align*}
It is also similar to get that
\begin{align*}
\left[\F{1}{\F{\mu \varepsilon^{10}}{(1+|y-\xi'|)^{N-12}}+\alpha(\F{\varepsilon}{\rho})^{N-2}}+\F{(1+|y-\xi'|)^{N-2}}{\alpha}\right] |\partial^2 \phi|\leq C\|\zeta\|_{**}+o(1).
\end{align*}
So we finally have
\begin{align*}
\|\phi\|_* \leq C\|\zeta\|_{**} +o(1),
\end{align*}
which is a contradiction.
\end{proof}

\begin{proof}[Proof of Proposition \ref{p3.1}]
Consider the space
\begin{align*}
\mathcal H=\left\{\phi\in H^2(\mathbb R^N) :\ \int_{\mathbb R^N} \chi Z_j \phi =0\quad \forall ~j=0,1,\cdots,N\right\}
\end{align*}
endowed with the inner product $(\phi,\psi)=\int_{\mathbb R^N} \Delta_g\phi\Delta_g\psi$. Problem (\ref{first}) expressed
in weak form is equivalent to that of finding a $\phi\in\mathcal H$ such that, for any $\psi\in\mathcal H$,
\begin{align*}
(\phi,\psi)=&\ -\int_{\mathbb R^N}\sum_{i,j} (a_N S_{\tilde g} \tilde g^{ij}+b_N \widetilde{\mathcal R}^{ij})\partial_i\phi\partial_j\psi+\F{N-4}{2}Q_{\tilde g}\phi\psi\\
&\ + \int_{\mathbb R^N} \F{N+4}{2} \tilde u_0^\F{8}{N-4}\phi\psi + \int_{\mathbb R^N}\zeta\psi.
\end{align*}
With the aid of Riesz's representation theorem, this equation can be rewritten in $\mathcal H$ in the operator form
\begin{align*}
\phi=K(\phi)+\tilde \zeta
\end{align*}
with certain $\tilde \zeta\in\mathcal H$ which depends linearly on $\zeta$ and $K$ is a compact operator in $\mathcal H$.
Fredholm's alternative guarantees unique solvability of this problem for any
$\tilde \zeta$ provided that the homogeneous equation $ \phi= K(\phi)$ has only the zero solution
in $\mathcal H$, which is equivalent to (\ref{first}) with $\zeta = 0$. Thus existence of
a unique solution follows from Lemma \ref{l3.2}. This finishes the proof.
\end{proof}

\begin{remark}\label{r3.1}
The result of Proposition \ref{p3.1} implies that the unique solution $\phi=T(\zeta)$ of (\ref{first}) defines a continuous linear map
from the weighted $L^\infty$ space $L^\infty_{**}$, equipped with norm $\|\cdot\|_{**}$, into the weighted $L^\infty$ space $L^\infty_*$, equipped with $\|\cdot\|_*$.
\end{remark}

It is important for later purposes to understand the differentiability of the operator
$T$ with respect to the variables $\xi'$ and $\lambda'$.

\BP\label{p3.2}
Assume $(\xi',\lambda')\in\Lambda$. We have
\begin{align*}
\|\nabla_{(\xi',\lambda')}T(\zeta)\|_* \leq C\|\zeta\|_{**}.
\end{align*}
\EP

\begin{proof}
Denote formally $Z=\partial_{\xi'}\phi$. We seek for an expression for $Z$. Then $Z$ satisfies the following equation:
\begin{align*}
P_g Z -\F{N+4}{2} \tilde u_0^\F{8}{N-4} Z =&\ \F{N+4}{2} \partial_{\xi'}(\tilde u_0^\F{8}{N-4}) \phi +\sum_{i=0}^N d_i \chi Z_i \\
&\ +\sum_{i=0}^N c_i \partial_{\xi'}(\chi Z_i) \qquad \text{in }\mathbb R^N,
\end{align*}
where $d_i=\partial_{\xi'}c_i$. Besides, from differentiating the orthogonality condition $\int_{\mathbb R^N} \phi \chi Z_j=0$,
we further get
\begin{align*}
\int_{\mathbb R^N} \phi~ \partial_{\xi'}(\chi Z_j) +\int_{\mathbb R^N} Z \chi Z_j=0.
\end{align*}\par

Choose $b_\ell$ such that
\begin{align*}
\sum_{\ell} b_\ell \int_{\mathbb R^N} \chi Z_\ell  Z_j =\int_{\mathbb R^N} \phi~ \partial_{\xi'}(\chi Z_j).
\end{align*}
Since this system is diagonal dominant with uniformly bounded coefficients, we see that it is uniquely solvable and that
\begin{align*}
|b_\ell| \leq C\mu\varepsilon^{10}\|\phi\|_*
\end{align*}
uniformly in $(\xi',\lambda')\in\Lambda$.\par

Let us now set
\begin{align}\label{3.1}
\eta=Z+\sum_{i=0}^N b_i \chi Z_i.
\end{align}
Then $\eta$ satisfies
\begin{equation*}
\begin{cases}
\D P_g\eta-\F{N+4}{2}\tilde u_0^\F{8}{N-4}\eta = \tilde\zeta + \sum_{i=0}^N d_i \chi Z_i \qquad \text{in }\mathbb R^N,\\
\D \int_{\mathbb R^N} \eta \chi Z_j=0\qquad \forall~j.
\end{cases}
\end{equation*}
where
\begin{align*}
\tilde \zeta = \sum_{i} b_i & \left(P_g (\chi Z_i)-\F{N+4}{2} \tilde u_0^\F{8}{N-4}(\chi Z_i)\right)+\F{N+4}{2} \partial_{\xi'}(\tilde u_0^\F{8}{N-4}) \phi \\
&\ +\sum_{i=0}^N c_i \partial_{\xi'}(\chi Z_i).
\end{align*}
Applying Lemma \ref{l3.2},
\begin{align*}
\|\eta\|_{*}\leq C\|\tilde \zeta\|_{**}.
\end{align*}
It can be directly checked that
\begin{align*}
|b_i|\left\|P_g(\chi Z_i)-\F{N+4}{2} \tilde u_0^\F{8}{N-4}(\chi Z_i) \right\|_{**} \leq  C\|\phi\|_*,
\end{align*}
\begin{align*}
&\ \left\|\F{N+4}{2} \partial_{\xi'}(\tilde u_0^\F{8}{N-4}) \phi\right\|_{**} \leq C\|\phi\|_*,\\
\intertext{and}
&\ \left\|c_i \partial_{\xi'}(\chi Z_i)\right\|_{**} \leq C\|\zeta\|_{**}.
\end{align*}
Thus $\|\tilde \zeta\|_{**}\leq C\|\zeta\|_{**}$, and then
\begin{align*}
\|\eta\|_{*}\leq C\|\zeta\|_{**}.
\end{align*}\par

Obviously $\|b_i Z_i\|_* \leq C\|\phi\|_* \leq C\|\zeta\|_{**}$. Therefore, we get by (\ref{3.1}) that
\begin{align*}
\|Z\|_* \leq C\|\zeta\|_{**}.
\end{align*}
The corresponding result for differentiation with respect to $\lambda'$ follows similarly. This concludes the proof.
\end{proof}

\section{Nonlinear Problem}

We recall that our aim is to solve Problem (\ref{2.1}). Rather than doing so directly,
we shall solve first the intermediate problem
\begin{equation}\label{4.1}
\begin{cases}
\D P_{\tilde g}\phi-\F{N+4}{2}\tilde u_0^\F{8}{N-4}\phi = -R+N(\phi)+\sum_{i=0}^N c_i \chi Z_i  \qquad  \text{in }\mathbb R^N,\\
\D \int_{\mathbb R^N} \phi \chi Z_j =0 \qquad \forall~j=0,1,\cdots,N.
\end{cases}
\end{equation}

\BP\label{p4.1}
There exists a unique solution to (\ref{4.1}) such that
\begin{align*}
\|\phi\|_* \leq \beta
\end{align*}
where $\beta>0$ is a large number independent of $\alpha$ and $\varepsilon$.
\EP

\begin{proof}
In terms of the operator $T$ defined in Remark \ref{r3.1}, Problem (\ref{4.1}) becomes
\begin{align*}
\phi = T[-R+ N(\phi)] : = A(\phi).
\end{align*}\par

For a given large number $\beta>0$, let us set
\begin{align*}
\mathcal S=\left\{\phi\in \mathcal H \cap L^\infty_*(\mathbb R^N):\ \|\phi\|_* \leq \beta \right\}.
\end{align*}\par

From Proposition \ref{p3.1}, we get
\begin{align*}
\|A(\phi)\|_{*} \leq C(\|R\|_{**}+\|N(\phi)\|_{**}).
\end{align*}
According to Lemma \ref{p4.3}, direct computation shows
\begin{align}\label{5.2}
\|R\|_{**} \leq C.
\end{align}
Here $C$ is independent of $\alpha$ and $\varepsilon$.
By the mean value theorem, we also easily have
\begin{align}\label{5.3}
\|N(\phi)\|_{**}\leq C\varepsilon^{4-\sigma}\|\phi\|_*^2.
\end{align}
Thus $A(\phi)\in \mathcal S$.\par

Furthermore, it is easy to check that for any $\phi_1, \phi_2\in \mathcal S$,
\begin{align*}
\|N(\phi_1)-N(\phi_2)\|_{**} \leq C\varepsilon^{4-\sigma} \|\phi_1-\phi_2\|_{*}.
\end{align*}
So
\begin{align*}
\|A(\phi_1)-A(\phi_2)\|_{*}\leq C\|N(\phi_1)-N(\phi_2)\|_{**}\leq C\varepsilon^{4-\sigma} \|\phi_1-\phi_2\|_{*},
\end{align*}
which implies that $A$ is a contraction mapping with the norm $\|\cdot\|_*$ inside $\mathcal S$.
Therefore the contraction mapping theorem yields the proposition.
\end{proof}

Our purpose in the remains of this section is to analyze the differentiability properties of the function $\phi$ defined in Proposition \ref{p4.1}

\BP\label{p4.2}
The function $(\xi',\lambda')\mapsto \phi(\xi',\lambda')$ provided by Proposition \ref{p4.1} is of class $C^1$ for the norm $\|\cdot\|_*$.
Moreover
\begin{align*}
\|\nabla_{(\xi',\lambda')} \phi\|_* \leq C.
\end{align*}
\EP

\begin{proof}
First, we come to the differentiability of $\phi_{(\xi',\lambda')}$. Consider the following map $H$:
$\Lambda\times\mathcal H\cap L^\infty_*(\mathbb R^N)\times \mathbb R^{N+1}\longrightarrow L^\infty_{**}(\mathbb R^N)\times \mathbb R^{N+1}$
of class $C^1$:
\begin{align*}
H((\xi',\lambda'),\phi,\bm c)=
\begin{pmatrix}
\D P_{\tilde g}(\tilde u_0+\phi)-\F{N-4}{2}(\tilde u_0+\phi)^\F{N+4}{N-4}-\sum_{i=0}^N c_{i}\chi Z_{i}\\
\D \int_{\mathbb R^N}\chi Z_{0}\phi\\
\vdots\\
\D \int_{\mathbb R^N}\chi Z_{N}\phi
\end{pmatrix}.
\end{align*}
Problem (\ref{4.1}) is then equivalent to $H((\xi',\lambda'),\phi,\bm c)=0$. We know that given $(\xi',\lambda')\in \Lambda$,
there is a unique solution $\phi_{(\xi',\lambda')}$. We will prove that the linear operator
\begin{align*}
\F{\partial H((\xi',\lambda'),\phi,\bm c)}{\partial(\phi,
c)}\bigg|_{((\xi',\lambda'),\phi_{(\xi',\lambda')},\bm c_{(\xi',\lambda')})} :\ \mathcal
H\cap L^\infty_*(\mathbb R^N) &\times\mathbb R^{N+1}\\
&\longrightarrow L^\infty_{**}(\mathbb R^N)\times\mathbb R^{N+1}
\end{align*}
is invertible for small $\varepsilon$ and $\alpha$. Then the $C^1$ regularity $(\xi',\lambda')\mapsto \phi(\xi',\lambda')$ follows from
the implicit function theorem. Indeed, we have
\begin{align*}
\F{\partial H((\xi',\lambda'),\phi,\bm c)}{\partial(\phi, \bm c)} &\bigg|_{((\xi',\lambda'),\phi_{(\xi',\lambda')},\bm c_{(\xi',\lambda')})}[\varphi,\bm d]\\
&=
\begin{pmatrix}
\D P_{\tilde g}\varphi - \F{N+4}{2} (\tilde u_0+\phi_{(\xi',\lambda')})^\F{8}{N-4}\varphi-\sum_{i=0}^N d_{i}\chi Z_{i}\\
\D \int_{\mathbb R^N}\chi Z_{0}\varphi\\
\vdots\\
\D \int_{\mathbb R^N}\chi_N Z_{N}\varphi
\end{pmatrix}.
\end{align*}
Since $\|\phi_{(\xi',\lambda')}\|_*\leq C$ , the same proof as that of Proposition \ref{p3.1} shows
that this operator is invertible for $\varepsilon$ and $\alpha$ small.\par

Since now $\phi=T[-R+N(\phi)]$, we have
\begin{align*}
\partial_{\xi'}\phi=\partial_{\xi'}T[-R+N(\phi)] + T[-\partial_{\xi'}R +\partial_{\xi'}N(\phi)].
\end{align*}
This implies, by Proposition \ref{p3.1} and \ref{p3.2},
\begin{align*}
\|\partial_{\xi'}\phi\|_* \leq C(\|R\|_{**}+\|N(\phi)\|_{**} + \|\partial_{\xi'}R\|_{**} +\|\partial_{\xi'}N(\phi)\|_{**}).
\end{align*}
Direct calculation shows
\begin{align*}
\partial_{\xi'} N(\phi)=&\ \F{N+4}{2}( \partial_{\xi'}\tilde u_0+\partial_{\xi'}\phi) \Big[(\tilde u_0+\phi)^\F{8}{N-4}-\tilde u_0^\F{8}{N-4}\Big]\\
&\ -\F{4(N+4)}{N-4}\tilde u_0^{-\F{N-12}{N-4}}(\partial_{\xi'}\tilde u_0) \phi.
\end{align*}
Hence we can obtain
\begin{align}\label{5.4}
\|\partial_{\xi'} N(\phi)\|_{**} \leq C \varepsilon^{4-\sigma}\|\phi\|_* (1+\|\partial_{\xi'}\phi\|_*)\leq C \varepsilon^{4-\sigma}(1+\|\partial_{\xi'}\phi\|_*).
\end{align}
It is easy to check
\begin{align}\label{5.5}
\|\partial_{\xi'}R\|_{**}\leq C.
\end{align}
Since we already have (\ref{5.2}) and (\ref{5.3}),
we can conclude
\begin{align*}
\|\partial_{\xi'}\phi\|_* \leq C.
\end{align*}\par

The corresponding result for differentiation with respect to $\lambda'$ can be gotten similarly. The proof is concluded.
\end{proof}

\section{Variational reduction}

As we have said, after Problem (\ref{4.1}) has been solved, we find a solution to
Problem (\ref{main equation}) if $(\xi',\lambda')$ is such that
\begin{align}\label{5.1}
c_i(\xi',\lambda') = 0\qquad   \text{for all } i.
\end{align}
This problem is indeed variational: it is equivalent to finding critical points of a
function of $(\xi',\lambda')$. To see this, we define
\begin{align*}
\mathcal F_{\tilde g}(\xi',\lambda')=E_{\tilde g} [\tilde u_0(\xi',\lambda')+\phi(\xi',\lambda')]
\end{align*}
where $\phi$ is the solution given by Proposition \ref{p4.1}.

\BL\label{l5.1}
$(\xi',\lambda')$ satisfies System (\ref{5.1}) if and only if $(\xi',\lambda')$ is a critical point of $\mathcal F_{\tilde g}$.
\EL

\begin{proof}
Since
\begin{align*}
P_{\tilde g}(\tilde u_0+\phi)-\F{N-4}{2} (\tilde u_0+\phi)^\F{N+4}{N-4} = \sum_{i=0}^N c_i \chi Z_i,
\end{align*}
we have
\begin{align*}
\partial_{(\xi',\lambda')} \mathcal F_{\tilde g}(\xi',\lambda')= \sum_{i=0}^N c_i \int_{\mathbb R^N}\chi Z_i [\partial_{(\xi',\lambda')}\tilde u_0+\partial_{(\xi',\lambda')}\phi],
\end{align*}
from which the necessity follows. In what follows we assume $\partial_{(\xi',\lambda')} \mathcal F_{\tilde g}(\xi',\lambda')=0$.
Then
\begin{align*}
\sum_{i=0}^N c_i \int_{\mathbb R^N}\chi Z_i \Big[\partial_{(\xi',\lambda')}\tilde u_0+\partial_{(\xi',\lambda')}\phi\Big] =0.
\end{align*}
Using Proposition  \ref{p4.2}, we can directly check that
\begin{align*}
\partial_{\xi'_i}\tilde u_0+\partial_{\xi'_i}\phi = Z_i + o(1),\\
\partial_{\lambda'}\tilde u_0+\partial_{\lambda'}\phi = Z_0 + o(1).
\end{align*}
Thus the above system for $c_i$ is diagonal dominant, which gives $c_i=0$ for all $i=0,\dots, N$. This concludes the proof.
\end{proof}

\section{Energy Expansion}

In this section we obtain an expansion of $\mathcal F_{\tilde g}$.

We first need to acquire a more refined estimate for $\phi$. Let $w(y)$ satisfies
\begin{equation}\label{5.7}
\begin{cases}
\D P_{\tilde g} w -\F{N+4}{2}\tilde u_0^\F{8}{N-4} w = R_1+\sum_{i}c_i \chi Z_i \qquad & \text{in }\mathbb R^N,\\
\D \int_{\mathbb R^N} \chi Z_i w =0,
\end{cases}
\end{equation}
where
\begin{align*}
R_1(y)=& -\mu\varepsilon^{10}\hat\chi(y)\Big[2\overline H_{ij}(\partial_{ijss}\tilde u_0)+ 2(\partial_s\overline H_{ij})( \partial_{ijs}\tilde u_0)
+ (\partial_{ss}\overline H_{ij})(\partial_{ij}\tilde u_0)\\
&\qquad\qquad -\F{b_N}{2}(\partial_{ss}\overline H_{ij})(\partial_{ij}\tilde u_0)-\F{b_N}{2}(\partial_{jss}\overline H_{ij})(\partial_{i}\tilde u_0)\Big].
\end{align*}
Here $\hat \chi$ is the cut-off function such that $\hat\chi(y)=1$ for $|y|\leq \F{\rho}{\varepsilon}$ and  $\hat\chi(y)=0$ for $|y|\geq \F{2\rho}{\varepsilon}$.
We define
\begin{align*}
\|w\|_*' &= \sup_{\mathbb R^N} \sum_{i=0}^2 \F{(1+|y-\xi'|)^{N-14+i}}{\mu\varepsilon^{10}}|\partial^i w(y)|,\\
\|R_1\|_{**}' &= \sup_{\mathbb R^N} \F{(1+|y-\xi'|)^{N-10}}{\mu\varepsilon^{10}}|R_1(y)|.
\end{align*}
A similar proof as that of Proposition \ref{p3.1} shows that there exists a unique solution $w$ to (\ref{5.7}) such that
\begin{align}\label{c}
\|w\|_*' \leq C\|R_1\|_{**}'\leq C.
\end{align}

We introduce again that
\begin{align*}
&\|\varphi_1\|_*'' =\sup_{\mathbb R^N}\sum_{i=0}^2\left\{ \F{1}{\F{\mu^2\varepsilon^{20}}{1+(|y-\xi'|)^{N-24+i}}+\alpha(\F{\varepsilon}{\rho})^{N-4+i}}+\F{(1+|y-\xi'|)^{N-4+i}}{\alpha}\right\}|\partial^i \varphi_1(y)|,\\
&\|\varphi_2\|_{**}'' = \sup_{\mathbb R^N}\Bigg\{ \F{\chi_{\{|y|\leq \F{\rho}{\varepsilon}\}}(1+|y-\xi'|)^{N-20}}{\mu^2\varepsilon^{20}} + \F{\chi_{\{\F{\rho}{\varepsilon}\leq |y|\leq \F{1}{\varepsilon}\}}(1+|y-\xi'|)^{N-1}}{\alpha\varepsilon}\\
&\qquad\qquad \qquad +\F{\chi_{\{|y|\geq \F{1}{\varepsilon}\}}(1+|y-\xi'|)^{N+\sigma}}{\alpha}\Bigg\}|\varphi_2(y)|.
\end{align*}

\BL
The function $\phi-w$ satisfies the estimate
\begin{align*}
\|\phi-w\|_*''\leq C.
\end{align*}
\EL

\begin{proof}
Obviously,
\begin{equation}\label{5.8}
\begin{cases}
\D P_{\tilde g} (\phi-w) -\F{N+4}{2}\tilde u_0^\F{8}{N-4} (\phi-w) = -R_2+N(\phi)+\sum_{i=0}^N c_i \chi Z_i & \text{in }\mathbb R^N,\\
\D \int_{\mathbb R^N} \chi Z_i (\phi-w) =0,
\end{cases}
\end{equation}
where
\begin{align*}
R_2(y)=& \Delta^2_{\tilde g} \tilde u_0 -\Delta^2 \tilde u_0-\mu\varepsilon^{10}\hat\chi\sum_{i,j,s}\Big[2\overline H_{ij}(\partial_{ijss}\tilde u_0)+ 2(\partial_s\overline H_{ij})( \partial_{ijs}\tilde u_0)\\
&\qquad+ (\partial_{ss}\overline H_{ij})(\partial_{ij}\tilde u_0) -\F{b_N}{2}(\partial_{ss}\overline H_{ij})(\partial_{ij}\tilde u_0)-\F{b_N}{2}(\partial_{jss}\overline H_{ij})(\partial_{i}\tilde u_0)\Big].
\end{align*}
It is easy to check
\begin{align*}
|R_2(y)|\leq
\begin{cases}
\D C \F{\mu^2 \varepsilon^{20}}{(1+|y-\xi'|)^{N-20}} \qquad & \D \text{for }|y|\leq \F{\rho}{\varepsilon},\medskip\\
\D C\F{\alpha\varepsilon}{(1+|y-\xi'|)^{N-1}} & \D \text{for }\F{\rho}{\varepsilon}\leq |y|\leq \F{1}{\varepsilon},\medskip \\
0& \D \text{for }|y|\geq \F{1}{\varepsilon}.
\end{cases}
\end{align*}
On the other hand, since $|N(\phi)|\leq C (1+|y-\xi'|)^{N-12}|\phi|^2$ and $\|\phi\|_* \leq C$, we have
\begin{align*}
|N(\phi)| \leq
\begin{cases}
\D C \F{\mu^2 \varepsilon^{20}}{(1+|y-\xi'|)^{N-16}}+\alpha(\F{\varepsilon}{\rho})^{N+4} \qquad & \D \text{for }|y|\leq \F{\rho}{\varepsilon},\medskip\\
\D C\F{\alpha}{(1+|y-\xi'|)^{N+4}} & \D \text{for }|y|\geq \F{\rho}{\varepsilon}.
\end{cases}
\end{align*}
Similar to the proof of Lemma \ref{l3.2}, we obtain
\begin{equation*}
\|\phi-w\|_*'' \leq C\|R_2\|_{**}'' + C\|N(\phi)\|_{**}'' \leq C.\qedhere
\end{equation*}
\end{proof}

\BP\label{l5.2}
The following expansion holds
\begin{align*}
2\mathcal F_{\tilde g}(\xi',\lambda')= &\ 2E-\int_{B_\F{\rho}{\varepsilon}}\sum_{i,j,k,\ell}\tilde h_{i\ell} \tilde h_{j\ell}(\partial_{ikk}\tilde u_0)(\partial_{j}\tilde u_0)+\int_{B_\F{\rho}{\varepsilon}}\Big(\sum_{i,j} \tilde h_{ij}(\partial_{ij}\tilde u_0)\Big)^2\\
&\ -\F{a_N}{4}\int_{B_\F{\rho}{\varepsilon}} \sum_{i,k,\ell,m}(\partial_\ell \tilde h_{mk})^2(\partial_i\tilde u_0)^2\\
&\ -\F{b_N}{4}\int_{B_\F{\rho}{\varepsilon}}\sum_{i,j,m,s}(\partial_j\tilde h_{ms})(\partial_i\tilde h_{sm})(\partial_i\tilde u_0)(\partial_j\tilde u_0)\\
&\ -\F{b_N}{2}\int_{B_\F{\rho}{\varepsilon}} \sum_{i,j,m,s}\Big[\tilde h_{ms}(\partial_{s}\tilde h_{ij})-\tilde h_{si}(\partial_s\tilde h_{mj}) +\tilde h_{sj}(\partial_i\tilde h_{ms})\\
&\hspace{15em} -\tilde h_{ms}(\partial_{i}\tilde h_{sj})\Big]\partial_m (\partial_i\tilde u_0\partial_j\tilde u_0)\\
&\ +\F{b_N}{2}\int_{B_\F{\rho}{\varepsilon}} \sum_{i,j,m,s}\tilde h_{is}(\partial_{mm}\tilde h_{js})(\partial_i\tilde u_0)(\partial_j\tilde u_0)\\
&\ +\F{N-4}{8(N-1)}\int_{B_\F{\rho}{\varepsilon}} \sum_{i,k,\ell,m}\Big[(\partial_{i\ell}h_{mk})^2+(\partial_\ell h_{mk})(\partial_{ii\ell}h_{mk})\Big]\tilde u_0^2\\
&\ -\F{N-4}{4(N-2)^2}\int_{B_\F{\rho}{\varepsilon}} \sum_{i,j,m,s}(\partial_{mm}\tilde h_{ij})(\partial_{ss}\tilde h_{ij})u_0^2\\
&\ +\int_{B_\F{\rho}{\varepsilon}} \sum_{i,j,s}\Big[2\tilde h_{ij}(\partial_{ijss}\tilde u_0)+ 2(\partial_s\tilde h_{ij})( \partial_{ijs}\tilde u_0)+ (\partial_{ss}\tilde h_{ij})(\partial_{ij}\tilde u_0)\nonumber\\
&\qquad\qquad -\F{b_N}{2}(\partial_{ss}\tilde h_{ij})(\partial_{ij}\tilde u_0)-\F{b_N}{2}(\partial_{jss}\tilde h_{ij})(\partial_{i}\tilde u_0)\Big]w\\
& +O\left(\mu^3 \varepsilon^\F{20N}{N-1}\right)+O\left(\mu^\F{2N}{N-4} \varepsilon^\F{20N}{N-4}\right)+O\left(\alpha\left(\F{\varepsilon}{\rho}\right)^{N-4}\right).\\
\end{align*}
\EP

\begin{proof}
Since $\phi$ is a solution to (\ref{4.1}),
\begin{align}\label{6.7}
&\ \int_{\mathbb R^N}\Delta_{\tilde g}(\tilde u_0+\phi)\Delta_{\tilde g}\phi+\sum_{i,j}(a_N S_{\tilde g}g^{ij}+b_N \widetilde{\mathcal R}^{ij})\partial_i(\tilde u_0+\phi) \partial_j\phi\nonumber\\
&\qquad\qquad\qquad +\F{N-4}{2}Q_{\tilde g}(\tilde u_0+\phi) \phi \nonumber\\
=&\ \F{N-4}{2}\int_{\mathbb R^N} (\tilde u_0+\phi)^\F{N+4}{N-4} \phi.
\end{align}
On the other hand
\begin{align}\label{6.8}
&\ \int_{\mathbb R^N}\Delta_{\tilde g}\tilde u_0\Delta_{\tilde g}\phi+\sum_{i,j}(a_N S_{\tilde g}g^{ij}+b_N \widetilde{\mathcal R}^{ij})\partial_i\tilde u_0 \partial_j\phi+\F{N-4}{2}Q_{\tilde g}\tilde u_0\phi \nonumber\\
=&\ \int_{\mathbb R^N} \Big\{P_{\tilde g}\tilde u_0-\Delta^2 \tilde u_0-\sum_{i,j,s}\Big[2h_{ij}(\partial_{ijss}\tilde u_0)+ 2(\partial_s\tilde h_{ij})( \partial_{ijs}\tilde u_0)+ (\partial_{ss} \tilde h_{ij})(\partial_{ij}\tilde u_0)\nonumber\\
&\qquad\qquad -\F{b_N}{2}(\partial_{ss} \tilde h_{ij})(\partial_{ij}\tilde u_0)-\F{b_N}{2}(\partial_{jss} \tilde h_{ij})(\partial_{i}\tilde u_0)\Big]\Big\}\phi \nonumber\\
&\ +\int_{\mathbb R^N} \sum_{i,j,s}\Big[2\tilde h_{ij}(\partial_{ijss}\tilde u_0)+ 2(\partial_s \tilde h_{ij})( \partial_{ijs}\tilde u_0)+ (\partial_{ss} \tilde h_{ij})(\partial_{ij}\tilde u_0)\nonumber\\
&\qquad\qquad -\F{b_N}{2}(\partial_{ss} \tilde h_{ij})(\partial_{ij}\tilde u_0)-\F{b_N}{2}(\partial_{jss} \tilde h_{ij})(\partial_{i}\tilde u_0)\Big](\phi-w)\nonumber\\
&\ +\int_{\mathbb R^N} \sum_{i,j,s}\Big[2\tilde h_{ij}(\partial_{ijss}\tilde u_0)+ 2(\partial_s \tilde h_{ij})( \partial_{ijs}\tilde u_0)+ (\partial_{ss} \tilde h_{ij})(\partial_{ij}\tilde u_0)\nonumber\\
&\qquad\qquad -\F{b_N}{2}(\partial_{ss} \tilde h_{ij})(\partial_{ij}\tilde u_0)-\F{b_N}{2}(\partial_{jss} \tilde h_{ij})(\partial_{i}\tilde u_0)\Big]w\nonumber\\
&\ +\F{N-4}{2}\int_{\mathbb R^N} \tilde u_0^\F{N+4}{N-4} \phi \nonumber\\
:=&\ I_1+I_2+I_3 +\F{N-4}{2}\int_{\mathbb R^N} \tilde u_0^\F{N+4}{N-4}\phi.
\end{align}
Since
\begin{align*}
&\ \bigg|P_{\tilde g}\tilde u_0-\Delta^2 \tilde u_0-\sum_{i,j,s}\Big[2\tilde h_{ij}(\partial_{ijss}\tilde u_0)+ 2(\partial_s\tilde h_{ij})( \partial_{ijs}\tilde u_0)+ (\partial_{ss}\tilde h_{ij})(\partial_{ij}\tilde u_0)\nonumber\\
&\qquad\qquad -\F{b_N}{2}(\partial_{ss}\tilde h_{ij})(\partial_{ij}\tilde u_0)-\F{b_N}{2}(\partial_{jss}\tilde h_{ij})(\partial_{i}\tilde u_0)\Big]\bigg|\\
\leq &\
\begin{cases}
\D C\F{\mu^2\varepsilon^{20}}{(1+|y-\xi'|)^{N-20}}\qquad  &\D \text{for }|y|\leq \F{\rho}{\varepsilon},\medskip\\
\D C\F{\alpha\varepsilon}{(1+|y-\xi'|)^{N-1}} &\D \text{for }|y|\geq \F{\rho}{\varepsilon},
\end{cases}
\end{align*}
we have
\begin{align}\label{6.4}
|I_1| \leq C\mu^3\varepsilon^{21}\rho^9|\log\varepsilon|+C\alpha^2\rho(\F{\varepsilon}{\rho})^{N-4}.
\end{align}
Since
\begin{align*}
&\ \bigg|\sum_{i,j,s}\Big[2\tilde h_{ij}(\partial_{ijss}\tilde u_0)+ 2(\partial_s\tilde h_{ij})( \partial_{ijs}\tilde u_0)+ (\partial_{ss}\tilde h_{ij})(\partial_{ij}\tilde u_0)\nonumber\\
&\qquad\qquad -\F{b_N}{2}(\partial_{ss}\tilde h_{ij})(\partial_{ij}\tilde u_0)-\F{b_N}{2}(\partial_{jss}\tilde h_{ij})(\partial_{i}\tilde u_0)\Big]\bigg|\\
\leq &\
\begin{cases}
\D C\F{\mu\varepsilon^{10}}{(1+|y-\xi'|)^{N-10}}\qquad  &\D \text{for }|y|\leq \F{\rho}{\varepsilon},\medskip\\
\D C\F{\alpha\varepsilon}{(1+|y-\xi'|)^{N-1}} &\D \text{for }|y|\geq \F{\rho}{\varepsilon}
\end{cases}
\end{align*}
and $\|\phi-w\|_*'' \leq C$, we obtain that
\begin{align}\label{6.5}
|I_2|\leq C\mu^3\varepsilon^{21}\rho^9|\log\varepsilon|+C\alpha^2\rho(\F{\varepsilon}{\rho})^{N-4}.
\end{align}
Since $\|w\|_*'\leq C$ and $h=0$ for $|y|\geq \F{1}{\varepsilon}$.
\begin{align}\label{6.6}
I_3=&\ \int_{|y|\leq \F{\rho}{\varepsilon}} \sum_{i,j,s}\Big[2\tilde h_{ij}(\partial_{ijss}\tilde u_0)+ 2(\partial_s\tilde h_{ij})( \partial_{ijs}\tilde u_0)+ (\partial_{ss}\tilde h_{ij})(\partial_{ij}\tilde u_0)\nonumber\\
&\qquad\qquad -\F{b_N}{2}(\partial_{ss}\tilde h_{ij})(\partial_{ij}\tilde u_0)-\F{b_N}{2}(\partial_{jss}\tilde h_{ij})(\partial_{i}\tilde u_0)\Big]w\nonumber\\
&\ + O\left(\alpha\mu\rho^{10} (\F{\varepsilon}{\rho})^{N-4}\right).
\end{align}
Thus, combining (\ref{6.7})-(\ref{6.6}), we have
\begin{align}
&\ 2\mathcal F_{\tilde g}(\xi',\lambda') -2E_{\tilde g}(\tilde u_0)\nonumber\\
=&\ \int_{|y|\leq \F{\rho}{\varepsilon}} \sum_{i,j,s}\Big[2\tilde h_{ij}(\partial_{ijss}\tilde u_0)+ 2(\partial_s\tilde h_{ij})( \partial_{ijs}\tilde u_0)+ (\partial_{ss}\tilde h_{ij})(\partial_{ij}\tilde u_0)\nonumber\\
&\qquad\qquad -\F{b_N}{2}(\partial_{ss}\tilde h_{ij})(\partial_{ij}\tilde u_0)-\F{b_N}{2}(\partial_{jss}\tilde h_{ij})(\partial_{i}\tilde u_0)\Big]w\nonumber\\
&\ +\F{N-4}{2}\int_{\mathbb R^N}\Bigg\{\F 4 N \left[(\tilde u_0+\phi)^\F{2N}{N-4}-\tilde u_0^\F{2N}{N-4}\right]\nonumber\\
&\qquad\qquad\qquad\qquad -\left[(\tilde u_0+\phi)^\F{N+4}{N-4}\tilde u_0-\tilde u_0^\F{N+4}{N-4}(\tilde u_0+\phi)\right]\Bigg\}\nonumber\\
&\ +O\left(\mu^3\varepsilon^{21}\rho^9|\log\varepsilon|\right)+O\left(\alpha^2\rho(\F{\varepsilon}{\rho})^{N-4}\right).
\end{align}
Using the fact that $\|\phi\|_*\leq C$, we have the pointwise estimate
\begin{align*}
&\ \Bigg|\F 4 N \left[(\tilde u_0+\phi)^\F{2N}{N-4}-\tilde u_0^\F{2N}{N-4}\right]-\left[(\tilde u_0+\phi)^\F{N+4}{N-4}\tilde u_0-\tilde u_0^\F{N+4}{N-4}(\tilde u_0+\phi)\right]\Bigg|
\leq C|\phi|^\F{2N}{N-4},
\end{align*}
which implies
\begin{align*}
&\ \int_{\mathbb R^N}\Bigg\{\F 4 N \left[(\tilde u_0+\phi)^\F{2N}{N-4}-\tilde u_0^\F{2N}{N-4}\right]
-\left[(\tilde u_0+\phi)^\F{N+4}{N-4}\tilde u_0-\tilde u_0^\F{N+4}{N-4}(\tilde u_0+\phi)\right]\Bigg\}\\
= &\ O(\mu^\F{2N}{N-4}\varepsilon^\F{20N}{N-4})+O\left(\alpha(\F{\varepsilon}{\rho})^N\right).
\end{align*}
Proposition \ref{p6.1} then gives the result.
\end{proof}

\section{Finding a Critical Point for Reduced Energy Functional}

We define a function $F$: $\mathbb R^N\times(0,\infty)\to\mathbb R$ as follows: given an pair $(\xi',\lambda')\in \mathbb R^N\times(0,\infty)$,
\begin{align*}
F(\xi',\lambda')=&\ -\int_{\mathbb R^N}\sum_{i,j,k,\ell}\overline H_{i\ell} \overline H_{j\ell}(\partial_{ikk}\tilde u_0)(\partial_{j}\tilde u_0)+\int_{\mathbb R^N}\left(\sum_{i,j} \overline H_{ij}(\partial_{ij}\tilde u_0)\right)^2\\
&\ -\F{a_N}{4}\int_{\mathbb R^N} \sum_{i,k,\ell,m}(\partial_\ell \overline H_{mk})^2(\partial_i\tilde u_0)^2\\
&\ -\F{b_N}{4}\int_{\mathbb R^N}\sum_{i,j,m,s}(\partial_j\overline H_{ms})(\partial_i\overline H_{sm})(\partial_i\tilde u_0)(\partial_j\tilde u_0)\\
&\ -\F{b_N}{2}\int_{\mathbb R^N} \sum_{i,j,m,s}\Big[\overline H_{ms}(\partial_{s}\overline H_{ij})-\overline H_{si}(\partial_s\overline H_{mj}) +\overline H_{sj}(\partial_i\overline H_{ms})\\
&\hspace{15em} -\overline H_{ms}(\partial_{i}\overline H_{sj})\Big]\partial_m (\partial_i\tilde u_0\partial_j\tilde u_0)\\
&\ +\F{b_N}{2}\int_{\mathbb R^N} \sum_{i,j,m,s}\overline H_{is}(\partial_{mm}\overline H_{js})(\partial_i\tilde u_0)(\partial_j\tilde u_0)\\
&\ +\F{N-4}{8(N-1)}\int_{\mathbb R^N} \sum_{i,k,\ell,m}\Big[(\partial_{i\ell}\overline H_{mk})^2+(\partial_\ell \overline H_{mk})(\partial_{ii\ell}\overline H_{mk})\Big]\tilde u_0^2\\
&\ -\F{N-4}{4(N-2)^2}\int_{\mathbb R^N} \sum_{i,j,m,s}(\partial_{mm}\overline H_{ij})(\partial_{ss}\overline H_{ij})u_0^2\\
&\ +\int_{\mathbb R^N} \sum_{i,j,s}\Big[2\overline H_{ij}(\partial_{ijss}\tilde u_0)+ 2(\partial_s\overline H_{ij})( \partial_{ijs}\tilde u_0)+ (\partial_{ss}\overline H_{ij})(\partial_{ij}\tilde u_0)\nonumber\\
&\qquad\qquad -\F{b_N}{2}(\partial_{ss}\overline H_{ij})(\partial_{ij}\tilde u_0)-\F{b_N}{2}(\partial_{jss}\overline H_{ij})(\partial_{i}\tilde u_0)\Big]\bar w
\end{align*}
where $\bar w(y) = w(y)/\mu\varepsilon^{10}$ and $w(y)$ is defined in (\ref{5.7}).

Proposition  \ref{p6.1} shows that the reduced energy functional ${\mathcal F}_{\tilde{g}}$ is close to $F$.

Our goal in this section is to show that the function $F(\xi',\lambda')$ has a strict (nondegenerate) local minimum point. First we have the following symmetry result.

\BP\label{p8.1}
The function $F(\xi',\lambda')$ satisfies $F(\xi',\lambda')=F(-\xi',\lambda')$ for all $(\xi',\lambda')\in \mathbb R^N\times(0,\infty)$.
Consequently, we have $\F{\partial}{\partial{\xi'}}F(0,\lambda')=0$ and $\F{\partial^2}{\partial\xi'\partial\lambda'}F(0,\lambda')=0$
for all $\lambda'>0$.
\EP

\begin{proof}
This follows immediately from the relation $\overline H_{ik}(-y)=\overline H_{ik}(y)$.
\end{proof}

To find a local minimum for $F$, we list some useful identities, which are all direct consequences of  definitions.

\BL
\begin{align*}
\partial_i \tilde u_0 &= -(N-4) \tilde u_0(y)\F{y_i-\xi'_i}{\lambda'^2+|y-\xi'|^2},\\
\partial_{ij} \tilde u_0 &=(N-2)(N-4)\tilde u_0(y)\F{(y_i-\xi'_i)(y_j-\xi'_j)}{(\lambda'^2+|y-\xi'|^2)^2}-(N-4)\tilde u_0(y)\F{\delta_{ij}}{\lambda'^2+|y-\xi'|^2},\\
\partial_{ijs} \tilde u_0 &= -N(N-2)(N-4)\tilde u_0(y)\F{(y_i-\xi'_i)(y_j-\xi'_j)(y_s-\xi'_s)}{(\lambda'^2+|y-\xi'|^2)^3}\\
&\quad +(N-2)(N-4)\tilde u_0(y)\F{\delta_{js}(y_i-\xi'_i)+\delta_{is}(y_j-\xi'_j)+\delta_{ij}(y_s-\xi'_s)}{(\lambda'^2+|y-\xi'|^2)^2}.
\end{align*}
\EL

\BL
\begin{align*}
\partial_\ell \overline H_{mk}&=2f'(|y|^2) H_{mk}y_\ell + f(|y|^2)(\partial_\ell H_{mk}),\\
\partial_{i\ell}\overline{H}_{mk}&= 2\delta_{i\ell}f'(|y|^2) H_{mk}+2f'(|y|^2)(\partial_i H_{mk})y_\ell+2f'(|y|^2)(\partial_\ell H_{mk})y_i\\
&\qquad +f(|y|^2)(\partial_{i\ell}H_{mk})+4f''(|y|^2)H_{mk} y_i y_\ell.
\end{align*}
\EL

\BL
\begin{align*}
\sum_i y_i\partial_i H_{ms}&=2H_{ms},& \sum_j y_j\partial_s H_{tj} &= -H_{st},\\
\sum_j y_j \partial_{ms}H_{ij}& =-\partial_s H_{im}-\partial_m H_{is},& \sum_\ell y_\ell (\partial_{i\ell}H_{mk})&=\partial_i H_{mk}.
\end{align*}
\EL

\BL\label{l9.1}
\begin{align*}
\int_{\partial B_1} \sum_{i,j} H_{ij}^2 &= \F{|S^{N-1}|}{2N(N+2)}\sum_{i,j,k,\ell}(W_{ikj\ell}+W_{i\ell j k})^2,\\
\int_{\partial B_1} \sum_{i,j,k} (\partial_k H_{ij})^2 &= \F{|S^{N-1}|}{N}\sum_{i,j,k,\ell}(W_{ikj\ell}+W_{i\ell j k})^2,\\
\int_{\partial B_1} \sum_{i,j,k,\ell} (\partial_{k\ell} H_{ij})^2 &= |S^{N-1}|\sum_{i,j,k,\ell}(W_{ikj\ell}+W_{i\ell j k})^2.
\end{align*}
\EL

\BL\label{l9.2}
\begin{align*}
\int_{\partial B_1} \sum_{i,j} H_{ij}^2y_py_q =&\ \F{2|S^{N-1}|}{N(N+2)(N+4)}\sum_{i,j,k}(W_{ikjp}+W_{ip j k})(W_{ikjq}+W_{iq j k})\\
&\qquad +\F{|S^{N-1}|}{2N(N+2)(N+4)}(W_{ikj\ell}+W_{i\ell j k})^2\delta_{pq},\\
\int_{\partial B_1} \sum_{t} H_{pt}H_{qt} =&\ \F{|S^{N-1}|}{2N(N+2)}\sum_{i,j,k}(W_{ikjp}+W_{ip j k})(W_{ikjq}+W_{iq j k}),\\
\int_{\partial B_1} \sum_{i,j,k} (\partial_{k} H_{ij})^2 y_p y_q =&\ \F{2|S^{N-1}|}{N(N+2)}\sum_{i,j,k}(W_{ikjp}+W_{ip j k})(W_{ikjq}+W_{iq j k})\\
&\qquad +\F{|S^{N-1}|}{N(N+2)}(W_{ikj\ell}+W_{i\ell j k})^2\delta_{pq},\\
\int_{\partial B_1} \sum_{i,j} (\partial_p H_{ij})(\partial_q H_{ij}) =&\ \F{|S^{N-1}|}{N}\sum_{i,j,k}(W_{ikjp}+W_{ip j k})(W_{ikjq}+W_{iq j k}),\\
\int_{\partial B_1} \sum_{i,j} H_{ij}(\partial_q H_{ij})y_p =&\ \F{|S^{N-1}|}{N(N+2)}\sum_{i,j,k}(W_{ikjp}+W_{ip j k})(W_{ikjq}+W_{iq j k}),\\
\int_{\partial B_1} \sum_{i,j,k,\ell} (\partial_{k\ell }H_{ij})^2 y_p y_q =&\ \F{|S^{N-1}|}{N}\sum_{i,j,k,\ell}(W_{ikj\ell}+W_{i\ell j k})^2\delta_{pq}.
\end{align*}
\EL

We mention that the proof of some identities in Lemma \ref{l9.1} and \ref{l9.2} can be found in \cite{B,BM}. The others may be proved by the same method.

\BL
It holds
\begin{align*}
\int_{\mathbb R^N}\left.\left\{\sum_{i,j,s,t}\overline H_{tj}\overline H_{it}(\partial_{iss}\tilde u_0)(\partial_j\tilde u_0)\right\}\right|_{\xi'=0}=0,
\end{align*}
and
\begin{align*}
&\ \int_{\mathbb R^N}\F{\partial^2}{\partial{\xi'_p}\partial{\xi'_q}}\left.\left\{\sum_{i,j,s,t}\overline H_{tj}\overline H_{it}(\partial_{iss}\tilde u_0)(\partial_j\tilde u_0)\right\}\right|_{\xi=0}\\
=&\ |S^{N-1}|\sum_{i,j,k}(W_{ikjp}+W_{ipjk})(W_{ikjq}+W_{iqjk})\lambda'^{N-4}\\
&\quad \cdot \F{(N-2)(N-4)^2}{N(N+2)}\int_0^\infty f(r^2)^2 \left[\F{Nr^{N+5}}{(\lambda'^2+r^2)^N}-\F{(N+2)r^{N+3}}{(\lambda'^2+r^2)^{N-1}}\right]\mathrm{d}r.
\end{align*}
\EL

\begin{proof}
Since
\begin{align*}
&\ \sum_{i,j,s,t}\overline H_{tj}\overline H_{it}(\partial_{iss}\tilde u_0)(\partial_j\tilde u_0)\\
=&\ f(|y|^2)^2 H_{tj}H_{it}(N-2)(N-4)^2 \tilde u_0^2\left[\F{N|y-\xi'|^2}{(\lambda'^2+|y-\xi'|^2)^4}-\F{(n+2)}{(\lambda'^2+|y-\xi'|^2)^3}\right]\\
&\qquad \cdot (y_i-\xi'_i)(y_j-\xi'_j)\\
=&\ N(N-2)(N-4)^2 f(|y|^2)^2 \F{\lambda'^{N-4}|y-\xi'|^2}{(\lambda'^2+|y-\xi'|^2)^{N}}\sum_{i,j,t}H_{it}H_{jt}\xi'_i\xi'_j\\
&\ -(N-2)(N+2)(N-4)^2 f(|y|^2)^2 \F{\lambda'^{N-4}}{(\lambda'^2+|y'-\xi'|^2)^{N-1}}\sum_{i,j,t}H_{it}H_{jt}\xi'_i\xi'_j
\end{align*}
and
\begin{align*}
&\ \F{\partial^2}{\partial{\xi'_p}\partial{\xi'_q}}\left.\left\{\sum_{i,j,s,t}\overline H_{tj}\overline H_{it}(\partial_{iss}\tilde u_0)(\partial_j \tilde u_0)\right\}\right|_{\xi'=0}\\
=&\ 2N(N-2)(N-4)^2 f(|y|^2)^2 \F{\lambda'^{N-4}|y|^2}{(\lambda'^2+|y|^2)^{N}}\sum_{t}H_{pt}H_{qt}\\
&\ -2(N-2)(N+2)(N-4)^2 f(|y|^2)^2 \F{\lambda'^{N-4}}{(\lambda'^2+|y|^2)^{N-1}}\sum_{t}H_{pt}H_{qt},
\end{align*}
the lemma follows directly by Lemma \ref{l9.1} and \ref{l9.2}.
\end{proof}

\BL
We have
\begin{align*}
\int_{\mathbb R^N}\left.\left\{\sum_{i,j,s,t}\overline{H}_{ij}\overline{H}_{st}\partial_{ij}\tilde u_0\partial_{st}\tilde u_0\right\}\right|_{\xi'=0}=0,\\
\int_{\mathbb R^n}\F{\partial^2}{\partial{\xi'_p}\partial{\xi'_q}}\left.\left\{\sum_{i,j,s,t}\overline{H}_{ij}\overline{H}_{st}\partial_{ij}\tilde u_0\partial_{st}\tilde u_0\right\}\right|_{\xi'=0}=0.
\end{align*}
\EL
\begin{proof}
The lemma easily follows from
\begin{align*}
&\ \sum_{i,j,s,t}\overline{H}_{ij}\overline{H}_{st}\partial_{ij}\tilde u_0\partial_{st}\tilde u_0\\
=&\ (N-2)^2(N-4)^2 f(|y|^2)^2 \F{\lambda'^{N-4}}{(\lambda'^2+|y-\xi'|^2)^{N}}\sum_{i,j,s,t}H_{ij}H_{st}\xi'_i\xi'_j\xi'_s\xi'_t.\tag*{\qedhere}
\end{align*}
\end{proof}

\BL
There hold
\begin{align*}
&\ \int_{\mathbb R^N}\left.\left\{\sum_{i,\ell,k,m}(\partial_\ell \overline H_{mk})^2 (\partial_i\tilde u_0)^2\right\}\right|_{\xi=0}\\
=&\ |S^{N-1}|\sum_{i,j,k,\ell}(W_{ikj\ell}+W_{i\ell jk})^2\lambda'^{N-4}\\
&\cdot \Bigg\{\F{2(N-4)^2}{N(N+2)}\int_0^\infty \left[r^2 f'(r^2)^2 +2 f(r^2)f'(r^2)\right]\F{r^{N+5}}{(\lambda'^2+r^2)^{N-2}}\mathrm{d}r\\
&\qquad\qquad + \F{(N-4)^2}{N}\int_0^\infty f(r^2)^2 \F{r^{N+3}}{(\lambda'^2+r^2)^{N-2}}\mathrm{d}r\Bigg\},
\end{align*}
and
\begin{align*}
&\ \int_{\mathbb R^N}\left.\F{\partial^2}{\partial \xi'_p \partial \xi'_q}\left\{\sum_{i,\ell,k,m}(\partial_\ell \overline H_{mk})^2 (\partial_i\tilde u_0)^2\right\}\right|_{\xi'=0}\\
=&\ |S^{N-1}|\sum_{i,j,k}(W_{ikjp}+W_{ipjk})(W_{ikjq}+W_{iqjk})\lambda'^{N-4}\\
&\cdot\Bigg\{\F{32(N-2)(N-4)^2}{N(N+2)(N+4)}\int_0^\infty \left[r^2f'(r^2)^2+2f(r^2)f'(r^2)\right]\\
&\hspace{15em}\cdot\left[\F{(N-1)r^{N+7}}{(\lambda'^2+r^2)^N}-\F{2r^{N+5}}{(\lambda'^2+r^2)^{N-1}}\right]\mathrm{d}r\\
&\qquad +\F{8(N-2)(N-4)^2}{N(N+2)}\int_0^\infty f(r^2)^2 \left[\F{(N-1)r^{N+5}}{(\lambda'^2+r^2)^N}-\F{2r^{N+3}}{(\lambda'^2+r^2)^{N-1}}\right]\mathrm{d}r\Bigg\}\\
&\ +|S^{N-1}|\sum_{i,j,k,\ell}(W_{ikj\ell}+W_{i\ell jk})^2\lambda'^{N-4}\delta_{pq}\\
&\cdot\Bigg\{\F{4(N-4)^2}{N(N+2)(N+4)}\int_0^\infty \left[r^2f'(r^2)^2+2f(r^2)f'(r^2)\right]\\
&\qquad \cdot \left[\F{2(N-1)(N-2)r^{N+7}}{(\lambda'^2+r^2)^n}-\F{(N-2)(N+8)r^{N+5}}{(\lambda'^2+r^2)^{N-1}}+\F{(N+4)r^{N+3}}{(\lambda'^2+r^2)^{N-2}}\right]\mathrm{d}r\\
&\quad + \F{2(N-4)^2}{N(N+2)}\int_0^\infty f(r^2)^2\\
&\qquad \cdot \left[\F{2(N-1)(N-2)r^{N+5}}{(\lambda'^2+r^2)^N}-\F{(N-2)(N+6)r^{N+3}}{(\lambda'^2+r^2)^{N-1}}+\F{(N+2)r^{N+1}}{(\lambda'^2+r^2)^{N-2}}\right]\mathrm{d}r\Bigg\}.
\end{align*}
\EL

\begin{proof}
Direct computation shows
\begin{align*}
&\ \sum_{i,\ell,k,m}(\partial_\ell \overline H_{mk})^2 (\partial_i\tilde u_0)^2\\
=&\ 4(N-4)^2\Big[|y|^2f'(|y|^2)^2+2f(|y|^2)f'(|y|^2)\Big]\lambda'^{N-4}\F{|y-\xi'|^2}{(\lambda'^2+|y-\xi'|^2)^{N-2}}\sum_{m,k}H_{mk}^2\\
&\ +(N-4)^2f(|y|^2)^2\lambda'^{N-4}\F{|y-\xi'|^2}{(\lambda'^2+|y-\xi'|^2)^{N-2}}\sum_{k,\ell,m}(\partial_\ell H_{mk})^2,
\end{align*}
\begin{align*}
&\ \left.\F{\partial^2}{\partial \xi'_p \partial \xi'_q}\left\{\sum_{i,\ell,k,m}(\partial_\ell \overline H_{mk})^2 (\partial_i\tilde u_0)^2\right\}\right|_{\xi'=0}\\
=&\ 16(N-1)(N-2)(N-4)^2\Big[|y|^2f'(|y|^2)^2+2f(|y|^2)f'(|y|^2)\Big]\F{|y|^2\lambda'^{N-4}}{(\lambda'^2+|y|^2)^{N}}\sum_{k,m}H_{km}^2y_py_q\\
&\ -32(N-2)(N-4)^2\Big[|y|^2f'(|y|^2)^2+2f(|y|^2)f'(|y|^2)\Big]\F{\lambda'^{N-4}}{(\lambda'^2+|y|^2)^{N-1}}\sum_{k,m}H_{km}^2y_py_q\\
&\ -8(N-2)(N-4)^2\Big[|y|^2f'(|y|^2)^2+2f(|y|^2)f'(|y|^2)\Big]\F{|y|^2\lambda'^{N-4}}{(\lambda'^2+|y|^2)^{N-1}}\sum_{k,m}H_{km}^2\delta_{pq}\\
&\ +8(N-4)^2\Big[|y|^2f'(|y|^2)^2+2f(|y|^2)f'(|y|^2)\Big]\F{\lambda'^{N-4}}{(\lambda'^2+|y|^2)^{N-2}}\sum_{k,m}H_{km}^2\delta_{pq}\\
&\ +4(N-1)(N-2)(N-4)^2f(|y|^2)^2\F{|y|^2\lambda'^{N-4}}{(\lambda'^2+|y|^2)^{N}}\sum_{k,\ell,m}(\partial_\ell H_{mk})^2y_py_q\\
&\ -8(N-2)(N-4)^2f(|y|^2)^2\F{\lambda'^{N-4}}{(\lambda'^2+|y|^2)^{N-1}}\sum_{k,\ell,m}(\partial_\ell H_{mk})^2y_py_q\\
&\ -2(N-2)(N-4)^2f(|y|^2)^2\F{|y|^2\lambda'^{N-4}}{(\lambda'^2+|y|^2)^{N-1}}\sum_{k,\ell,m}(\partial_\ell H_{mk})^2\delta_{pq}\\
&\ +2(N-4)^2f(|y|^2)^2\F{\lambda'^{N-4}}{(\lambda'^2+|y|^2)^{N-2}}\sum_{k,\ell,m}(\partial_\ell H_{mk})^2\delta_{pq}.
\end{align*}
Using Lemma \ref{l9.1} and \ref{l9.2}, we finish the proof.
\end{proof}

\BL
We have
\begin{align*}
&\ \int_{\mathbb R^N} \left.\left\{\sum_{i,j,m,s}(\partial_i \overline{H}_{ms})(\partial_j \overline{H}_{ms})\partial_i\tilde u_0\partial_j\tilde u_0\right\}\right|_{\xi'=0}\\
=&\ |S^{N-1}| (W_{ikj\ell}+W_{i\ell j k})^2\lambda'^{N-4}\\
&\quad\cdot \F{2(N-4)^2}{N(N+2)}\int_0^\infty \left[r^2f'(r^2)+f(r^2)\right]^2 \F{r^{N+3}}{(\lambda'^2+r^2)^{N-2}}\mathrm{d}r,
\end{align*}
\begin{align*}
&\int_{\mathbb R^n}\F{\partial^2}{\partial \xi'_p\partial \xi'_q}\left.\left\{\sum_{i,j,m,s}(\partial_i \overline{H}_{ms})(\partial_j \overline{H}_{ms})\partial_i\tilde u_0\partial_j\tilde u_0\right\}\right|_{\xi'=0}\\
=&\ |S^{N-1}|\sum_{i,j,k}(W_{ikjp}+W_{ipjk})(W_{ikjq}+W_{iqjk})\lambda'^{N-4}\\
&\ \ \cdot \Bigg\{\F{32(N-1)(N-2)(N-4)^2}{N(N+2)(N+4)}\int_0^\infty \left[r^2f'(r^2)+f(r^2)\right]^2\F{r^{N+5}}{(\lambda'^2+r^2)^N}\mathrm{d}r\\
&\qquad - \F{64(N-2)(N-4)^2}{N(N+2)(N+4)}\int_0^\infty \left[r^2f'(r^2)^2+f(r^2)f'(r^2)\right]\F{r^{N+5}}{(\lambda'^2+r^2)^{N-1}}\mathrm{d}r\\
&\qquad -\F{16(N-2)(N-4)^2}{N(N+2)}\int_0^\infty \left[r^2f(r^2)f'(r^2)+f(r^2)^2\right]\F{r^{N+3}}{(\lambda'^2+r^2)^{N-1}}\mathrm{d}r\\
&\qquad +\F{16(N-4)^2}{N(N+2)(N+4)}\int_0^\infty f'(r^2)^2\F{r^{N+5}}{(\lambda'^2+r^2)^{N-2}}\mathrm{d}r\\
&\qquad +\F{2(N-4)^2}{N}\int_0^\infty f(r^2)^2\F{r^{N+1}}{(\lambda'^2+r^2)^{N-2}}\mathrm{d}r\\
&\qquad +\F{8(N-4)^2}{N(N+2)}\int_0^\infty f(r^2)f'(r^2)\F{r^{N+3}}{(\lambda'^2+r^2)^{N-2}}\mathrm{d}r\Bigg\}\\
&\ +|S^{N-1}|\sum_{i,j,k,\ell}(W_{ikj\ell}+W_{i\ell jk})^2\lambda'^{N-4}\delta_{pq}\\
&\ \ \cdot \Bigg\{ \F{4(N-2)(N-4)^2}{N(N+2)} \int_0^\infty \left[r^2f'(r^2)+f(r^2)\right]^2 \\
&\hspace{13em}\cdot\left[\F{2(N-1)}{N+4}\F{r^{N+5}}{(\lambda'^2+r^2)^N}-\F{r^{N+3}}{(\lambda'^2+r^2)^{N-1}}\right]\mathrm{d}r\\
&\qquad -\F{16(N-2)(N-4)^2}{N(N+2)(N+4)}\int_0^\infty \left[r^2f'(r^2)^2+f(r^2)f'(r^2)\right]\F{r^{N+5}}{(\lambda'^2+r^2)^{N-1}}\mathrm{d}r\\
&\qquad +\F{4(N-4)^2}{N(N+2)(N+4)}\int_0^\infty f'(r^2)^2\F{r^{N+5}}{(\lambda'^2+r^2)^{N-2}}\mathrm{d}r\Bigg\}.
\end{align*}
\EL

\begin{proof}
Direct calculation shows
\begin{align*}
&\ \sum_{i,j,m,s}(\partial_i \overline{H}_{ms})(\partial_j \overline{H}_{ms})\partial_i\tilde u_0\partial_j\tilde u_0\\
=&\ \sum_{i,j,m,s}[2f'(|y|^2)y_i H_{ms} + f(|y|^2)\partial_i H_{ms}][2f'(|y|^2)y_j H_{ms} + f(|y|^2)\partial_j H_{ms}]\\
&\qquad \cdot (N-4)^2 \tilde u_0^2(y)\F{(y_i-\xi'_i)(y_j-\xi'_j)}{(\lambda'^2+|y-\xi'|^2)^2}\\
=&\ 4(N-4)^2\Big[|y|^2f'(|y|^2)+f(|y|^2)\Big]^2 \F{\lambda'^{N-4}}{(\lambda'^2+|y-\xi'|^2)^{N-2}} \sum_{m,s}H_{ms}^2\\
&\ -8(N-4)^2\Big[|y|^2f'(|y|^2)^2+f(|y|^2)f'(|y|^2)\Big]\F{\lambda'^{N-4}}{(\lambda'^2+|y-\xi'|^2)^{N-2}}\sum_{j,m,s}H_{ms}^2 y_j\xi'_j\\
&\ -4(N-4)^2\Big[|y|^2f(|y|^2)f'(|y|^2)+f(|y|^2)^2\Big]\F{\lambda'^{N-4}}{(\lambda'^2+|y-\xi'|^2)^{N-2}}\sum_{j,m,s}H_{ms}(\partial_j H_{ms}) \xi'_j\\
&\ +4(N-4)^2f'(|y|^2)^2\F{\lambda'^{N-4}}{(\lambda'^2+|y-\xi'|^2)^{N-2}}\sum_{i,j,m,s}H_{ms}^2 y_i y_j\xi'_i\xi'_j\\
&\ +(N-4)^2f(|y|^2)^2\F{\lambda'^{N-4}}{(\lambda'^2+|y-\xi'|^2)^{N-2}}\sum_{i,j,m,s}(\partial_i H_{ms})(\partial_j H_{ms})\xi'_i\xi'_j\\
&\ +4(N-4)^2f(|y|^2)f'(|y|^2)\F{\lambda'^{N-4}}{(\lambda'^2+|y-\xi'|^2)^{N-2}}\sum_{i,j,m,s}H_{ms}(\partial_j H_{ms})y_i\xi'_i\xi'_j,
\end{align*}
\begin{align*}
&\F{\partial^2}{\partial \xi'_p\partial \xi'_q}\left.\left\{\sum_{i,j,m,s}(\partial_i \overline{H}_{ms})(\partial_j \overline{H}_{ms})\partial_i\tilde u_0\partial_j\tilde u_0\right\}\right|_{\xi=0}\\
=&\ 16(N-1)(N-2)(N-4)^2\Big[|y|^2f'(|y|^2)+f(|y|^2)\Big]^2 \F{\lambda'^{N-4}}{(\lambda'^2+|y|^2)^{N}} \sum_{m,s}H_{ms}^2y_py_q\\
&\ -8(N-2)(N-4)^2\Big[|y|^2f'(|y|^2)+f(|y|^2)\Big]^2 \F{\lambda'^{N-4}}{(\lambda'^2+|y|^2)^{N-1}}\sum_{m,s}H_{ms}^2\delta_{pq}\\
&\ -32(N-2)(N-4)^2\Big[|y|^2f'(|y|^2)^2+f(|y|^2)f'(|y|^2)\Big]\F{\lambda'^{N-4}}{(\lambda'^2+|y|^2)^{N-1}}\sum_{m,s}H_{ms}^2 y_py_q\\
&\ -8(N-2)(N-4)^2\Big[|y|^2f(|y|^2)f'(|y|^2)+f(|y|^2)^2\Big]\F{\lambda'^{N-4}}{(\lambda'^2+|y|^2)^{N-1}}\\
&\qquad \cdot \sum_{m,s}\Big[H_{ms}(\partial_q H_{ms})y_p+H_{ms}(\partial_p H_{ms})y_q\Big]\\
&\ +8(N-4)^2f'(|y|^2)^2\F{\lambda'^{N-4}}{(\lambda'^2+|y|^2)^{N-2}}\sum_{m,s}H_{ms}^2y_py_q\\
&\ +2(N-4)^2f(|y|^2)^2\F{\lambda'^{N-4}}{(\lambda'^2+|y|^2)^{N-2}}\sum_{i,j,m,s}(\partial_p H_{ms})(\partial_q H_{ms})\\
&\ +4(N-4)^2f(|y|^2)f'(|y|^2)\F{\lambda'^{N-4}}{(\lambda'^2+|y|^2)^{N-2}}\sum_{m,s}\left[H_{ms}(\partial_q H_{ms})y_p+H_{ms}(\partial_p H_{ms})y_q\right].
\end{align*}
Lemma \ref{l9.1} and \ref{l9.2} then yield the result.
\end{proof}

\BL
It holds
\begin{align*}
\sum_{i,j,m,s}\Big[\overline H_{ms}(\partial_{s}\overline H_{ij})&-\overline H_{si}(\partial_s\overline H_{mj}) +\overline H_{sj}(\partial_i\overline H_{ms})\\
&-\overline H_{ms}(\partial_{i}\overline H_{sj})\Big]\partial_m (\partial_i\tilde u_0\partial_j\tilde u_0)=0.
\end{align*}
\EL

\begin{proof}
Direct computation shows
\begin{align*}
&\ \sum_{i,j,m,s}\overline H_{ms} (\partial_s \overline H_{ij}) \partial_{m}[\partial_i\tilde u_0\partial_j\tilde u_0]\\
=&\ 2\sum_{i,j,m,s}\overline H_{ms} (\partial_s \overline H_{ij}) \partial_{im} \tilde u_0\partial_j\tilde u_0\\
=&\ -2(N-4)^2 f(|y|^2)^2\F{\lambda'^{N-4}}{(\lambda'^2+|y-\xi'|^2)^{N-2}}\sum_{i,j} H_{ij}^2\\
&\ -2(N-4)^2  f(|y|^2)^2 \F{\lambda'^{N-4}}{(\lambda'^2+|y-\xi'|^2)^{N-2}} \sum_{i,j,m}H_{im}(\partial_m H_{ij})\xi'_j\\
&\ +4(N-2)(N-4)^2f(|y|^2)^2\F{\lambda'^{N-4}}{(\lambda'^2+|y-\xi'|^2)^{N-1}}\sum_{i,j,m}H_{im}H_{jm}\xi'_i\xi'_j\\
&\ +2(N-2)(N-4)^2f(|y|^2)^2\F{\lambda'^{N-4}}{(\lambda'^2+|y-\xi'|^2)^{N-1}}\sum_{i,j,m,s}H_{ms}(\partial_s H_{ij})\xi'_i\xi'_j\xi'_m
\end{align*}
and
\begin{align*}
&\ \sum_{i,j,m,s} \overline H_{is} (\partial_s \overline H_{jm})  \partial_{m}[\partial_i\tilde u_0\partial_j\tilde u_0]\\
=&\ \sum_{i,j,m,s} \overline H_{is} (\partial_s \overline H_{jm}) [\partial_{im}\tilde u_0 \partial_j\tilde u_0 + \partial_i\tilde u_0\partial_{jm}\tilde u_0] \\
=&\ -(N-4)^2 f(|y|^2)^2\F{\lambda'^{N-4}}{(\lambda'^2+|y-\xi'|^2)^{N-2}}\sum_{i,j} H_{ij}^2\\
&\ -(N-4)^2 f(|y|^2)^2 \F{\lambda'^{N-4}}{(\lambda'^2+|y-\xi'|^2)^{N-2}} \sum_{i,j,m}H_{im}(\partial_m H_{ij})\xi'_j\\
&\ +4(N-2)(N-4)^2f(|y|^2)^2\F{\lambda'^{N-4}}{(\lambda'^2+|y-\xi'|^2)^{N-1}}\sum_{i,j,m}H_{im}H_{jm}\xi'_i\xi'_j\\
&\ +2(N-2)(N-4)^2f(|y|^2)^2\F{\lambda'^{N-4}}{(\lambda'^2+|y-\xi'|^2)^{N-1}}\sum_{i,j,m,s}H_{is}(\partial_s H_{jm})\xi'_i\xi'_j\xi'_m.
\end{align*}
Thus
\begin{align*}
&\ \sum_{i,j,m,s}\left[\overline H_{ms} (\partial_s \overline H_{ij})-\overline H_{is} (\partial_s \overline H_{jm}) \right]\partial_{m}[\partial_i\tilde u_0\partial_j\tilde u_0]\\
=&\ -(N-4)^2 f(|y|^2)^2\F{\lambda'^{N-4}}{(\lambda'^2+|y-\xi'|^2)^{N-2}}\sum_{i,j} H_{ij}^2\\
&\  -(N-4)^2 f(|y|^2)^2 \F{\lambda'^{N-4}}{(\lambda'^2+|y-\xi'|^2)^{N-2}} \sum_{i,j,m}H_{im}(\partial_m H_{ij})\xi'_j.
\end{align*}
On the other hand,
\begin{align*}
&\ \sum_{i,j,m,s} \left[\overline H_{sj}  (\partial_i \overline H_{ms}) -\overline H_{ms}  (\partial_i \overline H_{js})\right] (\partial_m(\partial_i\tilde u_0)(\partial_j\tilde u_0))\\
=&\ \sum_{i,j,m,s} f(|y|^2)^2[H_{js}(\partial_i H_{ms})-H_{ms}(\partial_i H_{js})](\partial_{im}\tilde u_0\partial_j\tilde u_0+\partial_i\tilde u_0\partial_{jm}\tilde u_0)\\
=&\ (N-4)^2f(|y|^2)^2 \F{\lambda'^{N-4}}{(\lambda'^2+|y-\xi'|^2)^{N-2}} \sum_{i,j}H_{ij}^2 \\
&\ +(N-4)^2 f(|y|^2)^2\F{\lambda'^{N-4}}{(\lambda'^2+|y-\xi'|^2)^{N-2}} \sum_{i,j,s}H_{is}(\partial_i H_{js})\xi'_j.
\end{align*}
The lemma follows immediately.
\end{proof}

\BL
We have
\begin{align*}
\int_{\mathbb R^N} \left.\left\{\sum_{i,j,m,s}\overline{H}_{si}(\partial_{mm}\overline{H}_{js})\partial_i\tilde u_0\partial_j\tilde u_0\right\}\right|_{\xi'=0}=0,
\end{align*}
\begin{align*}
&\ \int_{\mathbb R^n}\F{\partial^2}{\partial\xi'_p\partial\xi'_q} \left.\left\{\sum_{i,j,m,s}\overline{H}_{si}(\partial_{mm}\overline{H}_{js})\partial_i\tilde u_0\partial_j\tilde u_0\right\}\right|_{\xi'=0}\\
=&\ |S^{N-1}|(W_{ikjp}+W_{ipjk})(W_{ikjq}+W_{iqjk})\lambda'^{N-4}\\
&\cdot \F{2(N-4)^2}{N(N+2)}\int_0^\infty \Big[(N+4)f(r^2)f'(r^2)+2r^2f(r^2)f''(r^2)\Big]\F{r^{N+3}}{(\lambda'^2+r^2)^{N-2}}\mathrm{d}r.
\end{align*}
\EL

\begin{proof}
It is directly checked that
\begin{align*}
&\ \sum_{i,j,m,s}\overline{H}_{si}(\partial_{mm}\overline{H}_{js})\partial_i\tilde u_0\partial_j\tilde u_0\\
=&\ (N-4)^2\Big[(2N+8)f(|y|^2)f'(|y|^2)+4|y|^2f(|y|^2)f''(|y|^2)\Big]\F{\lambda'^{N-4}}{(\lambda'^2+|y-\xi'|^2)^{N-2}}\\
&\qquad\cdot \sum_{i,j,s}H_{is}H_{js}\xi'_i\xi'_j,
\end{align*}
and
\begin{align*}
&\ \F{\partial^2}{\partial\xi'_p\partial\xi'_q} \left.\left\{\sum_{i,j,m,s}\overline{H}_{si}(\partial_{mm}\overline{H}_{js})\partial_i\tilde u_0\partial_j\tilde u_0\right\}\right|_{\xi=0}\\
=&\ 4(N-4)^2\Big[(N+4)f(|y|^2)f'(|y|^2)+2|y|^2f(|y|^2)f''(|y|^2)\Big]\F{\lambda'^{N-4}}{(\lambda'^2+|y-\xi'|^2)^{n-2}}\\
&\qquad\cdot \sum_{s}H_{ps}H_{qs}.
\end{align*}
we conclude the proof by using Lemma \ref{l9.1} and \ref{l9.2}.
\end{proof}

\BL
The following hold
\begin{align*}
&\ \int_{\mathbb R^N}\left.\left\{\sum_{i,k,\ell,m}\left[(\partial_{i\ell}\overline H_{mk})^2+(\partial_\ell \overline H_{mk})(\partial_{ii\ell}\overline H_{mk})\right]\tilde u_0^2\right\}\right|_{\xi'=0}\\
=&\ |S^{N-1}|\sum_{i,j,k,\ell}(W_{ikj\ell}+W_{i\ell jk})^2\lambda'^{N-4}\\
&\ \cdot \Bigg\{\F{2}{N(N+2)}\int_0^\infty \Big[3(N+8)f'(r^2)^2+2(N+8)f(r^2)f''(r^2)\\
&\qquad\qquad +2(N+18)r^2f'(r^2)f''(r^2)+4r^4f''(r^2)^2\\
&\qquad\qquad +4r^2f(r^2)f'''(r^2)+4r^4f'(r^2)f'''(r^2)\Big]\F{r^{N+3}}{(\lambda'^2+r^2)^{N-4}}\mathrm{d}r\\
&\quad +\F{2}{N}\int_0^\infty \Big[4r^2f'(r^2)^2 +(N+8)f(r^2)f'(r^2)+2r^2 f(r^2)f''(r^2)\Big]\F{r^{N+1}}{(\lambda'^2+r^2)^{N-4}}\mathrm{d}r\\
&\quad +\int_0^\infty f(r^2)^2\F{r^{N-1}}{(\lambda'^2+r^2)^{N-4}}\mathrm{d}r \Bigg\},
\end{align*}
and
\begin{align*}
&\int_{\mathbb R^N}\F{\partial^2}{\partial \xi'_p\partial \xi'_q}\left.\left\{\sum_{i,k,\ell,m}\left[(\partial_{i\ell}\overline{H}_{mk})^2+(\partial_\ell \overline{H}_{mk})(\partial_{ii\ell}\overline{H}_{mk})\right]\tilde u_0^2\right\}\right|_{\xi=0}\\
=&\ |S^{N-1}|(W_{ikjp}+W_{ipjk})(W_{ikjq}+W_{iqjk})\lambda'^{N-4}\\
&\ \ \Bigg\{\F{8(N-3)(N-4)}{N(N+2)(N+4)}\\
&\quad \cdot\int_0^\infty\Big[12(N+8)f'(r^2)^2+16r^4f''(r^2)^2+8(N+18)r^2f'(r^2)f''(r^2)\\
&\qquad\quad +8(N+8)f(r^2)f''(r^2)+16r^2f(r^2)f'''(r^2)+16r^4f'(r^2)f'''(r^2)\Big]\\
&\qquad \cdot \F{r^{N+5}}{(\lambda'^2+r^2)^{N-2}}\\
&\quad +\F{8(N-3)(N-4)}{N(N+2)}\\
&\quad \cdot\int_0^\infty\Big[8r^2f'(r^2)^2+2(N+8)f(r^2)f'(r^2)+4r^2f(r^2)f''(r^2)\Big]\F{r^{N+3}}{(\lambda'^2+|x|^2)^{N-2}}\Bigg\}\\
&\ +|S^{N-1}|(W_{ikj\ell}+W_{i\ell jk})^2\lambda'^{N-4}\delta_{pq}\\
&\ \ \Bigg\{\F{N-4}{N(N+2)}\\
&\quad \cdot\int_0^\infty\Big[12(N+8)f'(r^2)^2+16r^4f''(r^2)^2+8(N+18)r^2f'(r^2)f''(r^2)\\
&\qquad +8(N+8)f(r^2)f''(r^2)+16r^2f(r^2)f'''(r^2)+16r^4f'(r^2)f'''(r^2)\Big]\\
&\qquad \cdot \left[\F{2(N-3)}{N+4}\F{r^{N+5}}{(\lambda'^2+r^2)^{N-2}}-\F{r^{N+3}}{(\lambda'^2+r^2)^{N-3}}\right]\\
&\quad +\F{2(N-4)}{N}\int_0^\infty\Big[8r^2f'(r^2)^2+2(N+8)f(r^2)f'(r^2)+4r^2f(r^2)f''(r^2)\Big]\\
&\qquad \cdot \left[\F{2(N-3)}{N+2}\F{r^{N+3}}{(\lambda'^2+r^2)^{N-2}}-\F{r^{N+1}}{(\lambda'^2+r^2)^{N-3}}\right]\\
&\quad +2(N-4)\int_0^\infty f(r^2)^2\left[\F{2(N-3)}{N}\F{r^{N+1}}{(\lambda'^2+r^2)^{N-2}}-\F{r^{N-1}}{(\lambda'^2+r^2)^{N-3}}\right]\Bigg\}.
\end{align*}
\EL

\begin{proof}
We have
\begin{align*}
&\ \sum_{i,\ell}(\partial_{i\ell}\overline H_{mk})^2\\
=&\ \sum_{i,\ell}\Big[2\delta_{i\ell}f'(|y|^2) H_{mk}+2f'(|y|^2)y_\ell(\partial_i H_{mk})+2f'(|y|^2)y_i(\partial_\ell H_{mk})\\
&\qquad\quad +f(|y|^2)\partial_{i\ell}H_{mk}+4f''(|y|^2)y_\ell y_i H_{mk}\Big]^2\\
=&\ 4N f'(|y|^2)^2H_{mk}^2 +8|y|^2f'(|y|^2)^2\sum_i(\partial_i H_{mk})^2+f(|y|^2)^2\sum_{i,\ell}(\partial_{i\ell}H_{mk})^2\\
&\ +16f'(|y|^2)^2\sum_i y_i H_{mk}(\partial_i H_{mk})+4f'(|y|^2)f(|y|^2)\sum_i H_{mk}(\partial_{ii}H_{mk})\\
&\ +8f'(|y|^2)^2\sum_{i,\ell}y_i y_\ell(\partial_i H_{mk})(\partial_\ell H_{mk})+8f'(|y|^2)f(|y|^2)\sum_{i,\ell} y_\ell(\partial_i H_{mk})(\partial_{i\ell}H_{mk})\\
&\ +16f''(|y|^2)^2|y|^2 H_{mk}^2+80f'(|y|^2)f''(|y|^2)|y|^2H_{mk}^2\\
&\ +8f(|y|^2)f''(|y|^2)\sum_{i,\ell}y_i y_\ell H_{mk}(\partial_{i\ell}H_{mk})\\
=&\ \Big[(4N+64)f'(|y|^2)^2+16|y|^4f''(|y|^2)^2+80|y|^2f'(|y|^2)f''(|y|^2)+16f(|y|^2)f''(|y|^2)\Big]H_{mk}^2\\
&\ +8\Big[|y|^2f'(|y|^2)^2+f(|y|^2)f'(|y|^2)\Big]\sum_i(\partial_i H_{mk})^2+f(|y|^2)^2\sum_{i,\ell}(\partial_{i\ell}H_{mk})^2,
\end{align*}
\begin{align*}
&\ \sum_{i,\ell}(\partial_\ell \overline H_{mk})(\partial_{ii\ell}\overline H_{mk})\\
=&\ \sum_{\ell}\bigg\{\Big[(2N+8)f'(|y|^2)+4|y|^2f''(|y|^2)\Big](\partial_\ell H_{mk})\\
&\qquad\qquad+\Big[(4N+24)f''(|y|^2)+8|y|^2f'''(|y|^2)\Big]y_\ell H_{mk}\bigg\}\\
&\qquad \cdot \Big[2f'(|y|^2)y_\ell H_{mk} + f(|y|^2)\partial_\ell H_{mk}\Big]\\
=&\ \Big[8(N+4)f'(|y|^2)^2+8(N+8)|y|^2f'(|y|^2)f''(|y|^2)+8(N+6)f(|y|^2)f''(|y|^2)\\
&\qquad +16|y|^2f(|y|^2)f'''(|y|^2)+16|y|^4f'(|y|^2)f'''(|y|^2)\Big]H_{mk}^2\\
&\ +\Big[2(N+4)f(|y|^2)f'(|y|^2)+4|y|^2 f(|y|^2)f''(|y|^2)\Big]\sum_\ell (\partial_\ell H_{mk})^2,
\end{align*}
and
\begin{align*}
&\F{\partial^2}{\partial \xi'_p\partial \xi'_q}\left.\left\{\sum_{i,k,\ell,m}\left[(\partial_{i\ell}\overline{H}_{mk})^2+(\partial_\ell \overline{H}_{mk})(\partial_{ii\ell}\overline{H}_{mk})\right]\tilde u_0^2\right\}\right|_{\xi'=0}\\
=&\ 4(N-3)(N-4)\Big[12(N+8)f'(|y|^2)^2+16|y|^4f''(|y|^2)^2+8(N+18)|y|^2f'(|y|^2)f''(|y|^2)\\
&\quad +8(N+8)f(|y|^2)f''(|y|^2)+16|y|^2f(|y|^2)f'''(|y|^2)+16|y|^4f'(|y|^2)f'''(|y|^2)\Big]\\
&\qquad \cdot\F{\lambda'^{N-4}}{(\lambda'^2+|y|^2)^{N-2}}\sum_{k,m}H_{km}^2y_py_q\\
&-2(N-4)\Big[12(N+8)f'(|y|^2)^2+16|y|^4f''(|y|^2)^2+8(N+18)|y|^2f'(|y|^2)f''(|y|^2)\\
&\quad +8(N+8)f(|y|^2)f''(|y|^2)+16|y|^2f(|y|^2)f'''(|y|^2)+16|y|^4f'(|y|^2)f'''(|y|^2)\Big]\\
&\qquad \cdot\F{\lambda'^{N-4}}{(\lambda'^2+|y|^2)^{N-3}}\sum_{k,m}H_{km}^2\delta_{pq}\\
&\ +4(N-3)(N-4)\Big[8|y|^2f'(|y|^2)^2+2(N+8)f(|y|^2)f'(|y|^2)+4|y|^2f(|y|^2)f''(|y|^2)\Big]\\
&\qquad \cdot\F{\lambda'^{N-4}}{(\lambda'^2+|y|^2)^{N-2}}\sum_{k,\ell,m}(\partial_\ell H_{km})^2y_py_q\\
&\ -2(N-4)\Big[8|y|^2f'(|y|^2)^2+2(N+8)f(|y|^2)f'(|y|^2)+4|y|^2f(|y|^2)f''(|y|^2)\Big]\\
&\qquad \cdot\F{\lambda'^{N-4}}{(\lambda'^2+|y|^2)^{N-3}}\sum_{k,\ell,m}(\partial_\ell H_{km})^2\delta_{pq}\\
&\ +4(N-3)(N-4)f(|y|^2)^2\F{\lambda'^{N-4}}{(\lambda'^2+|y|^2)^{N-2}}\sum_{i,k,\ell,m}(\partial_{i\ell} H_{km})^2y_py_q\\
&\ -2(N-4)f(|y|^2)^2\F{\lambda'^{N-4}}{(\lambda'^2+|y|^2)^{N-3}}\sum_{i,k,\ell,m}(\partial_{i\ell} H_{km})^2\delta_{pq}.
\end{align*}
Lemma \ref{l9.1} and \ref{l9.2} then give the result.
\end{proof}

\BL
We have
\begin{align*}
&\ \int_{\mathbb R^N}\left.\left\{\sum_{i,j,m,s}(\partial_{mm}\overline{H}_{ij})(\partial_{ss}\overline{H}_{ij})\tilde u_0^2\right\}\right|_{\xi'=0}\\
=&\ |S^{N-1}|\sum_{i,j,k,\ell}(W_{ikj\ell}+W_{i\ell jk})^2\lambda'^{N-4}\\
&\ \cdot \F{1}{2N(N+2)}\int_0^\infty \Big[2(N+4)f'(r^2)+4r^2f''(r^2)\Big]^2\F{r^{N+3}}{(\lambda'^2+r^2)^{N-4}}\mathrm{d}r
\end{align*}
and
\begin{align*}
&\int_{\mathbb R^N}\F{\partial^2}{\partial \xi_p\partial \xi_q}\left.\left\{\sum_{i,j,m,s}(\partial_{mm}\overline{H}_{ij})(\partial_{ss}\overline{H}_{ij})\tilde u_0^2\right\}\right|_{\xi'=0}\\
=&\ |S^{N-1}|\sum_{i,j,k}(W_{ikjp}+W_{ipjk})(W_{ikjq}+W_{iqjk})\lambda'^{N-4}\\
&\quad\cdot\F{8(N-3)(N-4)}{N(N+2)(N+4)} \int_0^\infty \left[2(N+4)f'(r^2)+4r^2f''(r^2)\right]^2\F{r^{N+5}}{(\lambda'^2+r^2)^{N-2}}\mathrm{d}r\\
&\ +|S^{N-1}|\sum_{i,j,k,\ell}(W_{ikj\ell}+W_{i\ell jk})^2\lambda'^{N-4}\delta_{pq}\\
&\quad\cdot \F{N-4}{N(N+2)}\int_0^\infty \left[2(N+4)f'(r^2)+4r^2f''(r^2)\right]^2 \\
&\qquad\qquad\qquad\qquad\qquad\cdot\left[\F{2(N-3)}{N+4}\F{r^{N+5}}{(\lambda'^2+r^2)^{N-2}}-\F{r^{N+3}}{(\lambda'^2+r^2)^{N-3}}\right]\mathrm{d}r.
\end{align*}
\EL

\begin{proof}
Since
\begin{align*}
\sum_m\partial_{mm}\overline H_{ij}&=\sum_m \bigg(4f'(|y|^2)y_m \partial_m H_{ij}+2f'(|y|^2)H_{ij}+4f''(|y|^2)y_m^2 H_{ij}\\
&\qquad\quad +f(|y|^2)\partial_{mm}H_{ij}\bigg)\\
&= \Big[2(N+4)f'(|y|^2)+4|y|^2f''(|y|^2)\Big]H_{ij},
\end{align*}
and
\begin{align*}
&\F{\partial^2}{\partial \xi'_p\partial \xi'_q}\left.\left\{\sum_{i,j,m,s}(\partial_{mm}\overline{H}_{ij})(\partial_{ss}\overline{H}_{ij})u_{(\xi',\varepsilon)}^2\right\}\right|_{\xi'=0}\\
=&\ 4(N-3)(N-4)\Big[2(N+4)f'(|y|^2)+4|y|^2f''(|y|^2)\Big]^2\F{\lambda'^{N-4}}{(\lambda'^2+|y|^2)^{N-2}}\sum_{i,j}H_{ij}^2y_p y_q\\
&\ -2(N-4)\Big[2(N+4)f'(|y|^2)+4|y|^2f''(|y|^2)\Big]^2\F{\lambda'^{N-4}}{(\lambda'^2+|y|^2)^{N-3}}\sum_{i,j}H_{ij}^2\delta_{pq},
\end{align*}
we obtain the result by Lemma \ref{l9.1} and \ref{l9.2}.
\end{proof}

\BL\label{l8.2}
There hold
\begin{align*}
\left.\sum_{i,j,s}(\partial_{ss} \overline H_{ij})(\partial_{ij}\tilde u_0)w\right|_{\xi'=0} &=0,
\qquad &\F{\partial^2}{\partial\xi'_p\partial\xi'_q}\left.\sum_{i,j,s}(\partial_{ss} \overline H_{ij})(\partial_{ij}\tilde u_0)w\right|_{\xi'=0} &=0,\\
\left.\sum_{i,j,s}(\partial_s \overline{H}_{ij})(\partial_{sij}\tilde u_0)w\right|_{\xi'=0} &=0,
\qquad &\F{\partial^2}{\partial\xi'_p\partial\xi'_q}\left.\sum_{i,j,s}(\partial_s \overline{H}_{ij})(\partial_{sij}\tilde u_0)w\right|_{\xi'=0} &=0,\\
\left.\sum_{i,j,s}\overline{H}_{ij}(\partial_{ssij}\tilde u_0)w\right|_{\xi'=0} &=0,
\qquad &\F{\partial^2}{\partial\xi'_p\partial\xi'_q}\left.\sum_{i,j,s}\overline{H}_{ij}(\partial_{ssij}\tilde u_0)w\right|_{\xi'=0} &=0,\\
\left.\sum_{i,j,m}(\partial_{jmm}\overline{H}_{ij})(\partial_i u)w\right|_{\xi'=0} &=0,
\qquad &\F{\partial^2}{\partial\xi'_p\partial\xi'_q}\left.\sum_{i,j,m}(\partial_{jmm}\overline{H}_{ij})(\partial_i u)w\right|_{\xi'=0} &=0.
\end{align*}
\EL
\begin{proof}
Direct computation shows
\begin{align*}
&\ \sum_{i,j,s}(\partial_{ss} \overline H_{ij})(\partial_{ij}\tilde u_0)\\
=&\ \Big[2N f'(|y|^2)+8f'(|y|^2)+4|y|^2f''(|y|^2)\Big]\sum_{i,j}H_{ij}(\partial_{ij}\tilde u_0)\\
=&\ (N-2)(N-4)\Big[(2N+8)f'(|y|^2)+4|y|^2f''(|y|^2)\Big]\\
&\qquad \cdot \F{\lambda'^{\F{N-4}{2}}}{(\lambda'^2+|y-\xi'|^2)^\F{N}{2}}\sum_{i,j}H_{ij}\xi'_i\xi'_j,
\end{align*}
\begin{align*}
&\ \sum_{i,j,s}(\partial_s \overline{H}_{ij})(\partial_{sij}\tilde u_0)\\
=&\ \sum_{i,j,s}\Big[2f'(|y|^2)y_s H_{ij}+f(|y|^2)(\partial_s H_{ij})\Big](N-2)(N-4)\tilde u_0\\
&\ \cdot\bigg[\F{\delta_{is}(y_j-\xi'_j)+\delta_{js}(y_i-\xi'_i)+\delta_{ij}(y_s-\xi'_s)}{(\lambda'^2+|y-\xi'|^2)^2}-N\F{(y_s-\xi'_s)(y_i-\xi'_i)(y_j-\xi'_j)}{(\lambda'^2+|y-\xi'|^2)^3}\bigg]\\
=&\ -2N(N-2)(N-4)|y|^2f'(|y|^2)\F{\lambda'^{\F{N-4}{2}}}{(\lambda'^2+|y-\xi'|^2)^\F{N+2}{2}}\sum_{i,j}H_{ij}\xi'_i\xi'_j\\
&\ +2N(N-2)(N-4)f'(|y|^2)\F{\lambda'^{\F{N-4}{2}}}{(\lambda'^2+|y-\xi'|^2)^\F{N+2}{2}}\sum_{i,j,s}H_{ij}y_s\xi'_s\xi'_i\xi'_j\\
&\ +N(N-2)(N-4)f(|y|^2)\F{\lambda'^{\F{N-4}{2}}}{(\lambda'^2+|y-\xi'|^2)^\F{N+2}{2}}\sum_{i,j,s}(\partial_s H_{ij})\xi'_s\xi'_i\xi'_j,
\end{align*}
\begin{align*}
&\ \sum_{i,j,s}\overline{H}_{ij}(\partial_{ssij}\tilde u_0)\\
=&\ N(N-2)(N+2)(N-4)f(|y|^2)\F{|y-\xi'|^2\lambda'^{\F{N-4}{2}}}{(\lambda'^2+|y-\xi'|^2)^\F{N+4}{2}}\sum_{i,j}H_{ij}\xi'_i\xi'_j\\
&\ -N(N-2)(N-4)(N+4)f(|y|^2)\F{\lambda'^{\F{N-4}{2}}}{(\lambda'^2+|y-\xi'|^2)^\F{N+2}{2}}\sum_{i,j}H_{ij}\xi'_i\xi'_j,
\end{align*}
\begin{align*}
\sum_{i,j,m}(\partial_{jmm}\overline{H}_{ij})(\partial_i u)=0.
\end{align*}
Recall the equation for $w$ and note that if $\xi'=0$, $R_1(y)=0$, so $w|_{\xi'=0}\equiv 0$. The Lemma is easily concluded.
\end{proof}

 Combining the above identities,
we have the following proposition.

\BP\label{p9.16}
It holds that
\begin{align*}
&\ F(0,\lambda')\\
=&\ \F{N-4}{2}|S^{N-1}|\sum_{i,j,k,\ell}(W_{ikj\ell}+W_{i\ell jk})^2\lambda'^{N-4} \\
&\cdot \Bigg\{-\F{a_N(N-4)}{N(N+2)}\int_0^\infty \left[r^2 f'(r^2)^2 +2 f(r^2)f'(r^2)\right]\F{r^{N+5}}{(\lambda'^2+r^2)^{N-2}}\mathrm{d}r\\
&\qquad - \F{a_N(N-4)}{2N}\int_0^\infty f(r^2)^2 \F{r^{N+3}}{(\lambda'^2+r^2)^{N-2}}\mathrm{d}r\\
&\qquad -\F{b_N(N-4)}{N(N+2)}\int_0^\infty \left[r^2f'(r^2)+f(r^2)\right]^2 \F{r^{N+3}}{(\lambda'^2+r^2)^{N-2}}\mathrm{d}r\\
&\qquad +\F{1}{2N(N-1)(N+2)}\int_0^\infty \Big[3(N+8)f'(r^2)^2+2(N+8)f(r^2)f''(r^2)\\
&\hspace{9em} +2(N+18)r^2f'(r^2)f''(r^2)+4r^4f''(r^2)^2\\
&\hspace{9em} +4r^2f(r^2)f'''(r^2)+4r^4f'(r^2)f'''(r^2)\Big]\F{r^{N+3}}{(\lambda'^2+r^2)^{N-4}}\mathrm{d}r\\
&\quad +\F{1}{2N(N-1)}\int_0^\infty \Big[4r^2f'(r^2)^2 +(N+8)f(r^2)f'(r^2)+2r^2 f(r^2)f''(r^2)\Big]\\
&\hspace{25em}\cdot\F{r^{N+1}}{(\lambda'^2+r^2)^{N-4}}\mathrm{d}r\\
&\quad +\F{1}{4(N-1)}\int_0^\infty f(r^2)^2\F{r^{N-1}}{(\lambda'^2+r^2)^{N-4}}\mathrm{d}r\\
&\quad -\F{1}{N(N-2)^2(N+2)}\int_0^\infty \Big[(N+4)f'(r^2)+2r^2f''(r^2)\Big]^2\F{r^{N+3}}{(\lambda'^2+r^2)^{N-4}}\mathrm{d}r\Bigg\}.
\end{align*}
\EP

Next we compute the Hessian of $F$ at $(0,\lambda')$. Because of Lemma \ref{l8.2}, it is obvious that
\begin{align*}
&\ \F{\partial^2}{\partial\xi'_p\partial\xi'_q}\Bigg\{\sum_{i,j,s}\Big[2\overline H_{ij}(\partial_{ijss}\tilde u_0)+ 2(\partial_s\overline H_{ij})( \partial_{ijs}\tilde u_0)+ (\partial_{ss}\overline H_{ij})(\partial_{ij}\tilde u_0)\nonumber\\
&\qquad\qquad\qquad -\F{b_N}{2}(\partial_{ss}\overline H_{ij})(\partial_{ij}\tilde u_0)-\F{b_N}{2}(\partial_{jss}\overline H_{ij})(\partial_{i}\tilde u_0)\Big]\bar w\Bigg\}\Bigg|_{\xi'=0}\\
&=0.
\end{align*}

\BP
The second order partial derivatives of $F(\xi',\lambda')$ at $(0,\lambda')$ are given by
\begin{align*}
&\ \F{\partial^2}{\partial\xi'_p\partial\xi'_q}F(0,\lambda')\\
=&\ \F{(N-4)^2}{N(N+2)}|S^{N-1}|\sum_{i,j,k}(W_{ikjp}+W_{ipjk})(W_{ikjq}+W_{iqjk})\lambda'^{N-4}\\
&\cdot \Bigg\{-(N-2)\int_0^\infty f(r^2)^2 \left[\F{Nr^{N+5}}{(\lambda'^2+r^2)^N}-\F{(N+2)r^{N+3}}{(\lambda'^2+r^2)^{N-1}}\right]\mathrm{d}r\\
&\hspace{1.6em} -\F{8a_N(N-2)}{N+4}\int_0^\infty \left[r^2f'(r^2)^2+2f(r^2)f'(r^2)\right]\\
&\hspace{15em}\cdot\left[\F{(N-1)r^{N+7}}{(\lambda'^2+r^2)^N}-\F{2r^{N+5}}{(\lambda'^2+r^2)^{N-1}}\right]\mathrm{d}r\\
&\hspace{1.6em} -2a_N(N-2)\int_0^\infty f(r^2)^2 \left[\F{(N-1)r^{N+5}}{(\lambda'^2+r^2)^N}-\F{2r^{N+3}}{(\lambda'^2+r^2)^{N-1}}\right]\mathrm{d}r\\
&\hspace{1.6em} -\F{8b_N(N-1)(N-2)}{N+4}\int_0^\infty \left[r^2f'(r^2)+f(r^2)\right]^2\F{r^{N+5}}{(\lambda'^2+r^2)^N}\mathrm{d}r\\
&\hspace{1.6em} +\F{16b_N(N-2)}{N+4}\int_0^\infty \left[r^2f'(r^2)^2+f(r^2)f'(r^2)\right]\F{r^{N+5}}{(\lambda'^2+r^2)^{N-1}}\mathrm{d}r\\
&\hspace{1.6em} +4b_N(N-2)\int_0^\infty \left[r^2f(r^2)f'(r^2)+f(r^2)^2\right]\F{r^{N+3}}{(\lambda'^2+r^2)^{N-1}}\mathrm{d}r\\
&\hspace{1.6em} -\F{4b_N}{N+4}\int_0^\infty f'(r^2)^2\F{r^{N+5}}{(\lambda'^2+r^2)^{N-2}}\mathrm{d}r\\
&\hspace{1.6em} -\F{b_N(N+2)}{2}\int_0^\infty f(r^2)^2\F{r^{N+1}}{(\lambda'^2+r^2)^{N-2}}\mathrm{d}r\\
&\hspace{1.6em} -2b_N\int_0^\infty f(r^2)f'(r^2)\F{r^{N+3}}{(\lambda'^2+r^2)^{N-2}}\mathrm{d}r\\
&\hspace{1.6em} +b_N\int_0^\infty \Big[(N+4)f(r^2)f'(r^2)+2r^2f(r^2)f''(r^2)\Big]\F{r^{N+3}}{(\lambda'^2+r^2)^{N-2}}\mathrm{d}r\\
&\hspace{1.6em} +\F{4(N-3)}{(N-1)(N+4)}\\
&\qquad \cdot\int_0^\infty\Big[3(N+8)f'(r^2)^2+4r^4f''(r^2)^2+2(N+18)r^2f'(r^2)f''(r^2)\\
&\qquad\qquad\quad +2(N+8)f(r^2)f''(r^2)+4r^2f(r^2)f'''(r^2)+4r^4f'(r^2)f'''(r^2)\Big]\\
&\hspace{25em} \cdot \F{r^{N+5}}{(\lambda'^2+r^2)^{N-2}}\mathrm{d}r\\
&\hspace{1.6em} +\F{2(N-3)}{N-1}\int_0^\infty\Big[4r^2f'(r^2)^2+(N+8)f(r^2)f'(r^2)+2r^2f(r^2)f''(r^2)\Big]\\
&\hspace{23em}\cdot\F{r^{N+3}}{(\lambda'^2+r^2)^{N-2}}\mathrm{d}r\\
&\hspace{1.6em} -\F{8(N-3)}{(N-2)^2(N+4)} \int_0^\infty \left[(N+4)f'(r^2)+2r^2f''(r^2)\right]^2\F{r^{N+5}}{(\lambda'^2+r^2)^{N-2}}\mathrm{d}r\Bigg\}\\
&\ +\F{(N-4)^2}{N(N+2)}|S^{N-1}|\sum_{i,j,k,\ell}(W_{ikj\ell}+W_{i\ell jk})^2\lambda'^{N-4}\delta_{pq}\\
&\cdot\Bigg\{-\F{a_N}{N+4}\int_0^\infty \left[r^2f'(r^2)^2+2f(r^2)f'(r^2)\right]\\
&\qquad \cdot \left[\F{2(N-1)(N-2)r^{N+7}}{(\lambda'^2+r^2)^N}-\F{(N-2)(N+8)r^{N+5}}{(\lambda'^2+r^2)^{N-1}}+\F{(N+4)r^{N+3}}{(\lambda'^2+r^2)^{N-2}}\right]\mathrm{d}r\\
&\hspace{1.6em} - \F{a_N}{2}\int_0^\infty f(r^2)^2\\
&\qquad \cdot \left[\F{2(N-1)(N-2)r^{N+5}}{(\lambda'^2+r^2)^N}-\F{(N-2)(N+6)r^{N+3}}{(\lambda'^2+r^2)^{N-1}}+\F{(N+2)r^{N+1}}{(\lambda'^2+r^2)^{N-2}}\right]\mathrm{d}r\\
&\hspace{1.6em} -\F{b_N(N-2)}{N+4} \int_0^\infty \left[r^2f'(r^2)+f(r^2)\right]^2 \\
& \hspace{12em}\cdot\left[\F{2(N-1)r^{N+5}}{(\lambda'^2+r^2)^N}-\F{(N+4)r^{N+3}}{(\lambda'^2+r^2)^{N-1}}\right]\mathrm{d}r\\
&\qquad +\F{4b_N(N-2)}{N+4}\int_0^\infty \left[r^2f'(r^2)^2+f(r^2)f'(r^2)\right]\F{r^{N+5}}{(\lambda'^2+r^2)^{N-1}}\mathrm{d}r\\
&\qquad -\F{b_N}{N+4}\int_0^\infty f'(r^2)^2\F{r^{N+5}}{(\lambda'^2+r^2)^{N-2}}\mathrm{d}r\\
&\hspace{1.6em}  +\F{1}{2(N-1)(N+4)}\\
&\quad \cdot\int_0^\infty\Big[3(N+8)f'(r^2)^2+4r^4f''(r^2)^2+2(N+18)r^2f'(r^2)f''(r^2)\\
&\qquad \qquad +2(N+8)f(r^2)f''(r^2)+4r^2f(r^2)f'''(r^2)+4r^4f'(r^2)f'''(r^2)\Big]\\
&\hspace{16em} \cdot \left[\F{2(N-3)r^{N+5}}{(\lambda'^2+r^2)^{N-2}}-\F{(N+4)r^{N+3}}{(\lambda'^2+r^2)^{N-3}}\right]\mathrm{d}r\\
&\quad +\F{1}{2(N-1)}\int_0^\infty\Big[4r^2f'(r^2)^2+(N+8)f(r^2)f'(r^2)+2r^2f(r^2)f''(r^2)\Big]\\
&\hspace{16em} \cdot \left[\F{2(N-3)r^{N+3}}{(\lambda'^2+r^2)^{N-2}}-\F{(N+2)r^{N+1}}{(\lambda'^2+r^2)^{N-3}}\right]\mathrm{d}r\\
&\quad +\F{N+2}{4(N-1)}\int_0^\infty f(r^2)^2\left[\F{2(N-3)r^{N+1}}{(\lambda'^2+r^2)^{N-2}}-\F{Nr^{N-1}}{(\lambda'^2+r^2)^{N-3}}\right]\mathrm{d}r\\
&\hspace{1.6em}-\F{1}{(N-2)^2(N+4)}\int_0^\infty \left[(N+4)f'(r^2)+2r^2f''(r^2)\right]^2 \\
&\hspace{15em}\cdot\left[\F{2(N-3)r^{N+5}}{(\lambda'^2+r^2)^{N-2}}-\F{(N+4)r^{N+3}}{(\lambda'^2+r^2)^{N-3}}\right]\mathrm{d}r
\Bigg\}.
\end{align*}
\EP

In summary, we have reduced the derivative of $F$ and the Hessian of $F$ to integrals in terms of an auxilliary function $f$. In the next two sections, we choose the auxilliary function $f$ so that $F$ has a strict local minimum at $(0,1)$.  More precisely, we have to choose a function (which is a polynomial) so that the following conditions are satisfied

\noindent
($F_1$) $
\frac{\partial F}{\partial \lambda^{'}} (0,1)=0$;

\noindent
($F_2$) the matrix $ (\F{\partial^2}{\partial\xi'_p\partial \xi'_q}F(0,1))$
is positive definite;

\noindent
($F_3$)  $\F{\partial^2}{\partial\lambda'^2}F(0,1)>0$.

We remark that by Proposition \ref{p8.1},  it holds that $\frac{\partial F}{\partial \xi_p^{'} } (0,1)=0, \frac{\partial^2 F}{\partial \lambda^{'} \xi_p^{'}} (0,1)=0$. Conditions ($F_1$)-($F_2$) ensure that $F$ has a nondegenerate local minimum at $(0,1)$.

Our first choice is a linear function.

\section{Linear function and the case of  $N \geq 52$}

In this section, we show that when $N \geq 52$ the choice of suitable linear fuction satisfies ($F_1$)-($F_2$). (Surprisingly, this dimension $52$ also agrees with the second order Yamabe problem by Brendle \cite{B} in which he also chose a linear function.) The computations are unfortunately complicated even in this case. Many of the computations below are carried out by Mathematica. Since all these computations only involve finding the roots of certain polynomials, the computing errors can be controlled.

Let  the auxiliary function be
\BEN
f(s)=\tau + s.
\EEN
Using the software \textit{Mathematica}, we get the following two  propositions.

\BP
Assume $N> 12$, we have
\begin{align*}
&\ F(0,\lambda')\\
=&\ \F{(N-4)^2}{4( N^2-4)(N-8)}|S^{N-1}|\sum_{i,j,k,\ell}(W_{ikj\ell}+W_{i\ell jk})^2\F{\Gamma\left[\F{N}{2}-3\right] \Gamma\left[\F{N}{2}+3\right]}{ \Gamma[N+1]}I(\lambda'),
\end{align*}
where $\Gamma$ denotes the usual $\Gamma-\text{function}$ and
\begin{align*}
&\ I(\lambda') \\
=&\ -\frac{   \lambda'^8}{(N-12) (N-10) }\left(N^5-4 N^4-80 N^3+208 N^2-32 N-192\right) \\
&\qquad\qquad\qquad -\frac{2   (N-2)\lambda'^6 \tau}{N-10 } \left(N^3-8 N^2+16\right)\\
&\qquad\qquad\qquad\qquad\qquad -\frac{   (N-2)^2 \lambda'^4 \tau^2}{ N+4}\left(N^2-4 N-4\right).
\end{align*}
\EP

\BP
Assume $N> 12$, we have
\begin{align*}
&\ \F{\partial^2}{\partial\xi'_p\partial\xi'_q}F(0,\lambda')\\
=&\ \F{(N-4)^2}{N (N-8)(N+4)(N^2-4)}|S^{N-1}|\frac{\Gamma \left(\frac{N}{2}-3\right) \Gamma \left(\frac{N}{2}+3\right)}{ \Gamma (N)}\\
&\qquad \cdot\Bigg\{\sum_{i,j,k}(W_{ikjp}+W_{ipjk})(W_{ikjq}+W_{iqjk})J_1(\lambda')\\
&\qquad\qquad\qquad\qquad +\sum_{i,j,k,\ell}(W_{ikj\ell}+W_{i\ell jk})^2\delta_{pq}J_2(\lambda')\Bigg\},
\end{align*}
where
\begin{align*}
&\ J_1(\lambda')\\
=&\ -\frac{6 \lambda'^6}{(N-10) } \left(N^5-4 N^4-44 N^3+112 N^2-32N-96\right) \\
&\qquad- 4 (N-2)^2 \lambda'^4  \left(N^2-4 N-4\right)\tau
\end{align*}
and
\begin{align*}
&\ J_2(\lambda')\\
=&\ -\frac{\lambda'^6 (N-2) \left(N^4-N^3-56 N^2+40 N+88\right)}{N-10} \\
&\qquad -\lambda'^4 (N-2)^2 \left(N^2-4 N-4\right) \tau.
\end{align*}
\EP

\begin{remark}
In the above two propositions, the assumption $N>12$ guarantees the integrations in two propositions are finite.
\end{remark}

\BL\label{l9.3}
Assume that $N\geq 52$. Then there exists a $\tau_N\in\mathbb R$ such that $I'(1)=0$, $I''(1)>0$, $I(1)<0$, $J_1(1)>0$ and $J_2(1)>0$.
\EL

\begin{proof}
$I'(1)=0$ is equivalent to
\begin{multline*}
\gamma_N + \beta_N \tau + \alpha_N \tau^2 := \\
   -\frac{8}{(N-12) (N-10) }\left(N^5-4 N^4-80 N^3+208 N^2-32 N-192\right) \\
 -\frac{12   (N-2) }{N-10 } \left(N^3-8 N^2+16\right) \tau\\
 -\frac{4   (N-2)^2  }{ N+4}\left(N^2-4 N-4\right) \tau^2 =0.
\end{multline*}
The assumption $N\geq 52$ guarantees that there exists a solution $\tau_N\in\mathbb R$. Indeed, the discriminant is (after simplified) that
\begin{multline*}
\Delta_N : = \\
\frac{16 (N-2)^2 }{(N-10)^2 (N-12)(N+4)} \Big(N^8-72 N^7+1200 N^6-7200 N^5+15616 N^4 \\
-31744 N^3+34048 N^2+47104 N-49152 \Big).
\end{multline*}
Using Mathematica, we may solve the algebraic equation
\begin{multline}\label{z.2}
\mathcal A(N):=N^8-72 N^7+1200 N^6-7200 N^5+15616 N^4-31744 N^3 \\
 +34048 N^2+47104 N-49152=0
\end{multline}
and find that the largest real solution is about $N\approx 51.1957$.
Thus a real $\tau_N$ satisfying $I'(1)=0$ must exist for $N\geq 52$. We choose
\BE\label{z.1}
\tau_N = \frac{-\beta_N+\sqrt{\Delta_N}}{2\alpha_N} <-\F{99}{50},
\EE
since
\begin{align*}
&\ \gamma_N + \beta_N (-\F{99}{50}) + \alpha_N(-\F{99}{50})^2 \\
=&\ \frac{49 N^6+26730 N^5-1865112 N^4+18883712 N^3-36441904 N^2-9453152 N+45467520}{625 (N-12) (N-10) (N+4)} \\
>&\ 0
\end{align*}
for $N\geq 52$. In fact, we can solve the algebraic equation
\begin{multline*}
49 N^6+26730 N^5-1865112 N^4+18883712 N^3 \\
  -36441904 N^2-9453152 N+45467520 = 0
\end{multline*}
and get the largest $N\approx 51.7253$, while $\alpha_N<0$ for $N\geq 52$.\par

For this $\tau_N$ and $N\geq 52$, we have
\begin{align*}
&\ I''(1)|_{\tau_N}=I''(1)-3I'(1)|_{\tau_N}\\
=&\ -\frac{24 (N-2) \left(N^3-8 N^2+16\right) \tau_N}{N-10} \\
&\qquad -\frac{32 \left(N^5-4 N^4-80 N^3+208 N^2-32 N-192\right)}{(N-12) (N-10)} \\
>&\ \frac{24 (N-2) \left(N^3-8 N^2+16\right)}{N-10}\F{99}{50} \\
&\qquad -\frac{32 \left(N^5-4 N^4-80 N^3+208 N^2-32 N-192\right)}{(N-12) (N-10)} \\
:=&\ \mathcal A_{1N} > 0,
\end{align*}
since the largest $N$ satisfying $\mathcal A_{1N}=0$ is $N\approx 47.248$.\par

Similarly, we can check that, for $N\geq 52$,
\begin{align*}
  J_1(1)|_{\tau_N} >&\ -\frac{6}{(N-10) } \left(N^5-4 N^4-44 N^3+112 N^2-32N-96\right) \\
&\qquad -4 (N-2)^2\left(N^2-4 N-4\right) (-\F{99}{50})\\
>&\ 0.
\end{align*}
Also we have that
\begin{align*}
\ J_2(1)|_{\tau_N}>&\ -\frac{ (N-2) \left(N^4-N^3-56 N^2+40 N+88\right)}{N-10} \\
&\qquad - (N-2)^2 \left(N^2-4 N-4\right) (-\F{99}{50}) \\
>&\ 0.
\end{align*}\par

Finally, we compute the discriminant of $I(1)$ and get
\begin{align*}
  -\frac{24 (N-2)^2 \left(N^7-22 N^6+156 N^5-400 N^4+672 N^3-448 N^2-1024 N+768\right)}{(N-12) (N-10)^2 (N+4)},
\end{align*}
which is checked always negative for $N\geq 52$. So $I(1)<0$.\par

The proof is complete.
\end{proof}

\BP\label{p8.3n}
For $\tau_N$ chosen in Lemma \ref{l9.3}, the function $F(\xi',\lambda')$ has a strict local minimum at $(0,1)$.
\EP

\begin{proof}
Since  $I'(1)=0$, we have $\F{\partial}{\partial\lambda'}F(0,1)=0$. In addition Proposition \ref{p8.1} shows
$\F{\partial}{\partial\xi'}F(0,1)=0$. Therefore, $(0,1)$ is a critical point of $F(\xi',\lambda')$.\par

Since $J_1(1)>0$ and $J_2(1)>0$, it follows from Lemma \ref{l9.3}
that the matrix $(\F{\partial^2}{\partial\xi'_p\partial \xi'_q}F(0,1))$
is positive definite. Lemma \ref{l9.3} again shows that  $I''(1)>0$, which implies $\F{\partial^2}{\partial\lambda'^2}F(0,1)>0$.
Consequently, $(0,1)$ is a strict local minimum point.
\end{proof}

\section{Fourth polynomials and the case of $ 25 \leq N \leq 52$}

 Our  ultimate goal is to reduce the dimension assumption $N\geq 52$ to $N\geq 25$.  Unlike \cite{BM}, where a cubic polynomial
  is chosen, we have to select a 4th order polynomial
  \begin{equation}
  f(s)=\tau -12000s +2411 s^2-135s^3+ s^4.
  \end{equation}

\begin{remark}
The coefficients in $f(s)$ are not unique. And they are chosen in order to verify the conditions ($F_1$)-($F_3$). We have also tried cubic and fifth polynomials. They give larger bounds on $N$.
\end{remark}

Using the software \textit{Mathematica}, we get
\BP\label{p9.23}
Assume $N\geq 25$,
\begin{align*}
&\ F(0,\lambda')\\
=&\ \F{N-4}{16( N^2-4)}|S^{N-1}|\sum_{i,j,k,\ell}(W_{ikj\ell}+W_{i\ell jk})^2\F{\Gamma\left[\F{N}{2}-9\right] \Gamma\left[\F{N}{2}+7\right]}{ \Gamma[N+1]}I(\lambda')
\end{align*}
where $\Gamma$ denotes the usual $\Gamma-\text{function}$ and
\begin{align*}
&\ I(\lambda')\\
=&\ -\frac{(N-4) (N+14) \lambda'^{20}}{(N-24) (N-22) (N-20)}\Big(N^6+42 N^5-768 N^4-17248 N^3+38768 N^2\\
&\hspace{24em}-2336 N-38400\Big) \\
&\ +\frac{270 (N-4) \lambda'^{18}}{(N-22) (N-20)} \Big(N^6+32 N^5-612 N^4-10768N^3+24672 N^2-640 N\\
&\hspace{28em}-25600\Big) \\
&\ -\F{(N-4)\lambda'^{16}}{(N-20) (N+12)}\Big( 23047 N^6+543484 N^5-10985408 N^4-146678256 N^3\\
&\hspace{13.5em} +351063488 N^2-16180224 N-363260160\Big)\\
&\ +\F{30 (N-4)\lambda'^{14}}{(N+10) (N+12)}\Big(22499 N^6+356784 N^5-7984044 N^4-76228592 N^3\\
&\hspace{14em} +193344928 N^2-4902016 N-209193984\Big) \\
&\ -\frac{(N-18) (N-4)\lambda'^{12}}{(N+8) (N+10) (N+12)}\Big(2 \tau N^6+9052921 N^6+4 \tau N^5+84049210 N^5 \\
&\qquad -384 \tau N^4-2265707776 N^4-1856 \tau N^3-14204127072 N^3+4704 \tau N^2\\
&\qquad +40189627120 N^2+7360 \tau N-1125216800 N-12800 \tau-45033619968\Big)\\
&\ +\frac{30 (N-18) (N-16) (N-4) \lambda'^{10}}{(N+8) (N+10) (N+12)}\Big(9 \tau N^5+1928800 N^5-54 \tau N^4\\
&\qquad -3857600 N^4-864 \tau N^3-293177600 N^3+1584 \tau N^2 +648076800 N^2\\
&\hspace{12em} +2880 \tau N+61721600 N-4608 \tau-740659200\Big) \\
&\ -\frac{2 (N-18) (N-16)(N-14) (N-4) \lambda'^8}{(N+6) (N+8) (N+10) (N+12)}\Big(2411 \tau N^5+72000000 N^5\\
&\qquad -14466 \tau N^4-288000000 N^4-135016 \tau N^3-5760000000 N^3\\
&\qquad +308608 \tau N^2+14976000000 N^2+347184 \tau N-2304000000 N\\
&\hspace{21em}-694368   \tau-13824000000\Big) \\
&\ +\frac{2400 (N-18) (N-16) (N-14) (N-12) (N-4) (N-2)  \tau \lambda'^6}{(N+6) (N+8) (N+10)(N+12)}\Big(N^3-8 N^2\\
&\hspace{30em}+16\Big)\\
&\ -\frac{(N-18) (N-16) (N-14) (N-12) (N-10) (N-4) (N-2)^2  \tau^2 \lambda'^4}{(N+4) (N+6) (N+8) (N+10) (N+12)}\Big(N^2\\
&\hspace{27em}-4 N-4\Big)\bigg\}.
\end{align*}
\EP

Also by \textit{Mathematica}, the following holds.
\BP\label{p8.2}
Assume $N\geq 25$,
\begin{align*}
&\ \F{\partial^2}{\partial\xi'_p\partial\xi'_q}F(0,\lambda')\\
=&\ \F{2(N-4)^2}{(N+2)(N+4)}|S^{N-1}|\frac{\Gamma \left(\frac{N}{2}-7\right) \Gamma \left(\frac{N}{2}+5\right)}{ \Gamma (N+1)}\\
&\ \cdot\Bigg\{\sum_{i,j,k}(W_{ikjp}+W_{ipjk})(W_{ikjq}+W_{iqjk})J_1(\lambda')\\
&\qquad\qquad\qquad\qquad +\sum_{i,j,k,\ell}(W_{ikj\ell}+W_{i\ell jk})^2\delta_{pq}J_2(\lambda')\Bigg\}
\end{align*}
where
\begin{align*}
&\ J_1(\lambda')\\
=&\ -\frac{6 \lambda'^{18}  (N+10) (N+12) (N+14) (N+16) }{(N-22) (N-20) (N-18) (N-16) (N-2) }\big(N^5-4 N^4-380 N^3 \\
  &\qquad  +784 N^2-32 N-768\big) \\
&\ +\frac{135 \lambda'^{16} (N+10) (N+12) (N+14)}{2(N-20) (N-18) (N-16) (N-2) } (19 N^5-76 N^4-5776 N^3+11872 N^2\\
&\  -16 N-12160)\\
&\ -\frac{\lambda'^{14}  (N+10) (N+12) }{4 (N-18) (N-16) (N-2) }(340883 N^5-1363532 N^4-80448388 N^3 \\
&\ +167226320 N^2-2575840   N-170104608)\\
&\ +\frac{135 e^{12}  (N+10) }{2 (N-16) (N-2) }(27321 N^5-109284 N^4-4808496 N^3+10054304 N^2\\
&\ +153424N-10644864)\\
&\ -\frac{2 \lambda'^{10}  }{(N-2) } (\tau N^5+8647921 N^5-4 \tau N^4-34591684 N^4-124 \tau N^3\\
&\ -1072342204 N^3+208 \tau N^2+2287434512 N^2+416 \tau N+60226528 N \\
&\ -640 \tau-2486027776)\\
&  +\frac{15 \lambda'^8 (N-14)  }{2 (N-2) (N+8)} (27 \tau N^5+9644000 N^5-108 \tau N^4-38576000 N^4\\
&\ -2160 \tau N^3-771520000N^3+4320 \tau N^2+1728204800 N^2 \\
&\ +5616 \tau N-30860800 N-10368 \tau-1851648000)\\
& -\frac{\lambda'^6 (N-14) (N-12)  }{(N-2)  (N+6) (N+8)} (2411 \tau N^5+108000000 N^5-9644 \tau N^4\\
&\ -432000000 N^4-106084 \tau N^3-4752000000 N^3+270032 \tau N^2 \\
&\ +12096000000 N^2+154304 \tau N-3456000000 N-462912\tau-10368000000)\\
&\ +\frac{6000 \lambda'^4 (N-14) (N-12) (N-10) (N-2)   \tau}{ (N+6) (N+8)}\left(N^2-4 N-4\right)
\end{align*}
and
\begin{align*}
&\ J_2(\lambda')\\
=&\ -\frac{\lambda'^{18} (N+10) (N+12) (N+14) }{2(N-22) (N-20) (N-18) (N-16) (N-2) }(N^6+40 N^5-612 N^4\\
&\ -14416 N^3+33088 N^2-1344 N-33536)\\
&\ +\frac{135 \lambda'^{16} (N+10) (N+12) }{8 (N-20) (N-18) (N-16) (N-2) }(7 N^6+214 N^5-3384 N^4-61888 N^3\\
&\ +145840 N^2-864N-154880)\\
&\ -\frac{\lambda'^{14}  (N+10) }{8 (N-18) (N-16) (N-2) }(69141 N^6+1566133 N^5-25714074 N^4 \\
&\ -353077328 N^3+877267784 N^2-17765904   N-935422272)\\
&\ +\frac{15 \lambda'^{12}  }{8 (N-16) (N-2) }(112495 N^6+1742122 N^5-30701832 N^4\\
&\ -297250720 N^3+795284368 N^2+714144N-887282432)\\
&\ -\frac{\lambda'^{10}  }{4 (N-2)(N+8)}(2 \tau N^6+9052921 N^6+8 \tau N^5+85669210 N^5\\
&\ -312 \tau N^4-1709000832 N^4-1568 \tau N^3-10658431936 N^3+4160 \tau N^2 \\
&\ +32630916432 N^2+5376 \tau N-7615136 N-10240 \tau-37788230528)\\
&\ +\frac{15 \lambda'^8 (N-14) }{8 (N-2)  (N+8)}(27 \tau N^5+5786400 N^5-108 \tau N^4-7715200 N^4-2160 \tau N^3 \\
&\ -648076800 N^3+4320 \tau N^2+1543040000 N^2+5616 \tau N \\
&\ +216025600 N-10368 \tau-1851648000)\\
&\ -\frac{\lambda'^6 (N-14) (N-12)  }{4  (N+6) (N+8)}(2411 \tau N^4+72000000 N^4-4822 \tau N^3\\
&\ -72000000 N^3-115728 \tau N^2-4032000000N^2+38576 \tau N\\
&\ +2880000000 N+231456 \tau +6336000000)\\
&\ +\frac{1500 \lambda'^4 (N-14) (N-12) (N-10) (N-2)  \tau}{ (N+6) (N+8)} \left(N^2-4 N-4\right)\\
   &\qquad \left(N^5+7 N^4-48 N^3+16 N^2+384 N-864\right).
\end{align*}
\EP

\begin{remark}
  Similarly, the assumption $N\geq 25$ is necessary in Proposition \ref{p9.23} and Proposition \ref{p8.2}, otherwise the integrations diverge.
\end{remark}

\BL\label{l8.1}
Assume that $N\geq 25$. Then there exists a $\tau_N\in\mathbb R$ such that $I'(1)=0$, $I''(1)>0$, $I(1)<0$, $J_1(1)>0$ and $J_2(1)>0$.
\EL

\begin{proof}
By the software \textit{Mathematica}, $I'(1)=0$ is equivalent to the following quadratic polynomial
\begin{align*}
  &\ \mathcal B_{1N} \tau^2 + \mathcal B_{2N} \tau + \mathcal B_{3N} \\
= &\ \big(N^{12}-144N^{11}+9112 N^{10}-332800 N^9+7744896 N^8-119529600 N^7 \\
  &\ +1232785280 N^6-8335672320 N^5+35049651968N^4-81419810816 N^3\\
  &\ +72565198848 N^2+41687285760 N-81749606400\big) \tau^2 \\
  &\ +\big(-27025 N^{12}+3697664N^{11}-218286488 N^{10}+7232839376 N^9 \\
  &\ -145917536656 N^8+1798916204864 N^7-12415144280448 N^6 \\
  &\ +29852476048896 N^5+160699860860416 N^4-1111179746377728 N^3\\
  &\ +1284360017780736 N^2+1838260974256128N-2852825357352960\big)\tau\\
  &\ +168227346 N^{12}-21204491480 N^{11}+1108805978944 N^{10}\\
  &\ -30273127750912 N^9+427225222251424N^8-1914627449132672 N^7 \\
  &\ -25810516485700352 N^6+354281472809361408 N^5-800700785348505600N^4\\
  &\ -5829654487294640128 N^3+15038003695249195008 N^2 \\
  &\ -1338454157952024576 N-14134493112544788480 = 0.
\end{align*}
After simplified, the discriminant is
\begin{align*}
  &(N-24) (N-22) (N-20) (N-18) (N-2)^2 (N+4) \big(57441241 N^{17} \\
  &-13316757740 N^{16}+1421208951488 N^{15} \\
  &-92295209252880 N^{14}+4060887517487792 N^{13}\\
  &-127528377218205952 N^{12}+2932837691854966528N^{11}\\
  &-49861018904126426112 N^{10}+624738394629537111040 N^9 \\
  &-5683787728574744100864N^8+36513302074683044208640 N^7\\
  &-158739757047539234324480 N^6+443941679903779546513408N^5\\
  &-788934839032708877123584 N^4+947159822427449128648704 N^3 \\
  &-70877049300753252876288N^2-1727615795557515443306496 N\\
  &+1156307714218965199749120\big)=0.
\end{align*}
By Mathematica, we may check that the biggest zero of the last term is $N\approx 24.9422$.
Since $N\geq 25$, the discriminant is positive and there exists a real $\tau_N>17000$ such that $I'(1)=0$.
In fact, the largest $N$ satisfying $17000\mathcal B_{1N}+17000\mathcal B_{2N}+\mathcal B_{3N}=0$ is about $N\approx 24.9982$
and $17000\mathcal B_{1N}+17000\mathcal B_{2N}+\mathcal B_{3N}\to-\infty$ as $N\to+\infty$.
Thus
\BE\label{z.3}
17000\mathcal B_{1N}+17000\mathcal B_{2N}+\mathcal B_{3N}<0 \qquad \text{when }N\geq 25.
\EE
On the other hand, the largest one to $\mathcal B_{1N}=0$ is $N=24$.
So
\BE\label{z.4}
 \mathcal B_{1N}>0\qquad \text{for }N\geq 25.
\EE
Therefore a $\tau_N>17000$ does exist from (\ref{z.3}) and (\ref{z.4}). \par

As for $I''(1)$, we have
\begin{align*}
 &\  I''(1)\Big|_{\tau_N} =  I''(1) - 3I'(1)\Big|_{\tau_N} \\
 =&\ \frac{(N-4) }{2 (N-24) (N-22) (N-20) (N+6) (N+8) (N+10) (N+12)} \\
 &\ \Big[\big(18713 N^{11}-2782116 N^{10}+180731080 N^9-6710115152 N^8 \\
 &\ +156124054384 N^7-2341982210432 N^6+22336952299904 N^5\\
 &\ -126901060502528N^4+358317913417728 N^3-193547489329152 N^2\\
 &\ -744698653507584 N+844209499668480\big)\tau_N \\
 &\   -239387740 N^{11}+33343448176 N^{10}-1965147182144 N^9\\
 &\ +62606999818560 N^8-1119986901305472 N^7+9710408237796864 N^6\\
 &\ + 2497983280251904 N^5-775332925643090944 N^4+5322471076556407808N^3\\
 &\ -8454330859402850304 N^2-1265620067537485824 N+8423950663121633280 \Big]\\
 :=&\ \frac{(N-4)\left[\mathcal B_{4N}\tau_N+\mathcal B_{5N}\right] }{2 (N-24) (N-22) (N-20) (N+6) (N+8) (N+10) (N+12)}.
\end{align*}
On account that $\mathcal B_{4N}>0$ and $17000\mathcal B_{4N}+\mathcal B_{5N}>0$ for $N\geq 25$, we know that $I''(1)\Big|_{\tau_N}>0$.\par

Come to $J_1(1)$. Direct computation shows that
\begin{align*}
  &\ \F{4 (N-22) (N-20) (N-18) (N-16) (N-2) (N+4) (N+6) (N+8)}{N}J_1(1)|_{\tau_N} \\
 =&\ (15158 N^{11}-1928488 N^{10}+106956432N^9-3385235008 N^8 \\
  &\ +67144207872 N^7-861032971776 N^6+7082375502592 N^5 \\
  &\ -35544691776512N^4+95816075389440 N^3-91338406885376 N^2\\
  &\ -68140868911104 N+124511268372480)\tau_N\\
   &\ -204822475 N^{11}+24606405950 N^{10}-1253228747592 N^9 \\
   &\ +34711856215872 N^8-549030905462832N^7+4489294448546432 N^6\\
   &\ -7165699950704832 N^5-169724096918379392 N^4+1156114211550506752N^3\\
   &\ -1776295972609390592 N^2-81756788717899776 N+1606525147689615360 \\
 :=&\ \mathcal B_{6N}\tau_N + \mathcal B_{7N}.
\end{align*}
In respect that $\mathcal B_{6N}>0$ and $17000\mathcal B_{6N} + \mathcal B_{7N}>0$ for $N\geq 25$, we obtain that $J_1(1)|_{\tau_N}>0$.\par

As for $J_2(1)$, we can check by Mathematica that
\begin{align*}
  &\ \F{8 (N-22) (N-20) (N-18) (N-16) (N-2) (N+6) (N+8)}{N}J_2(1)|_{\tau_N} \\
 =&\ (7579 N^{11}-964244 N^{10}+53478216 N^9-1692617504 N^8+33572103936 N^7 \\
 &\ -430516485888 N^6+3541187751296 N^5-17772345888256 N^4\\
 &\ +47908037694720 N^3-45669203442688 N^2-34070434455552 N \\
 &\ +62255634186240)\tau_N\\
   &\ -73690617 N^{11}+8757708965 N^{10}-439144014720 N^9\\
   &\ +11871532517728 N^8-179607080990760 N^7+1303031622196560 N^6\\
   &\ +730732776846400 N^5-81834334806699648 N^4+487269552686137472N^3\\
   &\ -729001679303608320 N^2-294593226517125120 N+897917009560289280 \\
 :=&\ \mathcal B_{8N}\tau_N + \mathcal B_{9N}.
\end{align*}
Since $\mathcal B_{8N}>0$ and $17000\mathcal B_{8N} + \mathcal B_{9N}>0$ when $N\geq 25$, it holds that $J_2(1)|_{\tau_N}>0$. \par

Finally, direct calculation gives that
\begin{align*}
  &\ -\F{(N-24) (N-22) (N-20) (N+4) (N+6) (N+8) (N+10) (N+12)}{N-4}I(1) \\
 =&\ (N^{12}-144 N^{11}+9112N^{10}-332800 N^9+7744896 N^8-119529600 N^7 \\
 &\ +1232785280 N^6-8335672320 N^5+35049651968 N^4-81419810816N^3 \\
 &\ +72565198848 N^2+41687285760 N-81749606400) \tau^2 \\
 &\ +(-19446 N^{12}+2634424 N^{11}-153791984 N^{10}+5031783008 N^9 \\
 &\ -100052361056 N^8+1212544446208 N^7-8182468308224 N^6 \\
 &\ +18630407205888 N^5+108938164472832 N^4-719625359429632 N^3 \\
 &\ +817980109406208 N^2+1176909695483904 N-1816259033825280)\tau\\
   &\ +94536729 N^{12}-11656568130 N^{11}+593919303092 N^{10}-15702815300040 N^9 \\
   &\ +211574113757984 N^8-815960309177792 N^7-13866450151972480 N^6 \\
   &\ +171265848827756032 N^5-341975328615126528 N^4 \\
   &\ -2799428862920112128 N^3+7040619868808921088 N^2\\
   &\ -607654647921180672 N-6549150991381954560\\
 :=&\ \mathcal B_{10}\tau^2 + \mathcal B_{11}\tau + \mathcal B_{12}.
\end{align*}
Because $\mathcal B_{11}^2-4\mathcal B_{10}\mathcal B_{12}<0$ and $\mathcal B_{10}>0$ for $N\geq 25$, $I(1)$ then must be negative.\par

The proof is complete.
\end{proof}

Similar to the proof of Proposition \ref{p8.3n}, we obtain the following

\BP\label{p8.3}
Let $N \geq 25$. For $\tau_N$ chosen in Lemma \ref{l8.1}, the function $F(\xi',\lambda')$ has a strict local minimum at $(0,1)$.
\EP

\section{Proof of the main theorem}

In this section we prove the main result of this paper by a gluing method.

\BP\label{p9.1}
Assume $N\geq 25$. Moreover, let $g$ be a smooth metric on $\mathbb R^N$ of the form $g(x)=e^{h(x)}$,
where $h(x)$ is a trace-free symmetric two-tensor on $\mathbb R^N$ such that
\BEN
|h(x)|+|\partial h(x)|+|\partial^2 h(x)|+|\partial^3 h(x)|+|\partial^4 h(x)|\leq\alpha
\EEN
for all $x\in\mathbb R^N$, $h(0)=0$, $h(x)=0$ for $|x|\geq 1$, and
\BEN
h_{ik}(x)=\mu\varepsilon^8 f(\varepsilon^{-2}|x|^2) H_{ik}(x)
\EEN
for $|x|\leq \rho$. If $\alpha$ and $\rho^{4-n}\mu^{-2}\varepsilon^{N-24}$ are sufficiently small, then there exists
a positive solution $u(x)$ to
\BEN
P_g u &=\F{N-4}{2}u^\F{N+4}{N-4}  \qquad \text{in }\mathbb R^N,\\
\int_{\mathbb R^N} u^\F{2N}{N-4} &< \int_{\mathbb R^N} \left(\F{1}{1+|y|^2}\right)^N\\
\sup_{|x| \leq \varepsilon} u &\geq C\varepsilon^\F{4-N}{2}.
\EEN
\EP

\begin{proof}
By Proposition \ref{p8.3}, the function $F(\xi',\lambda')$ has a strict local minimum at $(0,1)$ and $F(0,1)<0$. Hence,
we can find an open set $\mathcal M\subset\Lambda$ such that $(0,1)\in\mathcal M$ and
\BEN
F(0,1) < \inf_{\partial\mathcal M} F(\xi',\lambda')<0.
\EEN
Using Lemma \ref{l5.2}, we obtain
\BEN
2\mathcal F_{\tilde g}(\xi',\lambda')=2E+\mu^2\varepsilon^{20}F(\xi',\lambda')+O\left(\mu^3\varepsilon^\F{20N}{N-1}+\mu^\F{20N}{N-4}\varepsilon^\F{20N}{N-4}\right)+O(\alpha(\F{\varepsilon}{\rho})^{N-4}).
\EEN
Hence, if $\rho^{4-n}\mu^{-2}\varepsilon^{N-24}$ is sufficiently small, we have
\begin{align*}
\mathcal F_{\tilde g}(0,1) < \inf_{\partial\mathcal M} \mathcal F_{\tilde g}(\xi',\lambda')<E.
\end{align*}
Consequently, there exists a point $(\bar\xi', \bar\lambda')\in \mathcal M$ such that
\begin{align*}
\mathcal F_{\tilde g}(\bar\xi', \bar\lambda') < \inf_{\partial\mathcal M} \mathcal F_{\tilde g}(\xi',\lambda')<E.
\end{align*}
Then $\tilde u_0(y)+\phi(y)$ is a  solution of (\ref{b}) with $\|\phi\|_*\leq C$.

Note that since $\| \phi \|_{*} \leq C$ we have
\begin{equation}
|\phi (y)| \leq C \frac{\alpha}{ (|+|y-\xi^{'}|)^{N-4}} \leq C \alpha u_0
\end{equation}
which shows that  $ u_0 +\phi >0$ provided $\alpha$ is small. Thus $u_0(x)+\varepsilon^{-\F{N-4}{2}}\phi(x/\varepsilon)$ is the positive solution we need.
\end{proof}

\BP
Let $N\geq 25$. Then there exists a smooth metric $g$ on $\mathbb R^N$ with the following properties:
\begin{itemize}
\item[1)] $g_{ij}(x)=\delta_{ij}$ for $|x|\geq \F{1}{2}$,
\item[2)] $g$ is not conformally flat,
\item[3)] There exists a sequence of positive function $u_n$ $(n\in\mathbb R^N)$ such that
\begin{align*}
P_g u_n &= \F{N-4}{2} u_n^\F{N+4}{N-4},\\
\int_{\mathbb R^N} u_n^\F{2N}{N-4} &< \int_{\mathbb R^N}\left(\F{1}{1+|y|^2}\right)^N,\\
\sup_{|x| \leq 1} u_n &\to \infty.
\end{align*}
\end{itemize}
\EP

\begin{proof}
Choose a smooth cut-off function $\eta$ such that $\eta(r)=1$ for $r\leq 1$ and $\eta(r)=0$ for $r\geq 2$.
We define a trace-free symmetric two-tensor on $\mathbb R^N$ by
\begin{align*}
h_{ij}(x)=\sum_{n=N_0}^\infty \eta(4n^2|x-x_n|) 2^{-\F{25}{3}n} f(2^{2n}|x-x_n|) H_{ij}(x-x_n),
\end{align*}
where $x_n=(\F{1}{n},0,\ldots,0)$. Clearly $h(x)$ is $C^\infty$.\par

Moreover, if $N_0$ is sufficiently large, then we have $h(x)=0$ for $|x|\geq \F{1}{2}$ and $|h|+|\partial h|+|\partial^2 h|+|\partial^3 h|+|\partial^4 h|\leq \alpha$.
Provided that $n\geq N_0$ and $|x-x_n|\leq \F{1}{4n^2}$, we have
\begin{align*}
h_{ij}(x)= 2^{-\F{25}{3}n} f(2^{2n}|x-x_n|) H_{ij}(x-x_n).
\end{align*}
Hence we can apply Proposition \ref{p9.1} with $\mu=2^{-n/3}$, $\varepsilon=2^{-n}$, $\rho=\F{1}{4n^2}$. From this the assertion follows.
\end{proof}

\section{Appendix}

In this section we will give the proof of (\ref{2.11}) and (\ref{2.12}). The proof of (\ref{2.11}) can be found in \cite{LN} and
we repeat it here for the sake of convience.

\BL
Assume that $0<s<N$ and $t>s$. Then
\begin{align*}
\int_{\mathbb R^N} \F{1}{|x-y|^{N-s}}\F{1}{(1+|y|)^t}\mathrm{d}y \leq
\begin{cases}
C (1+|x|)^{s-t} & \text{if }t<N,\\
C (1+|x|)^{s-N} \Big[1+\log(1+|x|)\Big] \quad & \text{if }t=N,\\
C (1+|x|)^{s-N} & \text{if }t>N.
\end{cases}
\end{align*}
\EL

\begin{proof}
First, observe that the above integral is well defined since $t>s$. So we only need to consider the case that $|x|$ is large.
Next we decompose it as follows:
\begin{align*}
&\ \left(\int_{|y-x|\leq \F{|x|}{2}}+\int_{\F{|x|}{2}\leq |y-x|\leq 2|x|}+\int_{|y-x|\geq 2|x|}\right)\F{1}{|x-y|^{N-s}}\F{1}{(1+|y|)^t}\mathrm{d}y\\
:=&\ I_1+I_2+I_3.
\end{align*}\par

$I_1$ may be estimated as follows. Since $|y-x|\leq |x|/2$ implies $|y|\geq |x|/2$,
\begin{align*}
I_1 &\leq \int_{|y-x|\leq \F{|x|}{2}}\F{1}{|x-y|^{N-s}}\F{1}{(1+|x|/2)^t}\mathrm{d}y\\
&\leq \F{1}{(1+|x|/2)^t} \int_0^{|x|/2} \F{1}{r^{N-s}}r^{N-1}\mathrm{d}r \\
&\leq C|x|^{s-t}.
\end{align*}
$I_3$ may be estimated similarly. Because $|y-x|\geq 2|x|$, $|y-x|\leq |y|+|x|\leq |y|+|y-x|/2$. Thus $|y-x|\leq 2|y|$ and
\begin{align*}
I_3 &\leq \int_{|y-x|\geq 2|x|} \F{1}{|x-y|^{N-s}} \F{1}{(1+|x-y|/2)^t}\mathrm{d}y\\
&\leq \int_{2|x|}^\infty \F{1}{r^{N-s}} \F{1}{(1+r/2)^t} r^{N-1} \mathrm{d}r\\
&\ \leq C|x|^{s-t}.
\end{align*}
Finally, we observe that
\begin{align*}
I_2 &\leq \F{C}{|x|^{N-s}}\int_{\F{|x|}{2}\leq |y-x|\leq 2|x|} \F{1}{(1+|y|)^t}\mathrm{d} y\\
&\leq \F{C}{|x|^{N-s}} \left(\int_{|y|\leq 1}+\int_{1\leq |y|\leq 3|x|} \right)\F{1}{(1+|y|)^t}\mathrm{d} y\\
&\leq C \F{C}{|x|^{N-s}} \left(C+C\int_1^{3|x|} r^{N-t-1} \mathrm{d}r\right)\\
&\leq
\begin{cases}
C |x|^{s-t} & \text{if }t<N,\\
C |x|^{s-N} \log|x| \qquad & \text{if }t=N,\\
C |x|^{s-N} & \text{if }t>N.
\end{cases}
\end{align*}
Now it is easily seen that the lemma holds.
\end{proof}

Next we come to the proof of (\ref{2.12}).
\begin{proof}[Proof of (\ref{2.12})]
Now $t=N-k$. Let $L>0$ is a large number. If $|y|\leq L r$, we have
\begin{align*}
&\ \int_{B_r} \F{1}{|y-z|^{N-s}}\F{1}{(1+|z|)^{N-k}}\mathrm{d}z \\
\leq&\ \int_{\mathbb R^N} \F{1}{|y-z|^{N-s}}\F{1}{(1+|z|)^{N-k}}\mathrm{d}z \\
\leq &\ C(1+|y|)^{s+k-N} \leq C r^{k} (1+|y|)^{s-N}.
\end{align*}
For $|y|\geq Lr$, obvious $\F{(1+|y|)^{N-s}}{|y-z|^{N-s}} \leq C$ since $|z|\leq r$.
Thus, recalling that $k>0$, we get
\begin{align*}
&\ \int_{B_r} \F{1}{|y-z|^{N-s}}\F{1}{(1+|z|)^{N-k}}\mathrm{d}z \\
\leq&\ C(1+|y|)^{s-N} \int_{B_r} \F{1}{(1+|z|)^{N-k}}\mathrm{d}z\\
\leq &\ Cr^{k}(1+|y|)^{s-N}.
\end{align*}
This concludes the proof.
\end{proof}

\end{document}